\documentclass[12pt, reqno]{article}
\pdfoutput=1
\usepackage{amssymb,amsmath,amsthm,amsfonts,latexsym}
\usepackage{epsfig,esint}
\usepackage{mathrsfs}
\usepackage{mathtools,textcomp}
\usepackage{graphicx}
\usepackage{epic,eepic,wrapfig,color, ifthen} 
\usepackage{url,cite}
\usepackage[colorlinks=true, citecolor=blue, filecolor=cyan, linkcolor=red, urlcolor=magenta]{hyperref}
\usepackage{pmboxdraw}
\usepackage{verbatim}
\usepackage[normalem]{ulem}
\usepackage{enumerate}
\usepackage[bottom]{footmisc}
\usepackage{lipsum}

\let\OLDthebibliography\thebibliography
\renewcommand\thebibliography[1]{
  \OLDthebibliography{#1}
  \setlength{\parskip}{0pt}
  \setlength{\itemsep}{0pt plus 0.3ex}
}

\numberwithin{equation}{section}

\mathtoolsset{showonlyrefs}

\newcommand{\bk}{\color{black}}

\def\pf{\noindent{\bf Proof. }}
\def\qed{{\hfill $\Box$ \bigskip}}

\def\1{{\bf 1}}
\def\nn{\nonumber}

 \def\sK {{\cal K}}

\def\R {{\mathbb R}} 
\def\N{{\mathbb N}} 
\def\E {{\mathbb E}}  \def \P{{\mathbb P}}

\def\eps{\varepsilon}

\def\wt{\widetilde}

\def\Kge{\sK^{\rm loc}_\ge(\beta,\kappa)}
\def\Kle{\sK^{\rm loc}_\le(\beta,\kappa)}
\def\Kloc{\sK^{\rm loc}(\beta,\kappa)}
\def\KK{\sK(\beta,\kappa)}

\def\Log{\text{\rm Log}}
\def\Lg{\text{\rm Log}\,}

\theoremstyle{plain}
\newtheorem{thm}{Theorem}[section]
\newtheorem{lem}[thm]{Lemma}
\newtheorem{cor}[thm]{Corollary}
\newtheorem{remark}[thm]{Remark}
\newtheorem{prop}[thm]{Proposition}
\newtheorem{defn}[thm]{Definition}

\theoremstyle{definition}
\newtheorem*{eg*}{Example}
\newtheorem*{egs*}{Examples}
\newtheorem*{def*}{Definition}
\theoremstyle{remark}

\begin{document}
	
	\title{Heat kernel estimates for Schr\"odinger operators with supercritical killing potentials}	
	\date{}

	\author{
		{\bf Soobin Cho}
		\quad 	{\bf Panki Kim\thanks{This work was supported by the National Research Foundation of Korea(NRF) grant funded by the Korea government (MSIT) (No. RS-2023-00270314)}}
		\quad {\bf Renming Song\thanks{Research supported in part by a grant from the Simons Foundation (\#960480, Renming Song).}
		}
	}
	
	\maketitle
	
	\begin{abstract}
		In this paper, we study the Schr\"odinger operator $\Delta-V$, where $V$ is a supercritical non-negative potential belonging to a large class of functions containing functions of the form $b|x|^{-(2+2\beta)}$, $b, \beta>0$.
		We obtain two-sided estimates on the heat kernel $p(t, x, y)$ of $\Delta-V$, along with estimates for the corresponding Green function. Unlike the case of the fractional Schr\"odinger operator $-(-\Delta)^{\alpha/2}-V$, $\alpha\in (0, 2)$, with
		supercritical killing potential dealt with in \cite{CS24}, in the present case, the heat kernel $p(t, x, y)$ decays to 0 exponentially
		as $x$ or $y$ tends to the origin. 
		
		\medskip
		
		\medskip
		
		
	\end{abstract}
	\allowdisplaybreaks

	\section{Introduction}\label{s:intro}	
	
	The Schr\"odinger operator 
	$\Delta-V$
	is of vital importance in mathematical physics 
	and chemistry,
	 analysis and probability, 
	and its associated semigroup is usually called a Schr\"odinger semigroup or a 
	Feynman-Kac semigroup. The Schr\"odinger operator and its Feynman-Kac semigroup have been studied intensively,
	see, for instance, \cite{CZ, Simon} and the references therein for early contributions, and \cite{BHJ08, CW23a, CW23b}
	and the references therein for recent contributions. 
	
	The main concern of this paper is on two-sided heat kernel estimates for Schr\"odinger operators with 
	potentials behaving like 
	inverse powers $|x|^{-\gamma}$ near the origin. 
	Inverse power potentials, such as the  Coulomb potential $V(x)=c|x|^{-1}$,  are very important in mathematical physics 
 and chemistry. 
 Various potentials   $V(x)$ proportional to  central potentials  $|x|^{-\beta}$
  for  $\beta>0$, particularly for even integer $\beta$, are commonly studied, with applications in areas such as electrostatics, quantum systems, molecular interactions, and the stability of matter. For further details,   
 we refer the reader to  \cite{AZ13,PR62, VW54, MRSW, LeBe70, KP92, Ba80, IS09} and the references therein.

	When $\gamma\in (0, 2)$ and $d\ge 2$,  or when $\gamma\in (0, 1)$ and $d=1$, the function $V(x)=\kappa|x|^{-\gamma}$ belongs to the Kato class, that is, 
	$$
	\lim_{t\to0}\sup_{x\in \R^d}\int^t_0\int_{\R^d}(4\pi s)^{d/2}e^{-|x-y|^2/
	(4s)}|V(y)|dyds=0.
	$$
	A potential $V$ belonging to the Kato class is small,  in the sense of quadratic forms,  compared to $\Delta$, 
	so $\Delta-V$ can be regarded as a small (or subcritical) perturbation of $\Delta$. It is well known that, when $V$ belongs to the Kato class, the heat kernel of 
	$\Delta-V$ admits short time two-sided  Gaussian estimates, 
	see, for instance, \cite{BM}. We mention in passing that, when $d=1$ and $\gamma\in [1, 2)$, the function $V(x)=\kappa|x|^{-\gamma}$ does not belong to the Kato class.
	As far as we know, in the case $d=1$, heat kernel
	estimates of $\Delta-\kappa|x|^{-\gamma}$, $\gamma\in [1, 2)$, have not been studied in the literature. However, this
	is not the concern of this paper.
	
	The function $V(x)=\kappa|x|^{-2}$ does not belong to the Kato class. According to the Hardy inequality,
	the quadratic form
	$$
	\int_{\R^d}|\nabla u(x)|^2dx+\kappa\int_{\R^d}u^2(x)|x|^{-2}dx, \quad u\in C^\infty_c(\R^d),
	$$
	is non-negative if and only if $\kappa\ge -(d-2)^2/4$.  In this sense, $V(x)=\kappa|x|^{-2}$ is a critical potential. 
	The study of Schr\"odinger operators $\Delta-V$  with critical potentials $V(x)=\kappa|x|^{-2}$ goes back to
	Baras and Goldstein \cite{Bar-Gol84}, who proved the existence of 
	non-trivial nonnegative solutions of the 
	heat equation 
	$\partial_tu(t, x)=\Delta u(t,x)-\kappa|x|^{-2}u(t, x)$ in $\R^d$ for 
    $\kappa\in [-(d-2)^2/4, 0]$, 
	and non-existence of 
	such solutions for $\kappa<-(d-2)^2/4$. Vazquez and Zuazua \cite{VZ} studied the Cauchy problem
	and spectral properties of the operator in bounded subsets of $\R^d$.  Sharp two-sided estimates for the heat kernel of 
	$\Delta-\kappa|x|^{-2}$ were obtained by Liskevich and
	Sobol \cite[p. 365, Examples 3.8, 4.4 and 4.10]{Lis-Sob03} for the case $\kappa\in (-(d-2)^2/4, 0)$. Milman and Semenov
	established  sharp 	two-sided bounds for the case $\kappa\ge -(d-2)^2/4$, see  \cite[Theorem 1]{MS04} and \cite{MS05}. 
	See also  Barbatis, Filippas and Tertikas \cite {BFT04}, and  Filippas, Moschini and Tertikas \cite{FMT07} for related results.

		Despite their importance, there have been limited results for Schrödinger operators $\Delta - V$ with 
	supercritical potentials $V(x)=\kappa|x|^{-\gamma}$ $(\gamma>2)$ in the literature. In particular,  there have been no results regarding their heat kernel estimates.
	For a certain class of (quasi-)exactly solvable potentials $V$, 
	exact solutions to the associated elliptic problems have been studied 
	in
	Turbiner \cite{Tu16, Tu88},  Dong, Ma and Esposito \cite{DME99}, Dong \cite{Do00, Do01}, and Agboola and Zhang \cite{AZ13},  among others. 
	Li and Zhang \cite{LZ17} studied regularity of weak solutions of elliptic and parabolic equations associated with
	Schr\"odinger operators with supercritical potential. 
	In the supercritical case, to guarantee the corresponding 
	quadratic form to be non-negative definite,
	we have to assume $\kappa>0$.  
	 
	Our main goal is to obtain two-sided
	heat kernel estimates of 	Schr\"odinger operator with a general class of supercritical potentials which includes $\Delta-\kappa|x|^{-(2+2\beta)}$ for 
	$\beta>0$ and $\kappa>0$.  To state our heat kernel estimates for $\Delta-\kappa|x|^{-(2+2\beta)}$, 
	we introduce some notation first. Throughout this paper, we use $\R^d_0$ to denote $\R^d\setminus\{0\}$ and $q(t, x, y)$ to denote the Gaussian heat kernel
	$$
	q(t,x,y) = (4\pi t)^{-d/2} e^{- {|x-y|^2}/{(4t)}}
	$$
	and use the notation
	$$
		\Lg r := \log(e-1+r), \quad r\ge 0.
$$
	For $\beta, \kappa>0$, we define a constant
	$\eta_1=\eta_1(\beta, \kappa):=(2^{-(12+7\beta)/(2+\beta)} \beta \sqrt \kappa)^{1/(2+\beta)}$, 
	 a function $h_{\beta,\kappa}$ on $(0, \infty)$ by
\begin{align}\label{e:function-h}
		h_{\beta,\kappa}(r):=		(r\wedge 1)^{-(d-2-\beta)/2} e^{	 -\frac{\sqrt \kappa}{\beta (r\wedge 1)^\beta}},
\end{align}
	and a function $H_{d, \beta, \kappa}$ on $(0, \infty)\times (0, \infty)$ by
	\begin{align}\label{e:H}
		H_{d,\beta,\kappa}(t,r):= \begin{cases}
			\displaystyle 
			h_{\beta,\kappa}(r) &\mbox{ if $d\ge 3$}, \\[4pt]
			\displaystyle   1 \wedge \frac{h_{\beta,\kappa}(r)\1_{(0,1)}(r) + 			( \Lg r ) 
			\1_{[1, \infty)}(r) }{\Lg \sqrt t} &\mbox{ if $d=2$}, \\[10pt]
			\displaystyle  1 \wedge \frac{h_{\beta,\kappa}(r) \1_{(0,1)}(r) + r \1_{[1, \infty)}(r)}{\sqrt t} &\mbox{ if $d= 1$}.
		\end{cases}
	\end{align}
	Our main results, specialized to the case of 
	$\Delta-\kappa|x|^{-(2+2\beta)}$, $\beta>0$ and $\kappa>0$,
	are as follows.
	
	\begin{thm}\label{t:special-smalltime} Suppose $\beta>0$ and $\kappa>0$. The Schr\"odinger operator
		$\Delta-\kappa|x|^{-(2+2\beta)}$ admits a heat kernel $p^{\beta, \kappa}(t, x, y)$ satisfying the following small time estimates: 
		
		\noindent
		(i) There exist $c_1,c_2,c_3>0$ such that  for all $t \in (0, 		4]$ and $x,y \in \R^d_0$,
		\begin{align*}
			p^{\beta, \kappa}(t,x,y) &\le c_1\bigg( 1 \wedge  \frac{h_{\beta,\kappa}(|x| )}{h_{\beta,\kappa}( \eta_1 t^{1/(2+\beta)})}\bigg)\bigg( 1 \wedge  \frac{h_{\beta,\kappa}(|y|)}{h_{\beta,\kappa}( \eta_1 t^{1/(2+\beta)})}\bigg)  \nn\\
			&\quad \times \exp \bigg(- \frac{c_2t}{(|x| \vee |y| \vee t^{1/(2+\beta)})^{2+2\beta}}\bigg)q(c_3t,x,y).
		\end{align*}
		
			\noindent
		(ii)
		(a) If $d\ge 2$, then	 there exist  $c_4,c_5,c_6>0$ such that  for all $t \in (0,	4]$ and $x,y \in \R^d_0$,
		\begin{align}\label{e:special-smalltime-lower}
			p^{\beta, \kappa}(t,x,y) &\ge c_4\bigg( 1 \wedge  \frac{h_{\beta,\kappa}(|x| )}{h_{\beta,\kappa}( \eta_1 t^{1/(2+\beta)})}\bigg)\bigg( 1 \wedge  \frac{h_{\beta,\kappa}(|y|)}{h_{\beta,\kappa}( \eta_1 t^{1/(2+\beta)})}\bigg)  \nn\\
			&\quad \times \exp \bigg(- \frac{c_5t}{(|x| \vee |y| \vee t^{1/(2+\beta)})^{2+2\beta}}\bigg)q(c_6t,x,y).
		\end{align}
		
		(b) If  $d=1$, then \eqref{e:special-smalltime-lower} holds for all 		 $t \in (0,4]$ 
		and $x,y \in \R^1_0$ with $xy>0$.
	\end{thm} 
	
	\begin{thm}\label{t:special-largetime}
		Suppose $\beta>0$ and $\kappa>0$. The heat kernel $p^{\beta, \kappa}(t,x,y) $ of $\Delta-\kappa|x|^{-(2+2\beta)}$ satisfies
		the following large time estimates: 
		
		\noindent
		(i) If $d\ge 2$, then	 there exist $c_i>0$, $1\le i \le 4$, such that  for all $t \ge 4$ and $x,y \in \R^d_0$,
		\begin{equation}\label{e:special-largetime}
			c_1 q(c_2t,x,y){ H_{d,\beta,\kappa}(t,x) H_{d,\beta,\kappa}(t,y)} \le	{p^{\beta, \kappa}(t,x,y)}\le c_3  q(c_4t,x,y){ H_{d,\beta,\kappa}(t,x) H_{d,\beta,\kappa}(t,y)}.
		\end{equation}
		
		\noindent
		(ii) If  $d=1$, then \eqref{e:special-largetime}  holds for all $t \ge 4$ and $x,y \in \R^1_0$ with $xy>0$.
	\end{thm}
	In fact, in Theorems \ref{t:smalltime} and \ref{t:largetime}, we establish the 
	heat kernel estimates of 	Schr\"odinger operator with a general class of supercritical 
	 potentials including $\kappa|x|^{-(2+2\beta)}$.
	Theorems \ref{t:special-smalltime} and \ref{t:special-largetime} are immediate consequences of 
	Theorems \ref{t:smalltime} and \ref{t:largetime} respectively.

As a direct consequence of our global heat kernel estimates, we  obtain  
two-sided Green function estimates. 
In particular, the  Green function $	G^{\beta,\kappa}(x,y):=\int_0^\infty p^{\beta,\kappa}(t,x,y)dt$ of $\Delta - \kappa|x|^{-(2+2\beta)}$ has the following estimates.
	\begin{thm}\label{t:special-Green}
		Suppose   $d\ge 3$.  Then there exist constants $\eta_2>0$, $c_1\ge c_2>0$ and $C>1$ such that   the following Green function estimates hold:
	\begin{align*}
	&C^{-1}\bigg( 1\wedge \frac{h_{\beta,\kappa}(|x|\wedge |y|)}{h_{\beta,\kappa}(\eta_2|x-y|)}\bigg) |x-y|^{2-d} \exp \bigg( - \frac{c_1|x-y|}{(|x| \vee |y|)^{1+\beta}}\bigg)\\
	&\le 	G^{\beta,\kappa}(x,y)  \le  C\bigg( 1\wedge \frac{h_{\beta,\kappa}(|x|\wedge |y|)}{h_{\beta,\kappa}(\eta_2|x-y|)}\bigg) |x-y|^{2-d} \exp \bigg( - \frac{c_2|x-y|}{(|x| \vee |y|)^{1+\beta}}\bigg).
\end{align*}
	\end{thm}
	Theorem \ref{t:special-Green}, along with the Green function estimates for the cases $d=1$ or $d=2$, follows from Theorem \ref{t:Green} and Remark \ref{r:Green}.

	This paper is a companion paper of \cite{CS24} which established 
	two-sided heat kernel estimates for the
	supercritical fractional Schr\"odinger operators $-(-\Delta)^{\alpha/2}-\kappa|x|^{-\gamma}$, $\alpha\in (0, 2)$, $\gamma>\alpha$
	and $\kappa>0$. However, there are some dramatic differences between the heat kernel estimates for the local and non-local cases. For the supercritical fractional Schr\"odinger operators 
	$-(-\Delta)^{\alpha/2}-\kappa|x|^{-\gamma}$, the heat kernel  
	decays to 0 like $|x|^\gamma$ as $|x|\to 0$,  see \cite[Theorem 1.1]{CS24}. While in the local case, as one see from the theorems above, the heat kernel  $ p^{\beta, \kappa}(t,x,y)$ of of $\Delta-\kappa|x|^{-(2+2\beta)}$ decays to 0 exponentially
	as $|x|\to 0$. 
	
	The approach of this paper is probabilistic. The 
	Feynman-Kac semigroup, and the sharp two-sided Dirichlet heat estimates for the Laplacian in 
	exterior balls
	obtained in \cite{Zh03, GS11} play essential roles. 

	Organization of the paper: 
	In Section \ref{s:2}, we give some preliminaries. In Section \ref{s:3}, we 
define the class of supercritical killing potentials and
	establish  survival probability estimates for Schr\"odinger operators with potentials within this class. These  estimates will play important roles in subsequent sections. Small time heat kernel estimates are derived  in Section \ref{s:4}.   In Section \ref{s:counters}, we show the optimality of our assumptions on the potentials by presenting counterexamples that show the failure of small time heat  kernel estimates when the potential goes beyond our class.
	Large time heat kernel estimates are proved in Section \ref{s:large}, and Green function estimates are proved in Section \ref{s:7}.

\smallskip

\noindent {\bf Notations}: 
Throughout this paper,   lower case letters  $c_i$, $i=0,1,2,...$, which denote positive real constants, are fixed in each statement and proof and  the labeling of these constants starts anew in each proof unless they are specified to denote particular values. 
We use the symbol ``$:=$'' to denote a definition, 
which is read as ``is defined to be.''  
We write $a\wedge b:=\min\{a,b\}$, $a\vee b:=\max\{a,b\}$.
$|A|$ denotes the $d$-dimensional Lebesgue measure of $A$.  The notation $f(x) \asymp g(x)$ means that there exist constants $c_2 \ge c_1>0$ such that $c_1g(x)\leq f (x)\leq c_2 g(x)$ for a specified range of the variable $x$.

	\section{Preliminaries}\label{s:2} 
	
	In the remainder of this paper, we always assume that $\beta>0$ and $\kappa>0$. 
	We also assume that $V$ is a non-negative locally bounded function on $\R^d_0$ such that
	\begin{align}\label{e:V-bounded}
		\sup_{x\in \R^d \setminus B(0,1)}	V(x) <\infty
	\end{align}
        and that 
	\begin{align}\label{e:V-origin}	\int_{B(0,1)\setminus \{0\}} V(x) dx = \infty.	\end{align}

	Let $W=\{W_t, t\ge 0; \, \P_x, x \in \R^d\}$ be a  Brownian motion on $\R^d$ with generator $\Delta$. For an open subset $U$ of $\R^d$, let $\tau^W_U:=\inf\{t>0:\, W_t\notin U\}$ be  the first exit time from $U$ for $W$. The killed Brownian motion $W^U$ on $U$ is defined by $W^U_t=
	W_t$ 
	if $t<\tau^W_U$ and $W^U_t=\partial$ if $t\ge \tau^W_U$,
	where 
	$\partial$ is the cemetery point.  It is known that the process $W^U$ has a jointly continuous transition density 	$q^U:(0,\infty)\times U \times U \to [0,\infty)$.

	For any  non-negative Borel function $f$ on $\R^d$, define
	\begin{align}\label{e:Feynman-Kac-0}
		P_t f(x) =  \E_x \left[ e^{-  \int_0^t   
		V(W_s) ds}  
		f(W_t) \,\right], \quad t\ge 0, x\in \R^d.
	\end{align}
	Since $V$ is non-negative, we can extend  $(P_t)_{t\ge 0}$ to a strongly continuous Markov semigroup on  $L^2(\R^d)$.
	Let $X$ be the Markov process on $\R^d$ associated with the semigroup $(P_t)_{t\ge 0}$. The generator of $X$ is $\Delta-V$.

	We claim that $X$ can be  considered as a process on $\R^d_0$ and it can be obtained by killing the killed Brownian motion $W^{\R^d_0}$ 
	using the potential $V$: For this, it suffices to show that for all $t>0$ and $x\in\R^d_0$,  
	\begin{align}\label{e:Feynman-Kac-1}
		P_tf(x)= \E_x \Big[ e^{-  \int_0^t   
		V(W^{\R^d_0}_s)ds }   f(W^{\R^d_0}_t) \,\Big].
	\end{align} 
	When $d\ge 2$, we have
	$\P_x(\tau^W_{\R^d_0}=\infty)=1$. Hence, \eqref{e:Feynman-Kac-1} follows from \eqref{e:Feynman-Kac-0}. Suppose that $d=1$. In this case,  
	since $V$ satisfies \eqref{e:V-origin}, 
	 by \cite[Theorem 1]{Eng-Sch81} or \cite[Lemma 1.6]{Zanzotto97},   we have
	$$
	\P_0\Big(\, \int^t_0 
	V(W_s)ds 
	=\infty \mbox{ for all } t>0\Big)=1.
	$$
	Thus, by the strong Markov property,  for any $x\in \R^1_0$, 
	$$
	\P_x\Big(\, \int^t_0 
	V(W_s)ds
	=\infty \mbox{ for all } t>
	\tau^W_{\R^1_0}\Big)=1.
	$$
	Therefore, we deduce that for any $t>0$ and $x \in \R^1_0$,
	\begin{align*}
		\E_x\Big[e^{-\int^{t}_0 
		V(W_s)ds}f(W_t) : t<\tau^W_{\R^1_0}\Big]=\E_x\Big[e^{-\int^{t}_0
		V(W^{\R^1_0}_s)ds}f(W^{\R^1_0}_t) \Big],
	\end{align*} 
	proving that \eqref{e:Feynman-Kac-1} is also valid for $d=1$. From now on, we always regard $X$ as a process on $\R^d_0$.

	For an open subset $U$ of $\R^d_0$, we denote by $\tau_U=\tau^{X}_U$ the first exit time from $U$ for $X$ and by $X^{U}$ the killed process of $X$ on $U$. 
	Let $U$ be an open subset of $\R^d_0$ such that $\overline U \subset \R^d_0$. Note that by \eqref{e:V-bounded}, 
	$$
	\lim_{t\to 0}\sup_{x\in U}\left| \int^t_0\int_U q^U(s, x, y) 
	V(y) 
	dyds\right| \le \lim_{t\to 0} \, t\sup_{y\in U} 
	V(y)=0.
	$$
	Hence, by general theory, the semigroup $(P^{U}_t)_{t\ge 0}$ of $X^{U}$ can be represented by the following Feynman-Kac formula: For any $t\ge 0$, $x \in U$ and any non-negative Borel function $f$ on $U$,
	\begin{align}\label{e:Feynman-Kac-2}
		P^{U}_t f(x) &=  \E_x \left[ e^{-  \int_0^t   
		V(W^U_s)ds }
	  f(W^U_t) \,\right].
	\end{align}
	Define $p^{U}_0(t,x,y) := q^U(t, x, y)
	$ and, for $k \ge 1$, 
	\begin{align*}
		p^{U}_{k}(t,x,y) =&-  \int_0^t \int_U q^U(s, x, z)
		p^{U}_{k-1}(t-s,z,y) V(z) dzds.
	\end{align*}
	Set $p^{U}(t,x,y) := \sum_{k=0}^\infty p^{U}_{k}(t,x,y).$   
	It is well known (see, for instance, \cite{BM}) that 
	$p^{U}(t,x,y) $ is jointly continuous on $(0,\infty) \times U \times U$ and is the transition density for $P^{U}_t$. 
	
	Let $B(x,r)$ and $\overline B(x,r)$ denote the open and closed balls of radius $r$ centered at $x$, respectively. For each $t>0$ and $x,y \in \R^d_0$,  $p^{\overline B(0,1/n)^c}(t,x,y)$ is  increasing as $n$ increases.  Define
	$$
	p(t,x,y) :=\lim_{n \to \infty} p^{\overline B(0,1/n)^c} (t,x,y), \qquad t>0, \; x,y \in \R^d_0.
	$$
By the monotone convergence theorem,   $p(t,x,y)$ is the transition density of $X$. 
	Note that $p(t,x,y)=p(t,y,x) \le q(t,x,y)$ for all $t>0$ and $x,y \in \R^d_0$. Moreover, as an increasing limit of continuous functions, for any fixed $t>0$ and $z \in \R^d_0$, the maps $x\mapsto p(t,x,z)$ and $y \mapsto p(t,z,y)$ are lower semi-continuous. When $d=1$, since $X$ has continuous sample paths, we see that $p(t,x,y)=0$ for all $t>0$ and $x,y \in \R^1_0$ with $xy<0$.

	We end this section by recalling the two-sided heat kernel estimates for killed Brownian motions in 
	exteriors of balls. Recall that $	\Lg r = \log(e-1+r)$ for $r\ge 0$.
	For  $R>0$, we define 
	$\psi_{d,R} : (0,\infty) \times  [R, \infty) \to (0,1]$ by
	\begin{align}\label{e:def-psi-R}
\psi_{d,R} (t,r) := \begin{cases}
			\displaystyle 1 \wedge \frac{(r-R)\wedge R}{\sqrt t\wedge R} &\mbox{ if $d\ge 3$}, \\[10pt]
			\displaystyle   1 \wedge \frac{((r-R) \wedge R) \,\Lg ( (r-R)/R)}{(\sqrt t \wedge R)\, \Lg ( \sqrt t/R)} &\mbox{ if $d=2$}, \\[10pt]
			\displaystyle  1 \wedge \frac{r-R}{\sqrt t} &\mbox{ if $d= 1$}.
		\end{cases}
	\end{align}

	\begin{prop}\label{p:DHKE}
		(i) If $d\ge 2$, then	there  exist $c_i>0$, $1\le i\le 4$, depending only on $d$ such that for all  $R>0$, $t>0$ and $x,y \in \overline B(0,R)^c$,
		\begin{equation}\label{e:DHKE-1}
	\!\!	c_1 q(c_2t,x,y){\psi_{d,R}(t,|x|)\psi_{d,R}(t,|y|) } \!\le \! q^{\overline B(0,R)^c}(t,x,y)\!
		\le
			\! c_3	q(c_4t,x,y){\psi_{d,R}(t,|x|)\psi_{d,R}(t,|y|) }.
		\end{equation}
		(ii) If $d=1$, then   \eqref{e:DHKE-1} holds for all $R>0$, $t>0$ and $x,y \in \overline B(0,R)^c$ with $xy>0$.
	\end{prop}
	\pf By the scaling property, we can assume $R=1$.
	
	\noindent
	(i) For $d\ge 3$,  the result follows from \cite[Theorem 1.1]{Zh03}. For $d=2$,  the resut follows from \cite[Theorem 5.14 and 5.16]{GS11} and the fact that the function  $h$ therein can be chosen as $h(r)=\log(1+r)$. See also \cite[Example 1]{GS02}.
	
		\noindent
	(ii) For $t>0$ and $x,y>1$, we have 
$
		q^{\overline B(0,1)^c}(t,x,y) = q(t,x,y) - q(t,x, 2-y).
$
	From this, we obtain \eqref{e:DHKE-1}. Since $q^{\overline B(0,1)^c}(t,x,y) = q^{\overline B(0,1)^c}(t,-x,-y)$ for all $t>0$ and $x,y \in \overline B(0,1)^c$ by symmetry, this completes the proof.\qed

	\section{Survival probability estimates}\label{s:3}
	
In this section we establish sharp two-sided survival probability  estimates that
give the precise boundary behaviors of the heat kernel $p(t,x,y)$. 
We first introduce the classes of supercritical killing potentials 
that we will deal with.   
Recall that we always assume that
$V$ is a locally bounded non-negative Borel function on $\R^d_0$ satisfying \eqref{e:V-bounded} and \eqref{e:V-origin}.
	\begin{defn}
		\rm	 Let   $\beta,\kappa>0$ be  constants.
		
		\smallskip

		\noindent(i) We say that $V$ belongs to the class  $\Kge$ if there exist constants $\beta'\in (0,\beta)$ and $C_1>0$ such that
		\begin{align}\label{e:V-cond-1}
			V(x) \ge \frac{\kappa}{|x|^{2+2\beta}} -  \frac{C_1}{|x|^{2+\beta'}}  \quad \text{for all} \;\, x \in \R^d_0 \text{ with } |x|\le 1.
		\end{align}

			\noindent(ii) We say that $V$ belongs to the class  $\Kle$ if there exist constants $\beta'\in (0,\beta)$ and $C_2>0$ such that
		\begin{align}\label{e:V-cond-2}
			V(x) \le \frac{\kappa}{|x|^{2+2\beta}} +  \frac{C_2}{|x|^{2+\beta'}}  \quad \text{for all} \;\, x \in \R^d_0 \text{ with } |x|\le 1.
		\end{align}
		
		\noindent (iii)   We say that $V$ belongs to the class  $\Kloc$ if  $V \in \Kge \cap \Kle$, that is, there exist a Borel function $c:\overline{B(0,1)}\to [0,\infty)$ and constants $C>0$ and $\theta>\beta$  such that
		\begin{align*}
			V(x) = \frac{c(x)}{|x|^{2+2\beta}} \quad \text{and}  \quad 	|c(x)-\kappa|\le C|x|^\theta   \quad \text{for all} \;\, x \in \R^d_0 \text{ with } |x|\le 1.
		\end{align*}

		\noindent (iv) 	 We say that $V$ belongs to the class  $\KK$ if  $V \in \Kloc$ and there exist $C_3,\gamma>0$ such that  
		\begin{align}\label{e:V-cond-3}
			V(x) \le 	\frac{C_3}{|x|^{2+\gamma}}  \quad \text{for all} \;\, x \in \R^d_0 \text{ with } |x| > 1.
		\end{align}
	\end{defn}

	We first consider the case that  
	the supercritical killing potential is $\kappa|x|^{-2-2\beta}$ and give an explicit harmonic function for
	$\Delta -\kappa |x|^{-2-2\beta}$.
	\medskip

	For any $\beta,\kappa>0$, 
	recall that the function $h_{\beta,\kappa}$ defined in \eqref{e:function-h} and define another function $\wt h_{\beta,\kappa}$ on $(0,\infty)$ by
	\begin{align}\label{e:function-wth}
			\wt h_{\beta,\kappa}(r)&:= r^{-(d-2)/2} \sK_{(d-2)/(2\beta)} \Big( \frac{\sqrt \kappa }{\beta r^{\beta}}  \Big),
	\end{align}
	where $\sK_{\nu}$ is the modified Bessel function of the second 
	kind, see \cite{Wa44}.

	It is known that for any $\nu\in \R$,
	\begin{align}\label{e:asymptotic-K-nu}
		\lim_{r\to \infty}  \sqrt r e^r K_\nu(r)= \sqrt{{\pi}/{2}}.
	\end{align}
	See \cite[9.7.2]{AS64}. Further,	by the recurrence formula 
	\cite[p.79, (3)]{Wa44}, 
	for any $\nu\in \R$ and $r>0$,
	\begin{align}\label{e:recurrence-K-nu}
		K_\nu'(r) + \frac{\nu}{r} K_\nu(r) = -K_{\nu-1}(r)<0.
	\end{align}
Using \eqref{e:recurrence-K-nu} and the fact that $\sK_{\nu}$ is positive, 	
	we get that for all $r>0$,
	\begin{align*}
		\wt h_{\beta,\kappa}'(r) &= - \frac{\sqrt \kappa}{r^{d/2+\beta}}\left[ \frac{(d-2)r^\beta}{2\sqrt \kappa} \sK_{(d-2)/(2\beta)} \Big( \frac{\sqrt \kappa }{\beta r^{\beta}}  \Big) +  \sK_{(d-2)/(2\beta)}' \Big( \frac{\sqrt \kappa }{\beta r^{\beta}}  \Big) \right]\nn\\
		&=\frac{\sqrt \kappa}{r^{d/2+\beta}} \sK_{(d-2-2\beta)/(2\beta)} \Big( \frac{\sqrt \kappa }{\beta r^{\beta}} \Big) >0.
	\end{align*}
	Hence, $\wt h_{\beta,\kappa}$ is increasing in $(0,\infty)$
and by	\eqref{e:asymptotic-K-nu},
\begin{align}
\label{e:derivative-est}
{c_1}{r^{-1-\beta}} \wt h_{\beta,\kappa}(r)\le \wt h_{\beta,\kappa}'(r)  \le 
 {c_2}{r^{-1-\beta}}\wt h_{\beta,\kappa}(r)
 \quad \text{for all $r\in (0,1]$.}
\end{align} 
	  Using \eqref{e:asymptotic-K-nu}, 
	 we see that 
	 the two functions $\wt h_{\beta,\kappa}$ and  $h_{\beta,\kappa}$ 
	 are comparable on $(0,1]$. That is, for any fixed $\beta,\kappa>0$,    there exists $C>1$ such that
	 \begin{align}\label{e:h-kappa-compare} 
	 	C^{-1}h_{\beta,\kappa}(r)\le  \wt	h_{\beta,\kappa}(r) \le C h_{\beta,\kappa}(r) \quad \text{for all} \;\, r\in (0,1].
	 \end{align}   
	 Thus,	 $h_{\beta,\kappa}$ is almost increasing in $(0,\infty)$, namely,  there exists $C> 1$ such that $h_{\beta,\kappa}(s)\le Ch_{\beta,\kappa}(r)$ for all $0<s<r$. 
	 
	 For any $a,b>0$, we have 
	 \begin{align}\label{e:exp-poly}
	 	\sup_{r>0} r^a e^{-br}= (a/b)^a e^{-a}<\infty.
	 \end{align}
	 
The following lemma, which is a simple consequence of \eqref{e:exp-poly}, will be used in Section \ref{s:4} several times.
	 \begin{lem}\label{l:h-ratio}
	 	Let $\beta,\kappa>0$. For any $\eps \in (0,1)$ and $R\ge 1$, there exists $C>1$ such that
	 	\begin{equation}\label{e:h-ratio}
	 	C^{-1} \exp \Big(- \frac{(1+\eps)\sqrt \kappa}{\beta r^\beta} \Big)\le 	h_{\beta,\kappa}(r) \le C \exp \Big(- \frac{(1-\eps)\sqrt \kappa}{\beta r^\beta} \Big) \quad \text{ for all $0<r\le R$}.
	 	\end{equation}
	 \end{lem}
	 \pf  For $r\in (0,1]$, if $d-2-\beta\ge 0$, then  the lower bound in \eqref{e:h-ratio} is evident and the upper bound follows from \eqref{e:exp-poly}. If  $d-2-\beta<0$, then  the upper bound in \eqref{e:h-ratio} is evident and the lower bound follows from \eqref{e:exp-poly}. 
	  For $r \in (1,R]$, by taking $C$ larger than $e^{\sqrt \kappa/\beta}$,  we get \eqref{e:h-ratio}. \qed

	Let 
	\begin{align*}
		L^{\beta,\kappa}:=\Delta - \frac{\kappa}{|x|^{2+2\beta}}\qquad \text{ and }\qquad L^V:=\Delta - V.
	\end{align*}
	
	By \cite[Lemmas 2.1 and 2.2]{LZ17}, we have
\begin{lem}\label{l:wt-h-kappa}
	For any $\beta,\kappa>0$,	the function $x\mapsto \wt h_{\beta,\kappa}(|x|)$ is smooth in $\R^d_0$ and	$L^{\beta,\kappa} 
(\wt h_{\beta,\kappa} (|\cdot|))(x)
	=0$ for all $x \in \R^d_0$.
\end{lem}

The following explicit estimate of  the action of the operator $L^{\beta,\kappa}$ on $|x|^a \wt h_{\beta,\kappa} (|x|)$ will be 
an important ingredient in 
constructing  suitable barrier functions for $ L^V$.

\begin{lem}\label{l:wt-h-kappa-2}
Let $\beta,\kappa>0$. There exist $c_1,c_2>1$  such that for any $a\in \R$,
	\begin{align*}
 \frac{c_1a}{|x|^{2+\beta-a}} + \frac{a (a+d-2)  }{|x|^{2-a}} \le 
 \frac{L^{\beta,\kappa} (|\cdot|^a \wt h_{\beta,\kappa} (|\cdot|)) (x)}{\wt h_{\beta,\kappa}(|x|) } 
  \le  
\frac{c_2a}{|x|^{2+\beta-a}} + \frac{a (a+d-2)  }{|x|^{2-a}}  
	\end{align*}
	for all  $x \in \R^d_0$ with $|x|\le 1$.
\end{lem}
\pf Define $f(r)=r^a \wt h_{\beta,\kappa}(r)$ for $r>0$.  For any $x \in \R^d_0$ with $|x|=r$, using  Lemma \ref{l:wt-h-kappa}, we get
\begin{align}\label{e:wt-h-kappa-2}
		&L^{\beta,\kappa} (f (|\cdot|)) (x)= L^{\beta,\kappa} (f (|\cdot|))  (x)- |x|^a L^{\beta,\kappa} (\wt h_{\beta,\kappa} (|\cdot|))(x) \nn\\
		&= f''(r) + \frac{d-1}{r} f'(r)  - \frac{\kappa}{r^{2+2\beta}} f(r)  - r^a \Big( \, \wt h_{\beta,\kappa}''(r) + \frac{d-1}{r}  \wt h_{\beta,\kappa}'(r)  - \frac{\kappa}{r^{2+2\beta}}  \wt h_{\beta,\kappa}(r)\Big)\nn\\
		& =\frac{ a (a+d-2)}{r^2}  f(r)  + 2 a r^{a-1} \wt h_{\beta,\kappa}'(r).
\end{align}
Combining \eqref{e:wt-h-kappa-2} with \eqref{e:derivative-est}, we arrive at the result.
\qed

We now return to general supercritical potentials $V$ and construct appropriate barrier functions for $L^V$.
Suppose $V\in \Kge$  or $V\in \Kle$ 
with $\beta'\in (0,\beta)$  being the constant in \eqref{e:V-cond-1} or  \eqref{e:V-cond-2}.
Define 
\begin{align}\label{e:def-u-1}
u_1(r)=
(2-r^{(\beta-\beta')/2} )\wt h_{\beta,\kappa}(r) \quad \text{and} \quad  u_2(r)=(1+r^{(\beta-\beta')/2} )\wt h_{\beta,\kappa}(r), \quad\; r >0.
\end{align} 
\begin{lem}\label{l:generator}
\noindent (i) Suppose $V\in \Kge$ 
with $\beta'\in (0,\beta)$ being the constant in \eqref{e:V-cond-1}.
Then $u_1$ satisfies the following properties:
\smallskip

(a) 
$\wt h_{\beta,\kappa}(|x|)\le u_1(|x|) \le 2\wt h_{\beta,\kappa}(|x|)$ 
for all $x \in \R^d_0$ with $|x|\le 1$.

(b) There exist constants $R_1\in (0,1]$ and $C>0$ such that 
\begin{align*}
	L^V (u_1(|\cdot|) )(x)
	\le 0 \quad \text{for all $x \in \R^d_0$ with $|x|\le R_1$}.
\end{align*}

\noindent(ii) Suppose $V\in \Kle$ 
with $\beta'\in (0,\beta)$ being the constant in \eqref{e:V-cond-2}.
Then $u_2$ satisfies the following properties:
\smallskip

(a) 
$\wt h_{\beta,\kappa}(|x|)\le u_2(|x|) \le 2\wt h_{\beta,\kappa}(|x|)$ 
for all $x \in \R^d_0$ with $|x|\le 1$.

(b) There exist constants $R_2\in (0,1]$ and $C>0$ such that 
\begin{align*}
		L^V (u_2(|\cdot|))(x) 
	\ge 0 \quad \text{for all $x \in \R^d_0$ with $|x|\le R_2$}.
\end{align*}
\end{lem}
\pf (i) Set $a:=(\beta-\beta')/2$. (a) is immediate. For all $x \in \R^d_0$ with $|x|\le 1$, using \eqref{e:V-cond-1} in the first inequality below, and Lemmas \ref{l:wt-h-kappa} and \ref{l:wt-h-kappa-2} in the third, we get
\begin{align*}
		L^V(u_1(|\cdot|))(x) &\le 2\Delta (\wt h_{\beta,\kappa}(|\cdot|) )(x) - \Delta (|\cdot|^a \wt h_{\beta,\kappa}(|\cdot|) )(x)\nn\\
	&\quad  - \frac{\kappa}{|x|^{2+2\beta}} (2\wt h_{\beta,\kappa}(|x|) - |x|^{a} \wt h_{\beta,\kappa}(|x|) ) + \frac{C_1}{|x|^{2+\beta'}}(2\wt h_{\beta,\kappa}(|x|) - |x|^{a} \wt h_{\beta,\kappa}(|x|) )\nn\\
	&\le 	2L^{\beta,\kappa} (\wt h_{\beta,\kappa}(|\cdot|))(x) - L^{\beta,\kappa}( |\cdot|^a \wt h_{\beta,\kappa}(|\cdot|))(x)
	 + \frac{2C_1}{|x|^{2+\beta'}}\wt h_{\beta,\kappa}(|x|) \nn\\
	&\le -  \left( c_1a +   a (a+d-2) |x|^{\beta} - 2C_1|x|^{a}\right) \frac{\wt h_{\beta,\kappa}(|x|)}{|x|^{2+\beta-a}}.
\end{align*}
Taking $R_1\in (0,1]$ satisfying $ a |a+d-2| R_1^{\beta} + 2C_1 R_1^a \le c_1a$, we arrive at the  result.

(ii) Set $a:=(\beta-\beta')/2$ as in (i). (a) is evident. For all $x \in \R^d_0$ with $|x|\le 1$, using \eqref{e:V-cond-2}, and Lemmas \ref{l:wt-h-kappa} and \ref{l:wt-h-kappa-2}, we see that 
\begin{align*}
L^V(u_2(|\cdot|)) (x)
	&\ge 	\Delta( \wt h_{\beta,\kappa}(|\cdot|) ) (x) + \Delta (|\cdot|^a \wt h_{\beta,\kappa}(|\cdot|) ) (x)\nn\\
	&\quad  - \frac{\kappa}{|x|^{2+2\beta}} (\wt h_{\beta,\kappa}(|x|) + |x|^{a} \wt h_{\beta,\kappa}(|x|) ) - \frac{C_2}{|x|^{2+\beta'}}(\wt h_{\beta,\kappa}(|x|) + |x|^{a} \wt h_{\beta,\kappa}(|x|) )\nn\\
	&\ge  	L^{\beta,\kappa} (\wt h_{\beta,\kappa}(|\cdot|)) (x)+ L^{\beta,\kappa}( |\cdot|^a \wt h_{\beta,\kappa}(|\cdot|))  (x)
	- \frac{2C_2}{|x|^{2+\beta'}}\wt h_{\beta,\kappa}(|x|) \nn\\
	&\ge   \left( c_2a -   a (a+d-2) |x|^{\beta} - 2C_2|x|^{a}\right) \frac{\wt h_{\beta,\kappa}(|x|)}{|x|^{2+\beta-a}}.
\end{align*}
Choosing $R_2\in (0,1]$ such that $ a |a+d-2| R_2^{\beta} + 2C_2 R_2^a \le c_2a$, this yields the  result. \qed

Denote by $\zeta$ the lifetime of $X$.
Recall that $u_1$ and $u_2$ are the barrier functions defined in \eqref{e:def-u-1}.
Using Ito's formula and Lemma \ref{l:generator}, 
we obtain the following bounds of the probability of exiting balls before the lifetime in terms of the barrier functions.
\begin{lem}\label{l:survival-probability}
(i) 
Suppose $V\in \Kge$ with $\beta'\in (0,\beta)$ being the constant in \eqref{e:V-cond-1}.
Then
	\begin{align*}
\P_x \left(\tau_{B(0,R)} <\zeta \right) \le \frac{u_1(|x|)}{ u_1(R)}  \quad \text{ for all  $x\in \R^d_0$ with $|x|<R\le R_1$,}
	\end{align*}
	where $R_1\in (0,1]$ is the constant in Lemma \ref{l:generator}(i).
	
\noindent	(ii) 
Suppose $V\in \Kle$ with $\beta'\in (0,\beta)$ being the constant in \eqref{e:V-cond-2}.
Then
\begin{align*}
	\P_x \left(\tau_{B(0,R)} <\zeta \right) \ge \frac{u_2(|x|)}{ u_2(R)} \quad \text{ for all  $x\in \R^d_0$ with $|x|<R\le R_2$,}
\end{align*}
where $R_2\in (0,1]$ is the constant in Lemma \ref{l:generator}(ii).
\end{lem}
\pf  Let $R\in (0,1]$. Since $X$ has continuous sample paths,  for all $x \in \R^d_0$ with $|x|<R$,
\begin{align}\label{e:survival-probability-0}
	\P_x \left(\tau_{B(0,R)} <\zeta \right)= \P_x \Big(X_{\tau_{B(0,R)} } \in \partial B(0,R) \Big) = \lim_{\eps \to0}\P_x \Big(X_{\tau_{B(0,R)\setminus \overline B(0,\eps)} } \in \partial B(0,R) \Big).
\end{align}

(i)  Let $R\le R_1$ and take $\eps \in (0,R)$. For any  $x \in B(0,R) \setminus \overline{B(0,\eps)}$,  we have
 \begin{align}\label{e:survival-probability-1}
 	\E_x\Big[ u_1\big(|X_{\tau_{B(0,R)\setminus \overline B(0,\eps)} }|\big)\Big] &\ge u_1(R) \,\P_x \Big(X_{\tau_{B(0,R)\setminus \overline B(0,\eps)} } \in \partial B(0,R) \Big).
 \end{align}
On the other hand, using Ito's formula and Lemma \ref{l:generator}(i-b), we see that
\begin{align}\label{e:survival-probability-2}
	&\E_x\Big[u_1\big(|X_{\tau_{B(0,R)\setminus \overline B(0,\eps)} }|\big)\Big] =u_1(|x|) + \E_x \int_0^{\tau_{B(0,R)\setminus \overline B(0,\eps)}} 	L^V (u_1( |\cdot|)) (X_s)ds \le u_1(|x|).
\end{align}
Combining \eqref{e:survival-probability-0}, \eqref{e:survival-probability-1} and \eqref{e:survival-probability-2}, we arrive at the result.

 (ii) Let $R\le R_2$ and take $\eps \in (0,R)$. We have
\begin{align}\label{e:survival-probability-3}
	\E_x\Big[u_2\big(|X_{\tau_{B(0,R)\setminus \overline B(0,\eps)} }|\big)\Big] \le  u_2(R) \, \P_x \Big(X_{\tau_{B(0,R)\setminus \overline B(0,\eps)} } \in \partial B(0,R) \Big)  + u_2(\eps)
\end{align}
and, by Ito's formula and Lemma \ref{l:generator}(ii-b),
\begin{align}\label{e:survival-probability-4}
	\E_x\Big[u_2\big(|X_{\tau_{B(0,R)\setminus \overline B(0,\eps)} }|\big)\Big] \ge u_2(|x|).
\end{align}
By Lemma \ref{l:generator}(ii-a), we see that $\lim_{\eps \to 0}u_2(\eps)=0$. Thus, using \eqref{e:survival-probability-0},  \eqref{e:survival-probability-3} and \eqref{e:survival-probability-4}, we get  the result.\qed

In Lemma \ref{l:survival-probability-t-2} below, we will use Lemma \ref{l:survival-probability}, the inequality 
$\P_x (\zeta >  t ) \le  	\P_x ( \tau_{B(0,R)}<\zeta )   + \P_x ( \tau_{B(0,R)} >t )$ with appropriate $R=R_t$ 
to get the upper bound on the survival probability $\P_x (\zeta >  t )$. To accomplish this, the key is 
to get an upper bound on $\P_x ( \tau_{B(0,R)} >t )$.
The following lemma is the first step toward obtaining the upper bound on $\P_x ( \tau_{B(0,R)} >t )$.
\begin{lem}\label{l:survival-probability-t-1}
Suppose $V\in \Kge$ with $\beta' \in (0,\beta)$ and $C_1$ being the constants in \eqref{e:V-cond-1}.
For all $t>0$, $R\in (0,(\kappa/(4C_1))^{1/(2\beta-\beta')}\wedge 1]$ and  $x,y \in \R^d_0$ satisfying $|x|\vee  |y|<R$, we have
\begin{equation}\label{e:survival-probability-t-1-1}
	\limsup_{\delta \to0} \frac{\P_x (X_t \in B(y,\delta), \, \tau_{B(0,R)}> t)}{|B(0,\delta)|} \le  (4\pi t)^{-d/2}  \exp \bigg( - \frac{|x-y|^2}{4t}  - \frac{3\kappa t }{4R^{2+2\beta}} \bigg).
\end{equation}
	In particular,  we have
	\begin{align}\label{e:survival-probability-t-1-2}
\sup_{x\in B(0,R)\setminus\{0\}}		\P_x \left(\tau_{B(0,R)} > t \right)  \le \exp \bigg(-   \frac{3\kappa t}{4R^{2+2\beta}}  \bigg).
	\end{align}
\end{lem}
\pf From  \eqref{e:Feynman-Kac-2} and \eqref{e:V-cond-1}, since  $R\le (\kappa/(4C_1))^{1/(2\beta-\beta')} \wedge 1$, we see that
\begin{align*}
&	\limsup_{\delta \to0} \frac{\P_x (X_t \in B(y,\delta),  \, \tau_{B(0,R)}> t)}{|B(0,\delta)|}\nn\\
	&\le  \limsup_{\delta\to 0} \frac{\P_x\big( W^{B(0,R)}_t \in B(y,\delta) \big)}{|B(0,\delta)|}  \exp \bigg(- t   \inf_{r \in (0,R)} \bigg(\frac{\kappa}{r^{2+2\beta}}- \frac{C_1}{r^{2+\beta'}} \bigg)  \bigg) \nn\\
	&\le  (4\pi t)^{-d/2}  \exp \bigg( - \frac{|x-y|^2}{4t}  - \frac{\kappa t }{R^{2+2\beta}} + \frac{C_1t}{R^{2+\beta'}}\bigg)\nn\\
	&\le  (4\pi t)^{-d/2}  \exp \bigg( - \frac{|x-y|^2}{4t}  - \frac{3\kappa t }{4R^{2+2\beta}} \bigg),
\end{align*}
where in the second inequality we use the fact that the function $r\mapsto k/r^{2+2\beta}-C_1/r^{2+\beta'}$ is increasing
on the interval $(0, (((2+2\beta)\kappa)/((2+\beta')C_1))^{1/(2\beta-\beta')})$.
This proves \eqref{e:survival-probability-t-1-1}. For any $x\in B(0,R)\setminus\{0\}$, using  \eqref{e:survival-probability-t-1-1}, we get
\begin{align*}
	\P_x \left(\tau_{B(0,R)} > t \right) & \le \int_{B(0,R)}  (4\pi t)^{-d/2}  \exp \bigg( - \frac{|x-y|^2}{4t}  - \frac{3\kappa t }{4R^{2+2\beta}} \bigg) dy \nn\\
&\le \exp \bigg(-   \frac{3\kappa t}{4R^{2+2\beta}}  \bigg) \int_{\R^d} q(t,x,y)dy  = \exp \bigg(-   \frac{3\kappa t}{4R^{2+2\beta}}  \bigg),
\end{align*}
proving that  \eqref{e:survival-probability-t-1-2} holds. 
\qed

Define
\begin{align}\label{e:def-eta-0}
\eta_0=\eta_0(\beta,\kappa):= (2^{-4-\beta/(2+\beta)}\beta \sqrt \kappa)^{1/(2+\beta)}.
\end{align} 
As mentioned earlier, the key to the proof of the next lemma is get an upper bound for  $\P_x ( \tau_{B(0,R} >t )$ with $R=\eta_0 t^{1/(2+\beta)}$.

\begin{lem}\label{l:survival-probability-t-2}
Suppose $V\in \Kge$ with $\beta' \in (0,\beta)$ and $C_1>0$ being the constants in \eqref{e:V-cond-1}.
There exists $C>0$ such that 
	\begin{equation*}
		\P_x \left(\zeta >  t \right) \le \frac{Ch_{\beta,\kappa}(|x|)}{h_{\beta,\kappa}( \eta_0 t^{1/(2+\beta)} )} \quad \text{for all $t\in (0,1]$ and $x\in \R^d_0$ with $|x|<\eta_0 t^{1/(2+\beta)}$.}
	\end{equation*}
\end{lem}
\pf  Let $R_1\in (0,1]$ be the constant in Lemma \ref{l:generator}(i) and 
$$\wt R_1:=\frac{(\kappa/(4C_1))^{1/(2\beta-\beta')} \wedge R_1}{2^{1/(2+\beta)}}.$$
 Set $R(t):=\eta_0 t^{1/(2+\beta)}$. We deal with two cases separately.

\smallskip

\noindent
{\it
Case 1:  $R(t)\le \wt R_1$.}
Using \eqref{e:h-kappa-compare}, Lemmas \ref{l:survival-probability}(i) and \ref{l:generator}(i-a), 
we get
\begin{equation}\label{e:survival-probability-t-main-0}
		\P_x \left(\zeta >  t \right) \le  	\P_x \left( \tau_{B(0,R(t))}<\zeta \right)   + \P_x \left( \tau_{B(0,R(t))} >t \right)\le \frac{c_1h_{\beta,\kappa}(|x|)}{h_{\beta,\kappa}(R(t))} +  \P_x \left( \tau_{B(0,R(t))} >t \right).
\end{equation}
Define for $n \ge 0$, $
r_n:=2^{n/(2+\beta)}|x|$ and  $ \eps_n:={2^{n+1+\beta/(2+\beta) } |x|^{2+\beta}}/(\beta \sqrt \kappa \,t) $.
Observe that
\begin{align}\label{e:survival-probability-t-1}
  \frac{\kappa \eps_n t}{r_n^{2+2\beta}}= 	\frac{2\sqrt \kappa}{\beta r_{n-1}^\beta} \quad \text{for all $n \ge 1$}.
\end{align}
Let $N:=\min\{n \ge 1: r_{n+1} \ge R(t)\}$. 
Since $R(t)\le \wt R_1\le 2^{-1/(2+\beta}R_1$, we have  $r_N <R(t)\le 2^{1/(2+\beta)}R(t) \le R_1\le 1$.
Using $r_N<2^{1/(2+\beta)}R(t)$, we get that
\begin{align*}
\sum_{1\le i\le N-1}\eps_i<\eps_N	=\frac{2^{1+\beta/(2+\beta) } r_N^{2+\beta}}{\beta \sqrt \kappa \,t}  <  \frac{2^{2+\beta/(2+\beta) }\eta_0^{2+\beta}}{\beta \sqrt \kappa } =\frac14.
\end{align*}
 Using the strong Markov property, we see that for all $1\le k \le N$ and $z \in \R^d$ with $|z| = r_{k-1}$,
\begin{align*}
		&\P_z \Big( \tau_{B(0,R(t))}> \big(1-\sum_{1\le i \le k-1}\eps_i \big)t \Big) \nn\\
		&\le \P_z \left( \tau_{B(0,r_k)} > \eps_k t\right) + \E_z \bigg[ \P_{X_{\tau_{B(0,r_k)}}} \Big(   \tau_{B(0,R(t))}> \big(1-\sum_{1\le i \le k}\eps_i \big)t \Big)   ; \tau_{B(0,r_k)} \le  \tau_{B(0,R(t))} \wedge (\eps_k t)\bigg]  \nn\\
	&\le \P_z \left( \tau_{B(0,r_k)} > \eps_k t\right) +\P_z \left( \tau_{B(0,r_k)}< \zeta\right)  \sup_{|w|=r_k}\P_w \Big( \tau_{B(0,R(t))} > \big(1-\sum_{1\le i \le k}\eps_i \big)t \Big).
\end{align*}
Hence, 
\begin{align}\label{e:survival-probability-t-main-1}
&\P_x \left( \tau_{B(0,R(t))} >t \right) \le \sup_{|z_0| = r_0} \P_{z_0} \left( \tau_{B(0,R(t))} > t \right)  \nn\\
		&\le \sup_{|z_0| = r_0} \bigg[ \P_{z_0} \left( \tau_{B(0,r_1)} > \eps_1 t\right) +\P_{z_0} \left( \tau_{B(0,r_1)}< \zeta \right)  \sup_{|z_1|=r_1}\P_{z_1} \left( \tau_{B(0,R(t))} > \big(1- \eps_1\big)t \right) \bigg]\nn\\[2pt]
		 &\le \cdots \nn\\[3pt]
		  &\le \sum_{k=0}^{N-1}  \bigg( \prod_{0\le i\le {k-1}}  \sup_{|z_i|=r_i} \P_{z_i} \left( \tau_{B(0,r_{i+1})}< \zeta \right)  \bigg)  \sup_{|z_k|=r_k}\P_{z_k}  \left(  \tau_{B(0,r_{k+1})} > \eps_{k+1}t \right)\nn\\
		  &\quad +  \bigg( \prod_{0\le i\le {N-1}}  \sup_{|z_i|=r_i} \P_{z_i} \left( \tau_{B(0,r_{i+1})}< \zeta \right)  \bigg)  \sup_{|z_N|=r_N}\P_{z_N}  \Big( \tau_{B(0,R(t))}> \big(1-\sum_{1\le i \le N}\eps_i \big)t \Big) \nn\\
		  &=:I_{1}+I_{2}.
\end{align}
For each $0\le k \le N-1$, using Lemmas \ref{l:survival-probability}(i) and  \ref{l:generator}(i-a) and  \eqref{e:survival-probability-t-1-2},   we obtain
\begin{align*}
	&\bigg( \prod_{0\le i\le {k-1}}  \sup_{|z_i|=r_i} \P_{z_i} \left( \tau_{B(0,r_{i+1})}< \zeta \right)  \bigg)  \sup_{|z_k|=r_k}\P_{z_k}  \left( \tau_{B(0,r_{k+1})} > \eps_{k+1}t \right) \nn\\
		&\le \frac{u_1(r_0)}{ u_1 (r_{k})} \exp \bigg( - \frac{3\kappa \eps_{k+1}t}{4r_{k+1}^{2+2\beta}}  \bigg)
	\le \frac{c_2h_{\beta,\kappa}(r_0)}{ h_{\beta,\kappa} (r_{k})} \exp \Big( - \frac{3\kappa \eps_{k+1}t}{4r_{k+1}^{2+2\beta}}  \Big).
\end{align*}
Using \eqref{e:survival-probability-t-1}  and  the fact that $\sup_{0<r\le 1} r^{(d-4-\beta)/2} e^{-\sqrt \kappa /(2\beta r^\beta)} <\infty$, we see that 
\begin{align*}
	    \exp \Big( - \frac{3\kappa \eps_{k+1}t}{4r_{k+1}^{2+2\beta}}  \Big)  &=    h_{\beta,\kappa}(r_k) r_k^{(d-2-\beta)/2}  \exp \Big( - 	\frac{\sqrt \kappa}{2\beta r_{k}^\beta}  \Big)  \le  c_3h_{\beta,\kappa}(r_k) r_k.
\end{align*}
Thus, 
\begin{align}\label{e:survival-probability-t-main-2}
I_{1}&\le c_2c_3h_{\beta,\kappa}(r_0)\sum_{k=0}^N r_k  \le c_4h_{\beta,\kappa}(r_0)r_N \le  c_4 h_{\beta,\kappa}(r_0).
\end{align}
On the other hand, using  $\sum_{1\le i \le N}\eps_i  < 2\eps_{N}<1/2$ and Lemma  \ref{l:survival-probability}(i) in the first inequality below, Lemma \ref{l:generator}(i-a) and \eqref{e:survival-probability-t-1-2} in the second,  $2^{-1/(2+\beta)} <r_N/R(t)< 2^{1/(2+\beta)}$ and $R(t)\le \wt R_1$ in the third, and $\sup_{0<r\le 1} r^{(d-2-\beta)/2} e^{-2^{\beta/(2+\beta)}\sqrt \kappa /(\beta r^\beta)} <\infty$ in the fourth, we get
\begin{align}\label{e:survival-probability-t-main-3}
I_{2}&\le \frac{u_1(r_0)}{u_1(r_{N})} \sup_{|z_N|=r_N}\P_{z_N}  \left( \tau_{B(0,R(t))} \ge  t/2 \right) \nn\\
&\le  c_5 h_{\beta,\kappa} (r_0)r_N^{(d-2-\beta)/2}\exp \Big(  \frac{\sqrt \kappa }{\beta r_{N}^{\beta}} - \frac{3\kappa t}{8 R(t)^{2+2\beta}} \Big) \nn\\
&\le  c_6 h_{\beta,\kappa} (r_0)  R(t)^{(d-2-\beta)/2}\exp \left(  \frac{2^{\beta/(2+\beta)}\sqrt \kappa }{\beta R(t)^{\beta}} - \frac{3\kappa t}{8R(t)^{2+2\beta}}\right)  \nn\\
&\le  c_7  h_{\beta,\kappa} (r_0)\exp \left(  \frac{2^{1+\beta/(2+\beta)}\sqrt \kappa }{\beta R(t)^{\beta}} - \frac{3\kappa t}{8 R(t)^{2+2\beta}}\right) \nn\\
&= c_7   h_{\beta,\kappa} (r_0)\exp \left(  - \frac{ 2^{2+\beta/(2+\beta)}\sqrt\kappa }{\beta  R(t)^{\beta}}\right) \le c_7h_{\beta,\kappa}(r_0).
\end{align}
Combining \eqref{e:survival-probability-t-main-0},  \eqref{e:survival-probability-t-main-1}, \eqref{e:survival-probability-t-main-2} and \eqref{e:survival-probability-t-main-3},  we arrive at the result in this case.

\smallskip

\noindent
{\it
Case 2: $R(t)>\wt R_1$.} If $|x|\ge \wt R_1$, then using the almost monotonicity of $h_{\beta,\kappa}$, we get
$$
\P_x(\zeta>t) \le 1 \le \frac{c_8h_{\beta,\kappa}(\eta_0) h_{\beta,\kappa}(|x|)}{h_{\beta,\kappa}(\wt R_1)h_{\beta,\kappa}(\eta_0 t^{1/(2+\beta)})}.
$$
If $|x|<\wt R_1$, then  by the 
conclusion of {\it Case 1} 
and the almost monotonicity of $h_{\beta,\kappa}$, we have
\begin{align*}
	\P_x(\zeta>t)  \le \P_x(\zeta >  (\wt R_1/\eta_0)^{2+\beta}) \le \frac{c_9h_{\beta,\kappa}(|x|)}{h_{\beta,\kappa}(\wt R_1)} \le \frac{c_{10} h_{\beta,\kappa}(\eta_0)h_{\beta,\kappa}(|x|)}{h_{\beta,\kappa}(\wt R_1)h_{\beta,\kappa}(\eta_0 t^{1/(2+\beta)})}.
\end{align*}
The proof is complete. \bk \qed

The last result in this section is a lower bound 
on the probability of exiting the ball $B(0,\eta_0(t/4)^{1/(2+\beta)} )$ before time $\zeta \wedge (t/3) $
in terms of $h_{\beta,\kappa}$.

\begin{lem}\label{l:survival-probability-t-3}
Suppose $V\in \Kle$ with $\beta' \in (0,\beta)$ and $C_2>0$ being the constants in \eqref{e:V-cond-2}.
There exist $t_0\in (0,1]$ and $C>0$ such that for all $t\in (0,t_0]$ and $x\in \R^d_0$ with $|x|<\eta_0(t/4)^{1/(2+\beta)}$,
	\begin{equation*}
		\P_x \left(\tau_{B(0,\eta_0(t/4)^{1/(2+\beta)} )} < \zeta \wedge (t/3) \right) \ge \frac{Ch_{\beta,\kappa}(|x|)}{h_{\beta,\kappa}( \eta_0 (t/4)^{1/(2+\beta)} )}.
	\end{equation*}
\end{lem}
\pf Define
\begin{align*}
	\wt V(x)= \1_{\overline{B(0,1)}}(x)(\kappa |x|^{-2-2\beta} + C_2|x|^{-2-\beta'}) +  \1_{\overline{B(0,1)}^c}(x) V(x).
\end{align*}
Then,  by \eqref{e:V-cond-2}, $\wt V \ge V$ in $\R^d_0$. Let $Y$ be the process corresponding to the generator $\Delta - \wt V$. Let $\zeta^Y$ denote the lifetime of $Y$ and  $\tau^Y_U:=\inf\{t>0: Y \notin U\}$ for an open set $U\subset \R^d_0$.  Since $\wt V \ge V$, we have  $\zeta^Y\le \zeta$.

Let $R_2$ be the constant in Lemma \ref{l:generator}(ii), $t\in (0,4(R_2/\eta_0)^{2+\beta} \wedge 1]$ and $x\in \R^d_0$ with $|x|<\eta_0(t/4)^{1/(2+\beta)}$. 
Using \eqref{e:Feynman-Kac-2}, we see that
\begin{align}\label{e:survival-probability-t-3-1}
	\P_x \left(\tau_{B(0,\eta_0(t/4)^{1/(2+\beta)} )} < \zeta \wedge (t/3) \right)& = \P_x \left(\tau^W_{B(0,\eta_0(t/4)^{1/(2+\beta)} )} <  t/3, \;  \tau^W_{B(0,\eta_0(t/4)^{1/(2+\beta)} )} <\zeta \right)\nn\\
	&\ge \P_x \left(\tau^W_{B(0,\eta_0(t/4)^{1/(2+\beta)} )} <  t/3, \;  \tau^W_{B(0,\eta_0(t/4)^{1/(2+\beta)} )} <\zeta^Y \right)\nn\\
	&= 	\P_x \left(\tau^Y_{B(0,\eta_0(t/4)^{1/(2+\beta)} )} < \zeta^Y \wedge (t/3) \right).
\end{align}
Besides, since $\wt V\in \Kloc$,  using Lemmas  \ref{l:survival-probability}(ii), \ref{l:generator}(ii-a) and \ref{l:survival-probability-t-2}, we get that 
\begin{align}\label{e:survival-probability-t-3-2}
	\P_x \left(\tau^Y_{B(0,\eta_0(t/4)^{1/(2+\beta)} )} < \zeta^Y \wedge (t/3) \right) 
	&\ge \P_x \left(\tau^Y_{B(0,\eta_0(t/4)^{1/(2+\beta)} )} < \zeta^Y \right) - \P_x \left(t/3 < \zeta^Y  \right)\nn\\
	&\ge \frac{ c_1h_{\beta,\kappa}(|x|)}{ h_{\beta,\kappa}( \eta_0 (t/4)^{1/(2+\beta)} )} - \frac{c_2h_{\beta,\kappa}(|x|)}{h_{\beta,\kappa}( \eta_0 (t/3)^{1/(2+\beta)} )}.
\end{align}
Since $h_{\beta,\kappa}$ is exponentially increasing at $0$, there exists $t_0 \in (0,4(R_2/\eta_0)^{2+\beta} \wedge 1]$ such that    $h_{\beta,\kappa}( \eta_0 (t/3)^{1/(2+\beta)} )\ge \frac{2c_2}{c_1}h_{\beta,\kappa}( \eta_0 (t/4)^{1/(2+\beta)} )$ for all $t\in (0,t_0]$. Combining  \eqref{e:survival-probability-t-3-1} and  \eqref{e:survival-probability-t-3-2}, we get the desired result. \qed

	\section{Small time heat kernel estimates}\label{s:4}
	
In this section, we derive short time explicit sharp  two-sided estimates of $p(t,x,y)$.  Recall that $\eta_0$ is defined as \eqref{e:def-eta-0}.	Define
	\begin{align*}
		\eta_1=\eta_1(\beta,\kappa):= 2^{-2/(2+\beta)}\eta_0.
	\end{align*}

		\begin{thm}\label{t:smalltime}
	(i)	Suppose $V\in \Kge$. For any $T\ge 1$, there exist $c_1,c_2,c_3>0$ such that  for all $t \in (0,T]$ and $x,y \in \R^d_0$, it holds that
	\begin{align}\label{e:smalltime-upper}
 p(t,x,y) &\le c_1\bigg( 1 \wedge  \frac{h_{\beta,\kappa}(|x| )}{h_{\beta,\kappa}( \eta_1 t^{1/(2+\beta)})}\bigg)\bigg( 1 \wedge  \frac{h_{\beta,\kappa}(|y|)}{h_{\beta,\kappa}( \eta_1 t^{1/(2+\beta)})}\bigg)  \nn\\
 &\quad \times \exp \bigg(- \frac{c_2t}{(|x| \vee |y| \vee t^{1/(2+\beta)})^{2+2\beta}}\bigg)q(c_3t,x,y).
\end{align}

 	\noindent (ii)	Suppose $V\in \Kle$ and let $T \ge 1$.
 
 \smallskip
 
 (a) If $d\ge 2$, then	 there exist  $c_4,c_5,c_6>0$ such that  for all $t \in (0,T]$ and $x,y \in \R^d_0$,
 \begin{align}\label{e:smalltime-lower}
 	p(t,x,y) &\ge c_4\bigg( 1 \wedge  \frac{h_{\beta,\kappa}(|x| )}{h_{\beta,\kappa}( \eta_1 t^{1/(2+\beta)})}\bigg)\bigg( 1 \wedge  \frac{h_{\beta,\kappa}(|y|)}{h_{\beta,\kappa}( \eta_1 t^{1/(2+\beta)})}\bigg)  \nn\\
 	&\quad \times \exp \bigg(- \frac{c_5t}{(|x| \vee |y| \vee t^{1/(2+\beta)})^{2+2\beta}}\bigg)q(c_6t,x,y).
 \end{align}
 
 (b) If  $d=1$, then \eqref{e:smalltime-lower} holds for all $t \in (0,T]$ and $x,y \in \R^1_0$ with $xy>0$.
	\end{thm} 
	
	The proofs of Theorem \ref{t:smalltime}(i) and (ii)  will be  given 
	in
	Subsections  \ref{ss:smalltime-upper} and \ref{ss:smalltime-lower} respectively.

\subsection{Small time upper estimates}\label{ss:smalltime-upper}

Throughout this subsection, we assume that 
$V\in \Kge$ with $\beta' \in (0,\beta)$ and $C_1>0$ being the constants in \eqref{e:V-cond-1}. We will simply write $h$ instead of $h_{\beta,\kappa}$. 
Let $R_1 \in (0,1]$ be the constant  in Lemma \ref{l:generator}(i) and  define
\begin{align}\label{e:R1p}
	R_1':=(\kappa/(4C_1))^{1/(2\beta-\beta')} \wedge R_1.
\end{align}

	Note that the map	$s\mapsto s^{-d/2}e^{-r^2/s}$ is increasing in $(0, 2r^2/d)$ and decreasing in $(2r^2/d,\infty)$. Hence, if $t\le R^2/(2d)$, then
\begin{align}\label{e:basic-upper-1}
	&	\sup_{0<s\le t} s^{-d/2} e^{ - {R^2}/({4s})}= 	 t^{-d/2} e^{ - {R^2}/({4t})},
\end{align}
and	if $t> R^2/(2d)$, then
\begin{align}\label{e:basic-upper-2}
		\sup_{0<s\le t} s^{-d/2} e^{ - {R^2}/{(4s)}}	&= 	 ( {R^2}/{(2d)})^{-d/2} e^{ - {d}/{2}}\le  			( {R^2}/{(2d)})^{-d/2}  e^{- {R^2}/{(4t)}}.
\end{align}

Recall that  $\zeta$ denotes the lifetime of $X$. Define 
\begin{align*}
	M_t:=\sup_{s\in [0,t] \, \cap \, [0,\zeta)} |X_s|, \quad t>0.
\end{align*}
Using the strong Markov property in the equality below, and Lemmas \ref{l:survival-probability}(i) and \ref{l:generator}(i-a), and  \eqref{e:Feynman-Kac-0} in the second inequality, we have that for all $t>0$, $R\in (0,R_1]$ and  $x,y \in \R^d_0$ satisfying $|x|<R$ and $|y| \not= R$,
\begin{align}\label{e:basic-upper-0}
	& \limsup_{\delta \to0} \frac{\P_x (X_t \in B(y,\delta),\, M_t\ge R)}{|B(0,\delta)|}\nn\\
	&= \limsup_{\delta\to 0}\frac{1}{|B(0,\delta)|}\E_x\left[ \P_{X_{\tau_{B(0,R)}}}\Big(X_{t-\tau_{B(0,R)}} \in B(y,\delta) \Big)  ;  \tau_{B(0,R)} < t\wedge \zeta \right]\nn\\
		&\le \P_x\left( \tau_{B(0,R)} <  \zeta \right) \limsup_{\delta\to 0}\frac{1}{|B(0,\delta)|} \sup_{0<s\le t,\, |z|=R} \P_{z}\Big(X_{s} \in B(y,\delta) \Big) \nn\\
		&\le \frac{c_1 h(|x|)}{h(R)} \limsup_{\delta\to 0}\frac{1}{|B(0,\delta)|} \sup_{0<s\le t,\, |z|=R} \P_{z}\Big(W_{s} \in B(y,\delta) \Big) \nn\\
			&= \frac{c_1 h(|x|)}{h(R)} \sup_{0<s\le t} (4\pi s)^{-d/2} e^{ -{(R-|y|)^2}/{(4s)}}.
\end{align}

Combining \eqref{e:basic-upper-0}, \eqref{e:basic-upper-1} and \eqref{e:basic-upper-2}, we obtain the next lemma, which holds for all $t>0$.

\begin{lem}\label{l:upper-basic-1}
There exists $C>0$ such that for all $t>0$, $R\in (0,R_1]$ and  $x,y \in \R^d_0$ satisfying 
$|x| \vee |y|<R$, it holds that
	\begin{align*}
 \limsup_{\delta \to0} \frac{\P_x (X_t \in B(y,\delta), \, M_t\ge R)}{|B(0,\delta)|} 
	\le  \frac{C h(|x|)}{h(R)}\bigg( 1 \wedge \frac{(R-|y|)^2}{t} \bigg)^{-d/2} t^{-d/2} \exp \bigg( - \frac{(R-|y|)^2}{4t}\bigg).
	\end{align*}
\end{lem}

Recall that $R_1'$ and $\eta_0$ are the constants defined in \eqref{e:R1p} and \eqref{e:def-eta-0}, respectively.

\begin{lem}\label{l:upper-basic-3}
	There exist $c_1,c_2>0$ such that for all $t\in (0,1]$, $R\in (0,R_1']$ and $x,y \in \R^d_0$ satisfying $|x|<\eta_0(t/2)^{1/(2+\beta)}$ and  $2^{1/\beta}|x|<  |y| < R$, it holds that
	\begin{align*}
		&\limsup_{\delta \to0} \frac{\P_x (X_t \in B(y,\delta),\, M_t< R)}{|B(0,\delta)|} \nn\\
		& \le  \frac{c_1 h(|x|)}{|y|^dh(\eta_0 (t/2)^{1/(2+\beta)})}  \exp \Big( - \frac{(1-2^{-1/\beta})\sqrt{\kappa} \, |y|}{\sqrt{2} R^{1+\beta}}\Big) + c_2 h(|x|) \exp \Big(  - \frac{3\kappa t}{8R^{2+2\beta}} +\frac{5\sqrt \kappa}{2\beta |y|^\beta}\Big).
	\end{align*}
\end{lem}
\pf  Using 	the strong Markov property in the equality below, 
Lemmas  \ref{l:survival-probability}(i), \ref{l:generator}(i-a),  \ref{l:survival-probability-t-2}, and the comparability
of $h$ and $\tilde{h}$ in the second inequality, we get that 
that for any sufficiently small $\delta>0$, 
\begin{align}\label{e:upper-basic-3-0}
	&\P_x\left(X_t \in B(y,\delta),\, M_t < R \right)\nn\\
	&= \E_x\left[ \P_{X_{\tau_{B(0,2^{-1/\beta}|y|)}}}\Big(X^{B(0,R)}_{t-\tau_{B(0,2^{-1/\beta}|y|)}} \in B(y,\delta) \Big) ;  \tau_{B(0,2^{-1/\beta}|y|)}< \zeta \wedge t \right] \nn\\
	&\le \P_x\left(t/2\le  \tau_{B(0,2^{-1/\beta}|y|)} < \zeta \right) \sup_{0<s\le t/2, \, |z| = 2^{-1/\beta}|y|} \P_z \left(X_{s}^{B(0,R)} \in B(y,\delta)\right)\nn\\
	&\quad +  \P_x\left( \tau_{B(0,2^{-1/\beta}|y|)}  < \zeta \wedge (t/2) \right) \sup_{t/2\le s\le t, \, |z| =2^{-1/\beta}|y|} \P_z \left(X^{B(0,R)}_{s} \in B(y,\delta) \right)\nn\\
		&\le  \frac{c_1 h(|x|)}{h(\eta_0 (t/2)^{1/(2+\beta)})}  \sup_{0<s\le t/2, \, |z| = 2^{-1/\beta}|y|} \P_z \left(X_{s}^{B(0,R)} \in B(y,\delta)\right)\nn\\
	&\quad + \frac{c_2h(|x|)}{h(2^{-1/\beta}|y|)}  \sup_{t/2\le s\le t, \, |z| =2^{-1/\beta}|y|} \P_z \left(X^{B(0,R)}_{s} \in B(y,\delta) \right).
\end{align}
Using \eqref{e:Feynman-Kac-2} and \eqref{e:V-cond-1} in the first inequality below, $R\le (\kappa/(4C_1))^{1/(2\beta-\beta')}
$ in the second, and \eqref{e:exp-poly} in the third, we see that
\begin{align}
	& \limsup_{\delta\to 0}\sup_{0<s\le t/2, \, |z| =2^{-1/\beta}|y|} \frac{\P_z \big(X^{B(0,R)}_{s} \in B(y,\delta)\big) }{|B(0,\delta)|} \nn\\
	&\le 	\limsup_{\delta\to 0} \sup_{0<s\le t/2, \, |z| = 2^{-1/\beta}|y|} \frac{\P_z \left(W_{s} \in B(y,\delta) \right) }{|B(0,\delta)|}\exp \Big(- s   \inf_{r \in (0,R)} \Big(\frac{\kappa}{r^{2+2\beta}}- \frac{C_1}{r^{2+\beta'}} \Big)  \Big) \nn\\
	&\le  \sup_{0<s\le t/2}\,(4\pi s)^{-d/2}\exp \Big( - \frac{((1-2^{-1/\beta})|y|)^2}{4s}- \frac{3\kappa s}{4R^{2+2\beta}}\Big)\nn\\
	&\le  \frac{c_2}{|y|^d}\sup_{s>0} \exp \Big( - \frac{((1-2^{-1/\beta})|y|)^2}{6s}- \frac{3\kappa s}{4R^{2+2\beta}}\Big)\label{e:upper-basic-3-1'}\\
	&=  \frac{c_2}{|y|^d}  \exp \Big( - \frac{(1-2^{-1/\beta})\sqrt{\kappa} \,|y|}{\sqrt{2}R^{1+\beta}}\Big).\label{e:upper-basic-3-1}
\end{align}
Following the arguments for \eqref{e:upper-basic-3-1'}, we also have
\begin{align*}
	&\limsup_{\delta \to 0}  \sup_{t/2\le s\le t, \, |z| = 2^{-1/\beta}|y|} \frac{\P_z \big(X^{B(0,R)}_{s} \in B(y,\delta) \big)}{|B(0,\delta)|}\nn\\
	&\le  \frac{c_2}{|y|^d}\sup_{t/2\le s\le t} \exp \Big( - \frac{((1-2^{-1/\beta})|y|)^2}{6s}- \frac{3\kappa s}{4R^{2+2\beta}}\Big)\le \frac{c_2}{|y|^d}\exp \Big(- \frac{3\kappa t}{8R^{2+2\beta}}\Big).
\end{align*}
Hence, using  \eqref{e:exp-poly}, we get
\begin{align}\label{e:upper-basic-3-2}
	&\limsup_{\delta \to 0}  \sup_{t/2\le s\le t, \, |z| = 2^{-1/\beta}|y|} \frac{\P_z \big(X^{B(0,R)}_{s} \in B(y,\delta) \big)}{|B(0,\delta)|}\nn\\
	&\le c_3 \begin{cases}
		\displaystyle  \frac{h(2^{-1/\beta}|y|)}{|y|^{(d+2+\beta)/2}}\exp \Big(  - \frac{3\kappa t}{8R^{2+2\beta}} +\frac{2\sqrt \kappa}{\beta |y|^\beta}\Big)  &\mbox{ if $|y|< 2^{1/\beta}$},\\
		\displaystyle \frac{ h(2^{-1/\beta}|y|)}{|y|^d}\exp \Big(  - \frac{3\kappa t}{8R^{2+2\beta}} \Big) &\mbox{ if $|y|\ge 2^{1/\beta}$}
	\end{cases}\nn\\
	&\le c_3h(2^{-1/\beta}|y|)\exp \Big(  - \frac{3\kappa t}{8R^{2+2\beta}} +\frac{5\sqrt \kappa}{2\beta |y|^\beta}\Big) .
\end{align}
Combining \eqref{e:upper-basic-3-0}, \eqref{e:upper-basic-3-1} and \eqref{e:upper-basic-3-2},  we arrive at  the result. \qed

\begin{lem}\label{l:upper-basic-5}
	There exist $c_1,c_2>0$ such that for all $t\in (0,1]$, $0<r<R\le R_1'$ and $x,y \in \R^d_0$ satisfying $|x|<\eta_0(t/2)^{1/(2+\beta)}$ and  $|x| \vee |y|<r$, it holds that
	\begin{align*}
	&	\limsup_{\delta \to0} \frac{\P_x (X_t \in B(y,\delta),\, r \le M_t< R)}{|B(0,\delta)|}\nn\\ &\le\frac{c_1h(|x|)}{(r-|y|)^dh( \eta_0 (t/2)^{1/(2+\beta)} )}   \exp \bigg( - \frac{\sqrt \kappa( r-|y|)}{\sqrt 8 R^{1+\beta}}\bigg)  + \frac{c_2  h(|x|)}{(r-|y|)^d} \exp \bigg( - \frac{3\kappa t}{8R^{2+2\beta}} + \frac{5\sqrt \kappa}{4\beta r^{\beta}}\bigg).
	\end{align*}
\end{lem}
\pf Using 	the strong Markov property, we see that for any sufficiently small $\delta>0$,
\begin{align}\label{e:upper-basic-5-1}
	&\P_x (X_t \in B(y,\delta),\, r \le M_t< R)\nn\\
	&= \E_x\left[ \P_{X_{\tau_{B(0,r)}}}\Big(X^{B(0,R)}_{t-\tau_{B(0, r)}} \in B(y,\delta) \Big) ; \tau_{B(0, r)} < \zeta \wedge t \right] \nn\\
	&\le \P_x\left(t/2\le  \tau_{B(0,r)} < \zeta \right) \sup_{0<s\le t/2, \, |z| =r} \P_z \left(X_{s}^{B(0,R)} \in B(y,\delta)\right)\nn\\
	&\quad +  \P_x\left( \tau_{B(0,r)} < \zeta \wedge (t/2) \right) \sup_{t/2\le s\le t, \, |z| = r} \P_z \left(X^{B(0,R)}_{s} \in B(y,\delta) \right).
\end{align}
Following the arguments for \eqref{e:upper-basic-3-1'},  using \eqref{e:Feynman-Kac-0}, \eqref{e:V-cond-1} and  \eqref{e:exp-poly},  we get
\begin{align}\label{e:upper-basic-5-2}
	& \limsup_{\delta\to 0}\sup_{0<s\le t/2, \, |z| = r} \frac{\P_z \big(X_{s}^{B(0,R)} \in B(y,\delta)\big)}{|B(0,\delta)|} \le  \sup_{0<s\le t/2}\frac{1}{(4\pi s)^{d/2}}\exp \Big( - \frac{( r-|y|)^2}{4s}- \frac{3\kappa s}{4R^{2+2\beta}}\Big)\nn\\
		&\le \frac{c_1}{(r-|y|)^d} \sup_{s>0}\,\exp \Big( - \frac{( r-|y|)^2}{6s}- \frac{3\kappa s}{4R^{2+2\beta}}\Big)= \frac{c_1}{(r-|y|)^d} \exp \Big( - \frac{\sqrt \kappa( r-|y|)}{\sqrt 8 R^{1+\beta}}\Big)
\end{align}
and
\begin{align}\label{e:upper-basic-5-3}
	&\limsup_{\delta \to 0}  \sup_{t/2\le s\le t, \, |z| = r} \frac{\P_z \big(X^{B(0,R)}_{s} \in B(y,\delta) \big)}{|B(0,\delta)|}\nn\\
	&\le \frac{c_1}{(r-|y|)^d} \sup_{t/2\le s\le t}\,\exp \Big( - \frac{( r-|y|)^2}{6s}- \frac{3\kappa s}{4R^{2+2\beta}}\Big)\le \frac{c_1}{(r-|y|)^d} \exp \Big( - \frac{3\kappa t}{ 8 R^{2+2\beta}}\Big)\nn\\
	&\le c_2 \begin{cases}
		\displaystyle \frac{ h(r)r^{(d-2-\beta)/2}}{(r-|y|)^d}\exp \Big(  - \frac{3\kappa t}{ 8 R^{2+2\beta}} +\frac{\sqrt \kappa}{\beta r^\beta}\Big)  &\mbox{ if $r< 1$},\\
		\displaystyle \frac{ h(r)}{(r-|y|)^d}\exp \Big(  - \frac{3\kappa t}{ 8R^{2+2\beta}} \Big) &\mbox{ if $r\ge 1$}
	\end{cases}\nn\\
	&\le  \frac{c_3h(r)}{(r-|y|)^d}\exp \Big(  - \frac{3\kappa t}{ 8R^{2+2\beta}} +\frac{5\sqrt \kappa}{4\beta r^\beta}\Big).
\end{align}
Combining \eqref{e:upper-basic-5-1}, \eqref{e:upper-basic-5-2} and  \eqref{e:upper-basic-5-3}, and using Lemmas  \ref{l:survival-probability}(i), \ref{l:generator}(i-a) and \ref{l:survival-probability-t-2},  we arrive at the result.
 \qed

We first use Lemma \ref{l:upper-basic-1} to establish interior upper bounds.

\begin{lem}\label{l:upper-HKE-0}
	There exist constants $c_1,c_2>0$ such that for all $t \in (0,1]$ and $x,y \in \R^d_0$ satisfying $|x| \wedge |y| \ge \eta_0(t/2)^{1/(2+\beta)}$, it holds that
	\begin{align}\label{e:upper-HKE-0}
		p(t,x,y) &\le c_1 t^{-d/2}\exp \left(- \frac{|x-y|^2}{16t}-\frac{c_2t}{(|x|\vee |y|)^{2+2\beta}}\right).
	\end{align}
\end{lem} 
\pf  By symmetry, we can assume that $|x|\le |y|$. If $|y| \ge R_1't^{1/(2+2\beta)}/2$, then  \eqref{e:upper-HKE-0} follows from $p(t,x,y) \le q(t,x,y)$. Suppose $|y|<R_1't^{1/(2+2\beta)}/2$. By the lower-semicontinuity of $p$, we have
\begin{align*}
	&p(t,x,y) \le \liminf_{\delta\to0} 
	\frac{1}{|B(y,\delta)|}	\int_{B(y,\delta)}p(t,x,v)dv\nn\\
	&\le \limsup_{\delta \to0} \frac{\P_x (X_t \in B(y,\delta),\, M_t< 2|y|)}{|B(0,\delta)|} +  \limsup_{\delta \to0} \frac{\P_x(X_t \in B(y,\delta),\,M_t \ge 2|y|)}{|B(0,\delta)|}\\
	&=:I_1+I_2.
\end{align*}
 By Lemma \ref{l:survival-probability-t-1}, we  have 
\begin{align*}
	I_1&\le  (4\pi t)^{-d/2}  \exp \Big( - \frac{|x-y|^2}{4t}  - \frac{3\kappa t }{4(2|y|)^{2+2\beta}}\Big).
\end{align*}
Further, using Lemma \ref{l:upper-basic-1} and the fact that $R_1/2>|y| > \eta_0 (t/2)^{1/2}$  in the first inequality below,   Lemma \ref{l:h-ratio} and the fact $|x|\le |y|$ in the second,  and $|x-y|\le 2|y|$ and $|y|^{2+\beta}>\eta_0^{2+\beta}t/2$ in the last, we obtain
\begin{align*}
	I_2&\le \frac{c_1h(|x|)}{ h(2|y|)}   t^{-d/2} \exp \Big( - \frac{|y|^2}{4t}\Big) \le c_2t^{-d/2} \exp \Big(  - \frac{|y|^2}{4t}- \frac{(2^{\beta/2}-1)\sqrt \kappa}{2^{3\beta/4}\beta |y|^\beta}\Big)\nn\\
			&\le c_2t^{-d/2} \exp \Big( - \frac{|x-y|^2}{16t}- \frac{(2^{\beta/2}-1)\sqrt \kappa \eta_0^{2+\beta} t}{2^{1+3\beta/4}\beta |y|^{2+2\beta}} \Big).
\end{align*}
 The proof is complete. \qed

\begin{lem}\label{l:upper-HKE-1}
	There exist constants $c_1,c_2>0$ such that for all $t \in (0,1]$ and $x,y \in \R^d_0$ satisfying $|x| \wedge |y| < \eta_0(t/2)^{1/(2+\beta)}$ and $|x| \vee |y| \ge \eta_0 t^{1/(2+\beta)}$, it holds that
	\begin{align*}
		p(t,x,y) &\le  \frac{c_1h(|x| \wedge |y|)}{h( \eta_0 (t/2)^{1/(2+\beta)})} t^{-d/2}\exp \left(-\frac{c_2|x-y|^2}{t}\right).
	\end{align*}
\end{lem} 
\pf  
Without loss of generality, we assume that $|x|\le |y|$. 
Set $R:=\eta_0(t/2)^{1/(2+\beta)}  \wedge R_1$. Note that $|y|-R \ge (1-2^{-1/(2+\beta)}) |y|$ and $|y-x| \le (1+2^{-1/(2+\beta)})|y|$. Thus, it holds that 
\begin{align}\label{e:upper-HKE-1-0}
	|y| - R \ge  (1-2^{-1/(2+\beta)})  \left( (\eta_0 t^{1/(2+\beta)}) \vee  \frac{|x-y|}{ 1+ 2^{-1/(2+\beta)}}\right).
\end{align} 
Using the lower-semicontinuity of $p$ and \eqref{e:basic-upper-0},  we have
\begin{align}\label{e:upper-HKE-1-1}
\!\!\!\!\!\!p(t,x,y)	\le \limsup_{\delta \to0} \frac{\P_x (X_t \in B(y,\delta) )}{|B(0,\delta)|}\le \frac{c_1 h(|x|)}{ h(R)}  \sup_{0<s\le t} (4\pi s)^{-d/2} \exp \Big(\!\! - \frac{(|y|-R)^2}{4s}\Big).
\end{align}
Further, using \eqref{e:basic-upper-1}, \eqref{e:basic-upper-2} and \eqref{e:upper-HKE-1-0}, we see that
\begin{align}\label{e:upper-HKE-1-2}
	&	\sup_{0<s\le t} s^{-d/2} \exp \Big( - \frac{(|y|-R)^2}{4s}\Big) \le  t^{-d/2} \Big( 1 \wedge \frac{(|y|-R)^2}{2dt}\Big)^{-d/2} \exp \Big( - \frac{(|y|-R)^2}{4t}\Big)\nn\\
		&\le  t^{-d/2} \Big( 1 \wedge \frac{c_2}{t^{\beta/(2+\beta)}}\Big)^{-d/2} \exp \Big( - \frac{c_3|x-y|^2}{t}\Big) \le c_4 t^{-d/2} \exp \Big( - \frac{c_3|x-y|^2}{t}\Big).
\end{align}
Combining \eqref{e:upper-HKE-1-1} with \eqref{e:upper-HKE-1-2}
and using $h( \eta_0 (t/2)^{1/(2+\beta)}) \asymp h(R)$ for $t\in (0,1]$,
we arrive at the result.\qed

		\begin{lem}\label{l:upper-HKE-2}
		There exist constants $c_1,c_2>0$ such that for all $t \in (0,1]$ and $x,y \in \R^d_0$ satisfying $|x| \wedge |y| < \eta_0(t/2)^{1/(2+\beta)}$ and $|x| \vee |y| < \eta_0 t^{1/(2+\beta)}$, it holds that
		\begin{align*}
		p(t,x,y) &\le  \frac{c_1h(|x| \wedge |y|)}{h( \eta_0 (t/2)^{1/(2+\beta)})} t^{-d/2}\exp \left(-c_2 t^{-\beta/(2+\beta)}\right).
	\end{align*}
	\end{lem} 
	\pf  
 Without loss of generality, we assume that $|x|\le |y|$. 
 If
  $t > 2^{-2}(R_1'/\eta_0)^{2+\beta}$, then by the semigroup property, $p(t/2,z,y)\le q(t/2,z,y)$ and Lemma \ref{l:survival-probability-t-2}, we get
 \begin{align*}
 	p(t,x,y) \le  \int_{\R^d_0} p(t/2,x,z) q(t/2,z,y) dz \le (2\pi t)^{-d/2} \P_x\left(\zeta >  t/2 \right) \le c_{1}t^{-d/2}h(|x|).
 \end{align*}
 Since  $\exp (- t^{-\beta/(2+\beta)})\ge c_2$ and $h(\eta_0(t/2)^{1/(2+\beta)}) \ge c_{3} $ by the almost monotonicity of $h$,
 the desired result follows immediately from the display above.

Assume now that 
$t \le 2^{-2}(R_1'/\eta_0)^{2+\beta}$.   Let 	$
 a:=2^{-1/(2+\beta)}
 $ and 
	$n_0\ge 1$ be  such that 
	\begin{align*}
		 a^{n_0}\eta_0 t^{1/(2+\beta)}\le  |y| < a^{n_0-1}\eta_0 t^{1/(2+\beta)}.
		\end{align*}
Note that		$|y|/a^{n_0+1} \le R_1'$.	We have
	\begin{align*}
		& p(t,x,y) \le \liminf_{\delta\to 0}
		\frac{1}{|B(y,\delta)|}\int_{B(y,\delta)}	p(t,x,v) dv\\ &\le  \limsup_{\delta\to 0}\frac{\P_x\left(X_t \in B(y,\delta),\, M_t \ge a^{-n_0-1}|y| \right) }{|B(0,\delta)|} +  \limsup_{\delta\to 0} \frac{\P_x\left( X_t \in B(y,\delta),\, M_t < a^{-1}|y| \right)}{|B(0,\delta)|}\\
		&\quad\;\;    + \sum_{1\le m \le n_0}\limsup_{\delta\to 0} \frac{\P_x\left(X_t \in B(y,\delta),\, a^{-m}|y| \le  M_t < a^{-m-1}|y| \right) }{|B(0,\delta)|}=:I_1+I_2+I_3.
	\end{align*}

	For $I_1$, applying Lemma \ref{l:upper-basic-1} and using the almost monotonicity of 
	$h$ and the assumptions $|y|<  \eta_0 t^{1/(2+\beta)}$ and $t\le 1$,   we get
	\begin{align}\label{e:small-upper1-I1}
		I_1 & \le\limsup_{\delta\to 0}\frac{\P_x\left(X_t \in B(y,\delta),\,  M_t \ge a^{-1}\eta_0 t^{1/(2+\beta)} \right) }{|B(0,\delta)|}  \nn\\
		&\le  \frac{c_3h(|x|)}{h(a^{-1}\eta_0 t^{1/(2+\beta)})} \Big( 1 \wedge \frac{(a^{-1}\eta_0 t^{1/(2+\beta)}-|y|)^2}{t} \Big)^{-d/2}  t^{-d/2} \exp \Big( - \frac{(a^{-1}\eta_0 t^{1/(2+\beta)}-|y|)^2}{4t}\Big)\nn\\
			&\le \frac{c_3h(|x|)}{h(a\eta_0 t^{1/(2+\beta)})}   \Big( 1 \wedge \frac{(2^{1/(2+\beta)}-1)^2\eta_0^2}{t^{\beta/(2+\beta)}} \Big)^{-d/2}  t^{-d/2} \exp \Big( - \frac{(2^{1/(2+\beta)}-1)^2\eta_0 ^2}{4t^{\beta/(2+\beta)}}\Big)\nn\\
			&\le \frac{c_4h(|x|)}{h(a\eta_0 t^{1/(2+\beta)})}  t^{-d/2} \exp \Big( - \frac{(2^{1/(2+\beta)}-1)^2\eta_0 ^2}{4t^{\beta/(2+\beta)}}\Big).
	\end{align}

	For $I_2$, first assume $|x| \ge 2^{-1/\beta}|y|$. Using 
	Lemma \ref{l:survival-probability-t-1} in the first inequality below,  $t>(|y|/\eta_0)^{2+\beta}$ in the second, $|y|\le 2^{1/\beta}|x| \wedge \eta_0t^{1/(2+\beta)}$ in the third  and Lemma \ref{l:h-ratio}  in the  last, we get
  	\begin{align}\label{e:small-upper1-I2-1}
  I_2&\le (4\pi t)^{-d/2}  \exp \Big( - \frac{|x-y|^2}{4t}  - \frac{3\kappa t }{4(|y|/a)^{2+2\beta}}\Big)\le  (4\pi t)^{-d/2} \exp \Big(- \frac{6\sqrt\kappa }{\beta |y|^\beta}\Big) \nn\\
  &\le  (4\pi t)^{-d/2}  \exp \Big(  - \frac{3\sqrt\kappa }{2\beta|x|^\beta} - \frac{3\sqrt\kappa }{\beta (\eta_0 t^{1/(2+\beta)})^\beta} \Big)\nn\\
  &\le  \frac{c_5h(|x|)}{ h( a \eta_0 t^{1/(2+\beta)})} t^{-d/2}  \exp \Big(  -\frac{\sqrt\kappa }{\beta (\eta_0 t^{1/(2+\beta)})^\beta}\Big) .
  \end{align}   
  In the case $|x| < 2^{-1/\beta}|y|$, 
  using Lemma \ref{l:upper-basic-3} in the first inequality below,  $t>(|y|/\eta_0)^{2+\beta}$ in the second, \eqref{e:exp-poly} and $t\le 1$ in the third, and  $|y|\le \eta_0 t^{1/(2+\beta)}$  in the fourth, we obtain
  	\begin{align}\label{e:small-upper1-I2-2}
  I_2	& \le  \frac{c_6 h(|x|)}{|y|^d h(a\eta_0 t^{1/(2+\beta)})}  \exp \Big( - \frac{(1-2^{-1/\beta})\sqrt{\kappa} |y|}{\sqrt {2}(|y|/a)^{1+\beta}}\Big)  +c_7h(|x|)\exp \Big(  - \frac{3\kappa t}{8(|y|/a)^{2+2\beta}} +\frac{5\sqrt \kappa}{2\beta |y|^\beta}\Big)\nn\\
  	& \le  \frac{c_6 h(|x|)}{|y|^dh(a\eta_0 t^{1/(2+\beta)})}  \exp \Big( - \frac{(1-2^{-1/\beta})\sqrt{\kappa}}{2^{1/2+(1+\beta)/(2+\beta)}|y|^{\beta}}\Big)  +c_7h(|x|)\exp \Big(  - \frac{\sqrt \kappa }{2\beta |y|^{\beta}} \Big)\nn\\
  	& \le  \frac{c_8h(|x|)}{h(a\eta_0 t^{1/(2+\beta)})}  t^{-d/2} \exp \Big( - \frac{c_9}{|y|^{\beta}}\Big) \le  \frac{c_8 h(|x|)}{h(a\eta_0 t^{1/(2+\beta)})}  t^{-d/2} \exp \Big( - \frac{c_9}{( \eta_0 t^{1/(2+\beta)})^\beta}\Big) .
  \end{align}

For $I_3$,	let $1\le  m\le n_0$.  Since $m\le n_0$, we have
\begin{align}\label{e:small-upper1-I3-ingre}
	 \frac{\kappa t}{(|y|/a^{m+1})^{2+2\beta}}  \ge  \frac{\kappa a^{-(m-1)(2+\beta)} (|y|/\eta_0)^{2+\beta}}{(|y|/a^{m+1})^{2+2\beta}} =    \frac{4\sqrt \kappa}{\beta (|y|/a^m)^{\beta}}.
\end{align}
Using Lemma \ref{l:upper-basic-5} in the first inequality below, \eqref{e:small-upper1-I3-ingre} in the second,   \eqref{e:exp-poly} in the third, and $t\le 1$ and $|y|/a^m<\eta_0t^{1/(2+\beta)}/a$ in the last,  we get
\begin{align*}
	&	\limsup_{\delta \to0} \frac{\P_x (X_t \in B(y,\delta),\, a^{-m}|y|\le M_t< a^{-m-1}|y|)}{|B(0,\delta)|}\nn\\ &\le\frac{c_{10}h(|x|)}{(a^{-m}-1)^d |y|^dh( a\eta_0 t^{1/(2+\beta)} )}   \exp \Big( - \frac{\sqrt \kappa( a^{-m}-1)|y|}{\sqrt 8(|y|/a^{m+1})^{1+\beta}}\Big) \nn\\
	&\quad + \frac{c_{11} h(|x|)}{(a^{-m}-1)^d |y|^d} \exp \Big( - \frac{3\kappa t}{8(|y|/a^{m+1})^{2+2\beta}} + \frac{5\sqrt \kappa}{4\beta (|y|/a^m)^{\beta}}\Big)\nn\\
	&\le\frac{c_{10}h(|x|)}{(1-a^m)^d\, (|y|/a^m)^d\,h( a\eta_0 t^{1/(2+\beta)} )}   \exp \Big( - \frac{a^{1+\beta}( 1-a)\sqrt \kappa}{\sqrt 8(|y|/a^m)^{\beta}}\Big) \nn\\
	&\quad + \frac{c_{11} h(|x|)}{(1-a^m)^d\, (|y|/a^m)^d} \exp \Big( - \frac{\sqrt \kappa}{4\beta(|y|/a^m)^{\beta}}\Big)\nn\\
		&\le\frac{c_{12}h(|x|) |y|/a^m}{(1-a)^d h( a\eta_0 t^{1/(2+\beta)} )}   \exp \Big( - \frac{a^{1+\beta}( 1-a)\sqrt \kappa}{2\sqrt 8(|y|/a^m)^{\beta}}\Big) + \frac{c_{13} h(|x|) |y|/a^m}{(1-a)^d} \exp \Big( - \frac{\sqrt \kappa}{8\beta(|y|/a^m)^{\beta}}\Big)\nn\\
&\le\frac{c_{14}h(|x|) |y|/a^m}{h( a\eta_0 t^{1/(2+\beta)} )}   t^{-d/2}\exp \Big( - \frac{c_{15}}{t^{\beta/(2+\beta)}}\Big).
\end{align*}
Thus, since $\sum_{1\le m\le n_0} a^{-m} \le a^{-n_0}/(1-a)$ and $|y|/a^{n_0}< \eta_0 t^{1/(2+\beta)}/a \le \eta_0/a$, we deduce that
\begin{equation}\label{e:small-upper1-I3}
 I_3
 \le   \frac{c_{14}h(|x|) |y|/a^{n_0}}{h( a\eta_0 t^{1/(2+\beta)} )}   t^{-d/2}\exp \Big(\! - \frac{c_{15}}{ t^{\beta/(2+\beta)}}\Big)\! \le\! \frac{c_{14}\eta_0 h(|x|) }{ah( a\eta_0 t^{1/(2+\beta)} )}   t^{-d/2}\exp \Big(\! - \frac{c_{15}}{ t^{\beta/(2+\beta)}}\Big).
\end{equation}

	Combining \eqref{e:small-upper1-I1}, \eqref{e:small-upper1-I2-1}, \eqref{e:small-upper1-I2-2} and \eqref{e:small-upper1-I3}, we get the  result. The proof is complete.	\qed

	\begin{lem}\label{l:upper-HKE-half}
		There exist constants $c_1,c_2,c_3>0$ such that for all $t \in (0,1]$ and $x,y \in \R^d_0$,
		\begin{align*}
			p(t,x,y) &\le c_1 \bigg( 1 \wedge  \frac{h(|x| \wedge |y|)}{h( \eta_0 (t/2)^{1/(2+\beta)})}\bigg)  t^{-d/2}\exp \left(- \frac{c_2|x-y|^2}{t} - \frac{c_3t}{(|x| \vee |y| \vee t^{1/(2+\beta)})^{2+2\beta}}\right).
		\end{align*}
	\end{lem} 
		\pf  
	Without loss of generality, we assume that $|x|\le |y|$. We consider the following three cases separately.
	
	\smallskip
	
	\noindent
{\it
	Case 1: $|x| \ge \eta_0 (t/2)^{1/(2+\beta)}$.} The result follows from Lemma \ref{l:upper-HKE-0}. 
	
	\noindent
{\it
	Case 2:  $|x|<\eta_0 (t/2)^{1/(2+\beta)}$ and $|y|\ge \eta_0 t^{1/(2+\beta)}$.} Observe that
	\begin{align*}
		\frac{|x-y|^2}{t} \ge \frac{(1-2^{-1/(2+\beta)})^2|y|^2}{t} \ge \frac{c_1}{t^{\beta/(2+\beta)}} \ge \frac{c_2t}{(|y| \vee t^{1/(2+\beta)})^{2+2\beta}}.
	\end{align*}
	Hence, applying Lemma \ref{l:upper-HKE-1}, we obtain
	\begin{align*}
		p(t,x,y) &\le  \frac{c_3h(|x| )}{h( \eta_0 (t/2)^{1/(2+\beta)})} t^{-d/2}\exp \Big(-\frac{c_4|x-y|^2}{t}\Big)\\
		& \le  \frac{c_3h(|x| )}{h( \eta_0 (t/2)^{1/(2+\beta)})} t^{-d/2}\exp \Big(-\frac{c_4|x-y|^2}{2t}-\frac{c_2c_4t}{2(|y| \vee t^{1/(2+\beta)})^{2+2\beta}}\Big).
	\end{align*}
	
	\noindent
{\it
	Case 3: $|x|<\eta_0 (t/2)^{1/(2+\beta)}$ and $|y|< \eta_0 t^{1/(2+\beta)}$.} In this case, it holds that
	\begin{align}\label{e:upper-HKE-small-small}
		\frac{|x-y|^2}{t} \le \frac{4|y|^2}{t}\le \frac{c_5t}{(|y| \vee t^{1/(2+\beta)})^{2+2\beta}} \le \frac{c_6}{t^{\beta/(2+\beta)}} .
	\end{align}
	Thus, using Lemma \ref{l:upper-HKE-2}, we get
	\begin{align*}
	p(t,x,y) &\le  \frac{c_7h(|x|)}{h( \eta_0 (t/2)^{1/(2+\beta)})} t^{-d/2}\exp \left(-c_8 t^{-\beta/(2+\beta)}\right)\\
	 &\le  \frac{c_7h(|x|)}{h( \eta_0 (t/2)^{1/(2+\beta)})} t^{-d/2}\exp \Big(-\frac{c_8|x-y|^2}{2c_6t}-\frac{c_8}{2t^{\beta/(2+\beta)}}\Big).
\end{align*}
The proof is complete. \qed

Recall that $\eta_1= 2^{-2/(2+\beta)}\eta_0$.

\begin{prop}\label{p:smalltime}
 There exist   $c_1,c_2,c_3>0$ such that \eqref{e:smalltime-upper} holds for all $t \in (0,1]$ and $x,y \in \R^d_0$.
\end{prop}
\pf
	Without loss of generality, we assume that $|x|\le |y|$. If $|y| \ge \eta_1 t^{1/(2+\beta)}$, then,  since $h$ is almost increasing,  the result follows from Lemma \ref{l:upper-HKE-half}. Assume $|y| < \eta_1t^{1/(2+\beta)}$. Using the semigroup property, we have
\begin{align*}
	p(t,x,y) = \bigg( \int_{B(0, 2\eta_1t^{1/(2+\beta)})} + \int_{B(0, 2\eta_1t^{1/(2+\beta)})^c}\bigg) p(t/2,x,z) p(t/2,z,y)dz=:I_1+I_2.
\end{align*}
Applying Lemma \ref{l:upper-HKE-half} and using the almost monotonicity of $h$, we get
\begin{align*}
	&I_1\le c_1 \int_{B(0, 2\eta_1t^{1/(2+\beta)})}  \Big( 1 \wedge  \frac{h(|x| \wedge |z|)}{h( \eta_1 t^{1/(2+\beta)})}\Big) q(c_2t,x,z ) \exp \Big(-\frac{c_3}{t^{\beta/(2+\beta)}}\Big)\nn\\
	&\qquad \qquad \qquad\qquad \qquad  \times  \Big( 1 \wedge  \frac{h(|y| \wedge |z|)}{h(\eta_1 t^{1/(2+\beta)})}\Big) q(c_2t,z,y ) \exp \Big(-\frac{c_3}{t^{\beta/(2+\beta)}}\Big)dz\nn\\
	&\le c_4 \Big( 1 \wedge  \frac{h(|x| )}{h(\eta_1t^{1/(2+\beta)})}\Big)\Big( 1 \wedge  \frac{h(|y|)}{h(\eta_1 t^{1/(2+\beta)})}\Big) \exp \Big(-\frac{2c_3}{t^{\beta/(2+\beta)}}\Big) \nn\\
	&\quad \times \int_{B(0, 2\eta_1t^{1/(2+\beta)})} q(c_2t,x,z ) q(c_2t,z,y )dz\nn\\
		&\le c_4 \Big( 1 \wedge  \frac{h(|x| )}{h( \eta_1t^{1/(2+\beta)})}\Big)\Big( 1 \wedge  \frac{h(|y|)}{h( \eta_1 t^{1/(2+\beta)})}\Big) (8\pi c_2t)^{-d/2} \exp \Big(- \frac{|x-y|^2}{8c_2t}-\frac{2c_3}{t^{\beta/(2+\beta)}}\Big).
\end{align*} 
Note that \eqref{e:upper-HKE-small-small} remains valid with different constants $c_5$ and $c_6$. 
Applying Lemma \ref{l:upper-HKE-half}, and using  the almost monotonicity of $h$ and \eqref{e:upper-HKE-small-small}, we also get
\begin{align*}
	I_2&\le c_1 \int_{B(0, 2\eta_1t^{1/(2+\beta)})^c}  \Big( 1 \wedge  \frac{h(|x| \wedge |z|)}{h( \eta_1t^{1/(2+\beta)})}\Big)    \Big( 1 \wedge  \frac{h(|y| \wedge |z|)}{h( \eta_1t^{1/(2+\beta)})}\Big)q(c_2t,x,z ) q(c_2t,z,y ) dz\nn\\
	&\le c_5 \Big( 1 \wedge  \frac{h(|x| )}{h(\eta_1t^{1/(2+\beta)})}\Big)\Big( 1 \wedge  \frac{h(|y|)}{h(\eta_1 t^{1/(2+\beta)})}\Big) (4\pi c_2t)^{-d/2} \exp \Big( - \frac{(\eta_1t^{1/(2+\beta)})^2}{4c_2t}\Big)\nn\\
	&\quad \times \int_{B(0, 2\eta_1t^{1/(2+\beta)})^c}  q(c_2t,x,z) dz\nn\\
	&\le c_6 \Big( 1 \wedge  \frac{h(|x| )}{h( \eta_1t^{1/(2+\beta)})}\Big)\Big( 1 \wedge  \frac{h(|y|)}{h( \eta_1 t^{1/(2+\beta)})}\Big) t^{-d/2} \exp \left(-\frac{c_7}{t^{\beta/(2+\beta)}}\right) \nn\\
		&\le c_6 \Big( 1 \wedge  \frac{h(|x| )}{h( \eta_1 t^{1/(2+\beta)})}\Big)\Big( 1 \wedge  \frac{h(|y|)}{h( \eta_1 t^{1/(2+\beta)})}\Big) t^{-d/2} \exp \Big(- \frac{c_8|x-y|^2}{t}-\frac{c_7}{2t^{\beta/(2+\beta)}}\Big) .
\end{align*} 
\qed

Using the semigroup property,  as a consequence of  Proposition \ref{p:smalltime}, we get the following upper bound of $p(t,x,y)$ for all $t>0$, which is sharp when  $d \ge 3$.

  \begin{cor}\label{c:UHK-general}
There exist  $c_1,c_2,c_3>0$  such that for all $t >0$ and $x,y \in \R^d_0$,
		\begin{align}\label{e:UHK-general}
		 p(t,x,y) &\le c_1\bigg( 1 \wedge  \frac{h_{\beta,\kappa}(|x| )}{h_{\beta,\kappa}( \eta_1 t^{1/(2+\beta)} \wedge 1)}\bigg)\bigg( 1 \wedge  \frac{h_{\beta,\kappa}(|y|)}{h_{\beta,\kappa}( \eta_1 t^{1/(2+\beta)} \wedge 1)}\bigg)  \nn\\
		&\quad \times \exp \bigg(- \frac{c_2t}{(|x| \vee |y| \vee t^{1/(2+\beta)})^{2+2\beta}}\bigg)q(c_3t,x,y).
		\end{align}
	\end{cor} 
	\pf Let $t>0$ and $x,y \in \R^d_0$. If $t\le 1$, then   \eqref{e:UHK-general} follows from Proposition \ref{p:smalltime}. Suppose that $t>1$.
	Using  the semigroup property and  \eqref{e:UHK-general} with $t=1/3$, we obtain
	\begin{align*}
		&p(t,x,y)= \int_{\R^d_0}\int_{\R^d_0} p(1/3,x,v) p(t-2/3,v,w) p(1/3, w,y) dvdw\\
		&\le c_1(1 \wedge h_{\beta,\kappa}(|x|))(1 \wedge h_{\beta,\kappa}(|y|))  \int_{\R^d_0}\int_{\R^d_0}   q(c_2,x,v)  q(t-2/3,v,w)  q(c_2, w,y)  dvdw\\
		&=  c_1(1 \wedge h_{\beta,\kappa}(|x|))(1 \wedge h_{\beta,\kappa}(|y|))   \big(4\pi(t-2/3+2c_2)\big)^{-d/2} \exp \Big( - \frac{|x-y|^2}{4(t-2/3+2c_2)} \Big) \nn\\
		&\le c_2(1 \wedge h_{\beta,\kappa}(|x|))(1 \wedge h_{\beta,\kappa}(|y|))  t^{-d/2} \exp \Big( - \frac{|x-y|^2}{12t} \Big) .
	\end{align*}Since
	\begin{align*}
		\exp \Big(- \frac{t}{(|x| \vee |y| \vee t^{1/(2+\beta)})^{2+2\beta}}\Big) \ge  \exp \left( - t^{-\beta/(2+\beta)}\right) \ge e^{-1},
	\end{align*}the proof is complete.\qed 

\noindent \textbf{Proof of  Theorem \ref{t:smalltime}(i).} 
It follows immediately from Corollary \ref{c:UHK-general}.
 
 \subsection{Small time lower estimates}\label{ss:smalltime-lower}

 Throughout this subsection, we assume that 
 $V\in \Kle$  with $\beta' \in (0,\beta)$ and $C_2>0$ being the constants in \eqref{e:V-cond-2}
 and continue to write $h$ instead of $h_{\beta,\kappa}$. 
 
 We begin with the following observation for the potential $V$.

 \begin{lem}\label{l:V-upper-bound} 
 	There exists $C>0$ such that for all $r>0$, it holds that
\begin{align*}
	\sup_{x\in  B(0,r/2)^c} V(x) \le \frac{C}{(r\wedge 1)^{2+2\beta}}.
\end{align*} 	
 \end{lem}
 \pf If $r\ge 2$, then the result follows from \eqref{e:V-bounded}. If  $r<2$, then by \eqref{e:V-cond-2} and  \eqref{e:V-bounded}, we get
 \begin{align*}
 		\sup_{x\in B(0,r/2)^c} V(x)& \le 	\sup_{R\in [r/2,1]} \frac{\kappa +C_2R^{2\beta-\beta'}}{R^{2+2\beta}} + 	\sup_{x\in B(0,1)^c} V(x)\le \frac{\kappa + C_2}{(r/2)^{2+2\beta}} + c_1\le \frac{\kappa+C_2+c_1}{(r/2)^{2+2\beta}} .
 \end{align*}
 \qed

 In the Feynman-Kac semigroup $(P_t)$, 
 two factors contribute to the exponential part 
 in the  lower estimate \eqref{e:smalltime-lower}: one comes from the killing potential $V$ and the other from the exponential decay in the Gaussian heat kernel  $q(t,x,y)$. We will use the following lemma to find the dominant term in the exponent.

 Recall that $\eta_0$ is defined as \eqref{e:def-eta-0} and $\eta_1= 2^{-2/(2+\beta)}\eta_0$.
 
 \begin{lem}\label{l:lower-tech-1}
(i) There exists $C\ge 1$ such that for all $t \in (0,1]$ and $x,y \in \R^d_0$ satisfying $|y|\ge |x| \ge \eta_1t^{1/(2+\beta)}$, it holds that
 	\begin{align}\label{e:lower-tech-1}
 	\frac{t}{(|x| \wedge 1)^{2+2\beta}} \le  	\frac{C|x-y|^2}{t}+\frac{ 2t}{|y|^{2+2\beta}} + 1.
 	\end{align}

 	\noindent (ii) For all $t \in (0,1]$ and $x,y \in \R^d_0$ satisfying $|x| < \eta_1t^{1/(2+\beta)}\le |y|$, it holds that \begin{align*}		\frac{|y|^2}{t} \le  \frac{4|x-y|^2}{t}+\frac{2^{4+2\beta}\eta_1^{4+2\beta}t}{|y|^{2+2\beta}}.	\end{align*}
 \end{lem}
 \pf (i) If $|x|\ge 1$ or $2^{-1/(2+2\beta)}|y| \le |x|<1$, then \eqref{e:lower-tech-1} holds trivially. Suppose $|x|<1$ and 
  $|y| > 2^{1/(2+2\beta)}|x|$. Then  we get
 \begin{align*}
 	\frac{|x-y|^2}{t} \ge \frac{(2^{1/(2+2\beta)}-1)^2|x|^2}{t} \ge \frac{(2^{1/(2+2\beta)}-1)^2\eta_1^{4+2\beta}t^{2}}{t|x|^{2+2\beta}} = \frac{c_1t}{|x|^{2+2\beta}}.
 \end{align*}

\noindent  (ii) Set $r_t:=\eta_1t^{1/(2+\beta)}$.  	If $ |y|\le 2r_t$, then\begin{align*} 	\frac{|y|^2}{t} \le \frac{ 4r_t^2}{t} = \frac{4\eta_1^{4+2\beta}t}{r_t^{2+2\beta}} \le  \frac{2^{4+2\beta}\eta_1^{4+2\beta}t}{|y|^{2+2\beta}},\end{align*} and if $|y|>2r_t$, then $|y-x|\ge  |y|-r_t \ge |y|/2$ so that $|y|^2/t \le 4|x-y|^2/t$.  
  \qed

We start with interior lower bounds. 
 
 \begin{lem}\label{l:lower-HKE-0}
 (i) Suppose $d\ge 2$.	There exist constants $c_1,c_2,c_3>0$ such that for all $t \in (0,1]$ and $x,y \in \R^d_0$ satisfying $|y| \ge |x| \ge \eta_1t^{1/(2+\beta)}$, it holds that
 	\begin{align}\label{e:lower-HKE-0}
 	 p(t,x,y) &\ge c_1 t^{-d/2}\exp \left(- \frac{c_2|x-y|^2}{t}-\frac{c_3 t}{|y|^{2+2\beta}}\right).
 	\end{align}
 	(ii) Suppose $d= 1$.	Then \eqref{e:lower-HKE-0} holds for all $t \in (0,1]$ and $x,y \in \R^d_0$ satisfying $xy>0$ and $|y| \ge |x| \ge \eta_1t^{1/(2+\beta)}$.
 \end{lem} 
 \pf    
 We only give the proof for $d\ge 3$. The proof for $d\le 2$ is similar.
 Using  \eqref{e:Feynman-Kac-2}, Lemma \ref{l:V-upper-bound} and  Proposition \ref{p:DHKE}, since  $|y| \ge |x| \ge  \eta_1t^{1/2}$, we obtain
 \begin{align*}
 	p(t,x,y)&\ge  p^{\overline B(0,|x|/2)^c}(t,x,y) \ge q^{\overline B(0,|x|/2)^c}(t,x,y) \exp \Big(-t \sup_{z\in  B(0,|x|/2)^c}V(z)\Big)\nn\\
 	&\ge c_1 \Big( 1 \wedge \frac{|x|/2}{\sqrt t}\Big) 
	\Big( 1 \wedge \frac{|x|/2}{\sqrt t}\Big) 
	t^{-d/2}\exp \Big(- \frac{c_2|x-y|^2}{t}-\frac{c_3 t}{(|x| \wedge 1)^{2+2\beta}}\Big)\nn\\
 		&\ge c_4 t^{-d/2}\exp \Big(- \frac{c_2|x-y|^2}{t}-\frac{c_3 t}{(|x| \wedge 1)^{2+2\beta}}\Big).
 \end{align*}
 By Lemma \ref{l:lower-tech-1}(i), this implies the desired result. \qed

 Note that, using the strong Markov property and the joint continuity of $p^{\overline B(0,\eps)^c}$, 
 we get that for all $r>\eps>0$, $t>0$, $s \in [t/2,t]$ and $x,y \in  \overline B(0,\eps)^c$, 
 \begin{align}\label{e:lower-2-10}
 	p^{\overline B(0,\eps)^c}(s,x,y)&\ge \E_x\left[ p^{\overline B(0,\eps)^c}(s-\tau_{B(0,r)\setminus \overline B(0,\eps)}, X^{\overline B(0,\eps)^c}_{\tau_{B(0,r)\setminus \overline B(0,\eps)}}, y); \tau_{B(0,r)\setminus \overline B(0,\eps)} <\zeta \wedge (t/3) \right]\nn\\
 	&\ge  \P_x \left(\tau_{B(0,r)\setminus \overline B(0,\eps)} <\zeta \wedge (t/3) \right)\inf_{t/6 \le s\le t,\,|z|=r} p^{\overline B(0,\eps)^c} (s, z,y)  .
 \end{align}

 \begin{lem}\label{l:lower-HKE-1}
 (i) Suppose $d\ge 2$.	There exist constants $c_1,c_2,c_3>0$ such that for all $t \in (0,1]$, $s \in [t/2,t]$ and $x,y \in \R^d_0$ satisfying $|x| < \eta_1t^{1/(2+\beta)}\le |y|$, it holds that
 	\begin{align}\label{e:lower-HKE-1}
 		p(s,x,y) &\ge  \frac{c_1h(|x|)}{h( \eta_1t^{1/(2+\beta)})} t^{-d/2}\exp \left(- \frac{c_2|x-y|^2}{t}-\frac{c_3 t}{|y|^{2+2\beta}}\right).
 	\end{align}
 	(ii) Suppose $d=1$. Then \eqref{e:lower-HKE-1} holds for all $t \in (0,1]$, $s \in [t/2,t]$ and $x,y \in \R^d_0$ satisfying $xy>0$ and $|x| < \eta_1t^{1/(2+\beta)}\le |y|$.
 \end{lem} 
 \pf We only give the proof of (i). The proof of (ii) is similar.
   Define $r_u:= \eta_1 u^{1/(2+\beta)}$ for $u>0$. Let  $t_0 \in (0,1]$ be the constant in Lemma \ref{l:survival-probability-t-3}. We consider the following three cases separately.
 
 \smallskip
 
 \noindent
 {\it Case 1: $t\le t_0$.}
By \eqref{e:lower-2-10}, for any sufficiently small $\eps>0$, 
  \begin{align}\label{e:lower-2-1}
 	p^{\overline B(0,\eps)^c}(s,x,y)\ge  \P_x \left(\tau_{B(0,r_t)\setminus \overline B(0,\eps)} <\zeta \wedge (t/3) \right)\inf_{t/6 \le s\le t,\,|z|=r_t} p^{\overline B(0,\eps)^c} (s, z,y)  .
 \end{align}
 Applying Lemma \ref{l:survival-probability-t-3}, we get
 \begin{align}\label{e:lower-2-2}
 \lim_{\eps\to0}\P_x \left(\tau_{B(0,r_t)\setminus \overline B(0,\eps)} <\zeta \wedge (t/3) \right) =\P_x \left(\tau_{B(0,r_t)} <\zeta \wedge (t/3) \right) \ge \frac{c_1h(|x|)}{h(r_t)}.
\end{align} 
 Further, by Lemmas \ref{l:lower-HKE-0} and \ref{l:lower-tech-1}(ii), since $|z-y|\le 2|y|$ for all $|z|=r_t$, we have
 \begin{align}\label{e:lower-2-3}
 	&\lim_{\eps \to 0} \inf_{t/6 \le s\le t,\,|z|=r_t} p^{\overline B(0,\eps)^c} (s, z,y)=  \inf_{t/6 \le s\le t,\,|z|=r_t} p (s, z,y)  \nn\\
 		& \ge c_2 t^{-d/2}\exp \Big(- \frac{c_3(2|y|)^2}{t}-\frac{c_4 t}{|y|^{2+2\beta}}\Big)\ge c_2 t^{-d/2}\exp \Big(- \frac{c_5|x-y|^2}{t}-\frac{c_6 t}{|y|^{2+2\beta}}\Big).
 \end{align}
  Combining   \eqref{e:lower-2-1}--\eqref{e:lower-2-3}, 
    we arrive at the result in this case.

 \noindent
 {\it Case 2: $t>t_0$ and $|x|<r_{t_0/2}$.} By Lemmas \ref{l:lower-HKE-0} and \ref{l:lower-tech-1}(ii), we see that for all $z \in B(0,2r_{t})\setminus B(0,r_{t})$,
 \begin{align}\label{e:lower-2-4}
 	&p(s-t_0/4, z,y) \ge c_7 (t/2-t_0/4)^{-d/2}\exp \Big(- \frac{c_8|z-y|^2}{t/2-t_0/4}-\frac{c_9 (t-t_0/4)}{|y|^{2+2\beta}}\Big)\nn\\
 	&\ge c_{10}t^{-d/2}\exp \Big(- \frac{c_{8}(3|y|)^2}{t/4}-\frac{c_9 t}{|y|^{2+2\beta}}\Big)\ge c_{10}\exp \Big(- \frac{c_{11}|x-y|^2}{t}-\frac{c_{12} t}{|y|^{2+2\beta}}\Big).
 \end{align}
 Besides, applying the result from {\it Case 1}, we get that for all $z \in B(0,2r_{t})\setminus B(0,r_{t})$,
 \begin{align}\label{e:lower-2-5}
 &	p(t_0/4,x,z)\ge  \frac{c_{13}h(|x|)}{h(r_{t_0/2})} t_0^{-d/2}\exp \Big(- \frac{c_{14}(3r_t)^2}{t_0}-\frac{c_{15} t_0}{r_t^{2+2\beta}}\Big)\nn\\
 	&\ge  \frac{c_{13}h(|x|)}{h(r_{t_0/2})} t_0^{-d/2}\exp \Big(- \frac{c_{14}(3r_1)^2}{t_0}-\frac{c_{15} t_0}{r_{t_0}^{2+2\beta}}\Big) \ge c_{16} h(|x|).
 \end{align}
 Using the semigroup property, \eqref{e:lower-2-4},   \eqref{e:lower-2-5}, $t\le 1$ and the almost monotonicity of $h$,  we deduce that
 	\begin{align}\label{e:lower-2-6}
 	p(s,x,y)& \ge\int_{B(0,2r_{t_0})\setminus B(0,r_{t_0})} p(t_0/4,x,z)p(s-t_0/4, z,y)dz\nn\\
 	&\ge c_{17} h(|x|)\exp \Big(- \frac{c_{11}|x-y|^2}{t}-\frac{c_{12} t}{|y|^{2+2\beta}}\Big)
	|B(0,2r_{t_0})\setminus B(0,r_{t_0})|
	\nn\\
 		&= \frac{c_{18}t_0^{d/2} h(r_{t_0})h(|x|)}{h(r_{t_0})} \exp \Big(- \frac{c_{11}|x-y|^2}{t}-\frac{c_{12} t}{|y|^{2+2\beta}}\Big)\nn\\
 			&\ge \frac{c_{19}t_0^{d/2} h(r_{t_0})h(|x|)}{h(r_{t})} t^{-d/2}\exp \Big(- \frac{c_{11}|x-y|^2}{t}-\frac{c_{12} t}{|y|^{2+2\beta}}\Big).
 \end{align}
 
 \noindent
{\it  Case 3:  $t>t_0$ and $|x|\ge r_{t_0/2}$.} Using  Lemma \ref{l:lower-HKE-0}, we obtain for all $z\in B(0,2r_{t})\setminus B(0,r_{t})$,
 \begin{equation}\label{e:lower-2-7-0}
 	p(t_0/4,x,z)\ge  c_{20} t_0^{-d/2}\exp \Big(- \frac{c_{21}(3r_t)^2}{t_0}-\frac{c_{22} t_0}{r_{t}^{2+2\beta}}\Big)\ge  c_{23} t_0^{-d/2}\exp \Big(- \frac{c_{21}(3r_1)^2}{t_0}-\frac{c_{22} t_0}{r_{t_0}^{2+2\beta}}\Big).
 \end{equation}
 Note that   $h(|x|) \le c_{24}h(r_{1})$ 
 by the almost monotonicity of $h$. Hence, by \eqref{e:lower-2-7-0}, we obtain
 \begin{align}\label{e:lower-2-7}
 		p(t_0/4,x,z)&\ge  c_{25}  \ge c_{26}h(|x|) \quad \text{for all  $z\in B(0,2r_{t})\setminus B(0,r_{t})$}.
 \end{align}
 Repeating the arguments for \eqref{e:lower-2-6}, using \eqref{e:lower-2-7} instead of \eqref{e:lower-2-5}, we get the  result.  \qed

\begin{lem}\label{l:lower-HKE-2}
(i) Suppose $d\ge 2$.	There exist constants $c_1,c_2>0$ such that for all $t \in (0,1]$ and $x,y \in \R^d_0$ satisfying $|x| \le |y|< \eta_1t^{1/(2+\beta)}$, it holds that
	\begin{align}\label{e:lower-HKE-2}
		p(t,x,y) &\ge  \frac{c_1h(|x|)h(|y|)}{h( \eta_1 t^{1/(2+\beta)})^2} t^{-d/2}\exp \left(-c_2t^{-\beta/(2+\beta)}\right).
	\end{align}
(ii) Suppose $d=1$. Then \eqref{e:lower-HKE-2} holds for all $t \in (0,1]$ and $x,y \in \R^d_0$ satisfying $xy>0$ and  $|x| \le |y|< \eta_1t^{1/(2+\beta)}$.
\end{lem} 
\pf Since the proofs are similar, we only present the proof for (i). Set $r_t:=\eta_1t^{1/(2+\beta)}$. By the semigroup property and  Lemma \ref{l:lower-HKE-1}, we get
\begin{align*}
&p(t,x,y)\ge 	\int_{B(0,2r_t)\setminus B(0,r_t)} p(t/2,x,z)  p(t/2,z,y) dz \nn\\
&\ge \frac{c_1^2h(|x|)h(|y|)}{h( r_t)^2} t^{-d}\exp \left(- \frac{2c_2 t}{r_t^{2+2\beta}}\right) 	\int_{B(0,2r_t)\setminus B(0,r_t)} \exp \left(- \frac{c_3(|x-z|^2 + |z-y|^2)}{t}\right) dz\\
&\ge \frac{c_1^2h(|x|)h(|y|)}{h( r_t)^2} t^{-d}\exp \left(- \frac{2c_2 t}{r_t^{2+2\beta}}- \frac{2c_3(3r_t)^2}{t}\right) 	\int_{B(0,2r_t)\setminus B(0,r_t)}  dz\\
&\ge \frac{c_1^2h(|x|)h(|y|)}{h( r_t)^2} t^{-(1+\beta)d /(2+\beta)}\exp \left(-c_4t^{-\beta/(2+\beta)}\right) \ge \frac{c_1^2h(|x|)h(|y|)}{h( r_t)^2} t^{-d/2}\exp \left(-c_4t^{-\beta/(2+\beta)}\right) .
\end{align*}
 \qed
 
 \smallskip

 Combining Lemmas \ref{l:lower-HKE-0}, \ref{l:lower-HKE-1} and \ref{l:lower-HKE-2}, we  obtain the full lower bounds for  $t \in (0,1]$.
  
 	\begin{prop}\label{p:smalltimel}	
		(a) If $d\ge 2$, then there exist  $c_4,c_5,c_6>0$  such that \eqref{e:smalltime-lower} holds  for all $t \in (0,1]$ and $x,y \in \R^d_0$.
		
			\noindent	
		(b) If  $d=1$, there  exist  $c_4,c_5,c_6>0$ such that \eqref{e:smalltime-lower} holds for all $t \in (0,1]$ and $x,y \in \R^1_0$ with $xy>0$.
	\end{prop} 
\pf
	  The proposition follows easily from Lemmas \ref{l:lower-HKE-0}, \ref{l:lower-HKE-1} and \ref{l:lower-HKE-2}.  \qed

  \noindent \textbf{Proof of  Theorem \ref{t:smalltime}(ii).} 
	 Since the proofs are similar, we only give the proof for $d\ge 2$. 
	 By Proposition \ref{p:smalltimel}, it suffices to prove \eqref{e:smalltime-lower} for 
$t\in (1,T]$.  Let 
 $x,y \in \R^d_0$,  $t \in (1,T]$ 
and $4\le N\in \N$ be such that $N-1< 3t \le N$.
Set 
	$
	z_x:=(1 + 2|x|^{-1}) x $ and  $	z_y:=(1 + 2|y|^{-1}) y .
	$
For all $v \in B(z_x,1)$ and $w \in B(z_y,1)$,  we have $|x-v| \vee |y-w| \le 3$, 
$|v| \wedge |w| \ge 1$ and $|v-w| \le |x-y|+6$. Thus, by Proposition \ref{p:smalltimel} and the almost monotonicity of $h_{\beta,\kappa}$, we get that for all $s \in [1/3,2/3]$, $v,v' \in B(z_x,1)$ and $w\in B(z_y,1)$,
\begin{align*}
	p(s,v,v') &\ge  c_1,\\
	p(s,x,v) &\ge  c_2\Big(1 \wedge \frac{h_{\beta,\kappa}(|x|)}{h_{\beta,\kappa}(\eta_1 s^{1/(2+\beta)})} \Big) \ge  c_3 \Big(1 \wedge \frac{h_{\beta,\kappa}(|x|)}{h_{\beta,\kappa}(\eta_1t^{1/(2+\beta)})} \Big),\\
	p(s,w,y) &\ge c_2\Big(1 \wedge \frac{h_{\beta,\kappa}(|y|)}{h_{\beta,\kappa}(\eta_1 s^{1/(2+\beta)})} \Big) \ge  c_3 \Big(1 \wedge \frac{h_{\beta,\kappa}(|y|)}{h_{\beta,\kappa}(\eta_1t^{1/(2+\beta)})} \Big),
\end{align*}
and
\begin{align*}
		p(s,v,w) &\ge c_4 s^{-d/2} \exp \Big( - \frac{c_5(|x-y|+6)^2}{s} \Big) \ge c_6 s^{-d/2} \exp \Big( - \frac{2c_5(|x-y|^2 + 36) }{s} \Big)\nn\\	& \ge c_7 t^{-d/2} \exp \Big( - \frac{2Nc_5|x-y|^2 }{t} \Big).
\end{align*}
Using these and  the semigroup property, we arrive at
	\begin{align*}
	&	p(t,x,y)\ge \int_{B(z_x,1)} \cdots \int_{B(z_x,1)}\int_{B(z_y,1)} p(1/3,x,v_1) p(1/3,v_1,v_2) \cdots p(1/3,v_{N-3},v_{N-2}) \\
	&\qquad \qquad \qquad \qquad\qquad \qquad\qquad  \times p(t-(N-1)/3,v_{N-2},w) p(1/3, w,y) \,dv_1 \cdots dv_{N-2} dw\\
		&\ge   c_8|B(0,1)|^{N-1} \Big(1 \wedge \frac{h_{\beta,\kappa}(|x|)}{h_{\beta,\kappa}(\eta_1t^{1/(2+\beta)})} \Big)  \Big(1 \wedge \frac{h_{\beta,\kappa}(|x|)}{h_{\beta,\kappa}(\eta_1t^{1/(2+\beta)})} \Big)  t^{-d/2}  \exp \Big(- \frac{2Nc_5|x-y|^2}{t}\Big) .
	\end{align*}
	Since 	$4 \le N<3T+1$,
the proof is complete. \qed

\section{Optimality of the 
	class $\Kloc$}\label{s:counters}

  In this section,  we consider the following two conditions:
 
 \medskip
 
 (1) There exists $C>0$ such that
 \begin{align}\label{e:V-cond-1-critical}
 	V(x) \le \frac{\kappa}{|x|^{2+2\beta}} -  \frac{C}{|x|^{2+\beta}}  \quad \text{for all} \;\, x \in \R^d_0 \text{ with } |x|\le 1.
 \end{align}
 
  (2) There exists $C>0$ such that
 \begin{align}\label{e:V-cond-2-critical}
 	V(x) \ge \frac{\kappa}{|x|^{2+2\beta}} + \frac{C}{|x|^{2+\beta}}  \quad \text{for all} \;\, x \in \R^d_0 \text{ with } |x|\le 1.
 \end{align}
 We will show that our conditions on $V$ in Theorem \ref{t:smalltime} are optimal by showing that when \eqref{e:V-cond-1-critical} holds, the assertion of Theorem \ref{t:smalltime}(i) is false, and that when
 \eqref{e:V-cond-2-critical} holds, the assertion of Theorem \ref{t:smalltime}(ii) is false.

Recall that $L^V=\Delta - V$ and that the function $\wt h_{\beta,\kappa}$ defined in \eqref{e:function-wth}.
\begin{lem}\label{l:generator-counterexample}	
	\noindent(i)  
	If \eqref{e:V-cond-1-critical} holds, then there exists a constant $a_1>0$  such that 
	\begin{align*}
		L^V (|\cdot|^{-a_1} \wt h_{\beta,\kappa}(|\cdot|))(x)\ge 0 \quad \text{for all $x \in \R^d_0$ with $|x|\le 1$}.
	\end{align*}

	\noindent (ii)  
	If  \eqref{e:V-cond-2-critical} holds, then there 
	exists a constant $a_2>0$  such that 
	\begin{align*}
	L^V (|\cdot|^{a_2} \wt h_{\beta,\kappa}(|\cdot|))(x)\le 0 \quad \text{for all $x \in \R^d_0$ with $|x|\le 1$}.
	\end{align*}
\end{lem}
\pf  Since the proofs are similar, we only give the proof of (i). 
 For all $a>0$ and $x \in \R^d_0$ with $|x|\le 1$, using \eqref{e:V-cond-1-critical} in the first inequality below, and Lemma \ref{l:wt-h-kappa-2} in the 
 second, we get
\begin{align*}
	L^V (|\cdot|^{-a} \wt h_{\beta,\kappa}(|\cdot|))(x) &\ge   L^{\beta,\kappa}( |\cdot|^{-a} \wt h_{\beta,\kappa}(|\cdot|))(x) + \frac{c_1\wt h_{\beta,\kappa}(|x|) }{|x|^{2+\beta+a}}\nn\\
	&\ge  \left( c_1-c_2a  - a (d-a-2) |x|^{\beta}  \right) \frac{\wt h_{\beta,\kappa}(|x|)}{|x|^{2+\beta+a}}.
\end{align*}
Taking $a_1>0$ to satisfy $|c_2+d-a_1-2|a_1 =c_1$, we arrive at the  result.\qed

Repeating the arguments of Lemma \ref{l:survival-probability}, using Lemma \ref{l:generator-counterexample} in place of Lemma \ref{l:generator}, we get
\begin{lem}\label{l:survival-probability-counterexample}
	(i) 
	If \eqref{e:V-cond-1-critical} holds, then
	\begin{align*}
		\P_x \left(\tau_{B(0,R)} <\zeta \right) \ge \frac{|x|^{-a_1} \wt h_{\beta,\kappa}(|x|)}{ R^{-a_1} \wt h_{\beta,\kappa}(R)}  \quad \text{ for all  $x\in \R^d_0$ with $|x|<R\le 1$,}
	\end{align*}
	where $a_1>0$ is the constant in Lemma \ref{l:generator-counterexample}(i).
	
	\noindent	(ii) 
	If \eqref{e:V-cond-2-critical} holds,  then
	\begin{align*}
		\P_x \left(\tau_{B(0,R)} <\zeta \right) \le \frac{|x|^{a_2} \wt h_{\beta,\kappa}(|x|)}{ R^{a_2} \wt h_{\beta,\kappa}(R)}  \quad \text{ for all $x\in \R^d_0$ with $|x|<R\le 1$,}
	\end{align*}
	where $a_2>0$ is the constant in Lemma \ref{l:generator-counterexample}(ii).
\end{lem}

Recall that the constant $\eta_0>0$ is defined as \eqref{e:def-eta-0}. Applying Lemma \ref{l:survival-probability-counterexample}(ii)  in place of Lemma \ref{l:survival-probability}(i) in the proof of 
Lemma \ref{l:survival-probability-t-2} (with $R_1=1$), 
We get the following upper bounds of the survival probability.

\begin{lem}\label{l:survival-probability-t-counterexample}
	If \eqref{e:V-cond-2-critical} holds, then there exists $C>0$ such that
	\begin{equation*}
		\P_x \left(\zeta >  t \right) \le \frac{C|x|^{a_2}h_{\beta,\kappa}(|x|)}{t^{a_2/(2+\beta)}h_{\beta,\kappa}( \eta_0 t^{1/(2+\beta)} )} \quad \text{for all $t\in (0,1]$ and $x\in \R^d_0$ with $|x|<\eta_0 t^{1/(2+\beta)}$},
	\end{equation*}
	where   $a_2>0$ is the constant in Lemma \ref{l:generator-counterexample}(ii).
\end{lem}

We now check what kind of 
survival probability estimates one can get from the heat kernel estimates in Theorem \ref{t:smalltime}.

\begin{lem}\label{l:consequence-HKE}
	(i) 
	If \eqref{e:smalltime-upper} 	holds with $T=1$,
	then there exists $C>0$ such that
		\begin{equation*}
		\P_x \left(\zeta >  t \right) \le \frac{Ch_{\beta,\kappa}(|x|)}{h_{\beta,\kappa}( \eta_1 t^{1/(2+\beta)} )} \quad \text{for all $t\in (0,1]$ and $x\in \R^d_0$ with $|x|<\eta_1 t^{1/(2+\beta)}$}.
	\end{equation*}
		(ii) 
		If the estimates in Theorem \ref{t:smalltime}(ii) 
			hold with $T=1$, 
		 then
		there exist $C,\lambda>0$ such that
	\begin{equation*}
		\P_x \left(\zeta >  t \right) \ge \frac{Ce^{-\lambda t^{-\beta/(2+\beta)}}h_{\beta,\kappa}(|x|)}{h_{\beta,\kappa}( \eta_1 t^{1/(2+\beta)} )}  \quad \text{for all $t\in (0,1]$ and $x\in \R^d_0$ with $|x|<\eta_1 t^{1/(2+\beta)}$}.
	\end{equation*}
\end{lem}
\pf (i) Using \eqref{e:smalltime-upper}, we get that for all $t\in (0,1]$ and $x\in \R^d_0$ with $|x|<\eta_1 t^{1/(2+\beta)}$,
\begin{align*}
		\P_x \left(\zeta >  t \right) = \int_{\R^d_0} p(t,x,y) dy \le  \frac{c_1h_{\beta,\kappa}(|x|)}{h_{\beta,\kappa}( \eta_1 t^{1/(2+\beta)} )} \int_{\R^d} q(c_2t,x,y) dy = \frac{c_1h_{\beta,\kappa}(|x|)}{h_{\beta,\kappa}( \eta_1 t^{1/(2+\beta)} )} .
\end{align*}
(ii) Since the proofs are similar, we only give the proof for $d\ge 2$. Using \eqref{e:smalltime-lower} and the almost monotonicity of $h$, we get that for all $t\in (0,1]$ and $x\in \R^d_0$ with $|x|<\eta_1 t^{1/(2+\beta)}$,
\begin{align*}
	\P_x \left(\zeta >  t \right) &\ge \int_{B(0,2\eta_1 t^{1/(2+\beta)}) \setminus B(0,\eta_1 t^{1/(2+\beta)})} p(t,x,y) dy \\
	&\ge  \frac{c_1h_{\beta,\kappa}(|x|)}{h_{\beta,\kappa}( \eta_1 t^{1/(2+\beta)} )}  t^{-d/2} \exp \bigg( - \frac{c_2(3\eta_1 t^{1/(2+\beta)})^2}{t}\bigg) \int_{B(0,2\eta_1 t^{1/(2+\beta)}) \setminus B(0,\eta_1 t^{1/(2+\beta)})}   dy\\
		&\ge  \frac{c_3h_{\beta,\kappa}(|x|)}{h_{\beta,\kappa}( \eta_1 t^{1/(2+\beta)} )}   \exp \bigg( - \frac{c_2(3\eta_1 t^{1/(2+\beta)})^2}{t}\bigg) .
\end{align*}
\qed

\begin{lem}\label{l:survival-probability-t-3-counterexample}
Assume \eqref{e:V-cond-1-critical}. If	  \eqref{e:smalltime-upper}  
	holds with $T=1$, 
	then there exist $t_0'\in (0,1]$ and $C>0$ such that for all $t\in (0,t_0']$ and $x\in \R^d_0$ with $|x|<\eta_1(t/4)^{1/(2+\beta)}$, it holds that
	\begin{equation*}
		\P_x \left(\tau_{B(0,\eta_1(t/4)^{1/(2+\beta)} )} < \zeta \wedge (t/3) \right) \ge \frac{C|x|^{-a_1}h_{\beta,\kappa}(|x|)}{ t^{-a_1/(2+\beta)}h_{\beta,\kappa}( \eta_1 (t/4)^{1/(2+\beta)} )},
	\end{equation*}
	where $a_1>0$ is the constant in Lemma \ref{l:generator-counterexample}(i).
\end{lem}
\pf  Using Lemmas \ref{l:survival-probability-counterexample}(i) and \ref{l:consequence-HKE}(i), we get that for all $t\in (0,4\eta_1^{-2-\beta} \wedge 1 ]$ and $x \in \R^d_0$ with $|x|<\eta_1(t/4)^{1/(2+\beta)}$,
\begin{align}\label{e:survival-probability-t-3-counterexample}
	&\P_x \left(\tau_{B(0,\eta_1(t/4)^{1/(2+\beta)} )} < \zeta \wedge (t/3) \right) \ge \P_x \left(\tau_{B(0,\eta_1(t/4)^{1/(2+\beta)} )} < \zeta \right) - \P_x \left( \zeta>t/3 \right)\nn\\
	&\ge \frac{ c_1|x|^{-a_1}h_{\beta,\kappa}(|x|)}{ (\eta_1 t^{1/(2+\beta)})^{-a_1} h_{\beta,\kappa}( \eta_1 (t/4)^{1/(2+\beta)} )} - \frac{c_2h_{\beta,\kappa}(|x|)}{h_{\beta,\kappa}( \eta_1 (t/3)^{1/(2+\beta)} )}.
\end{align}
Note that since $h_{\beta,\kappa}$ is exponentially increasing at $0$, there exists $t_0' \in (0,4\eta_1^{-2-\beta} \wedge 1 ]$ such that    $h_{\beta,\kappa}( \eta_1 (t/3)^{1/(2+\beta)} )\ge \frac{2c_2}{c_1}h_{\beta,\kappa}( \eta_0 (t/4)^{1/(2+\beta)} )$ for all $t\in (0,t_0']$. 
Thus, the desired result follows from \eqref{e:survival-probability-t-3-counterexample}.
\qed

\begin{thm}\label{t:counter-smalltime}
	(i)		
	If \eqref{e:V-cond-1-critical} holds, then \eqref{e:smalltime-upper} fails for any $c_1,c_2,c_3>0$.
	
	\noindent (ii) 		
	If \eqref{e:V-cond-2-critical} holds, then the estimates in Theorem \ref{t:smalltime}(ii) fail.
\end{thm}	
\pf (i) Suppose that \eqref{e:smalltime-upper}  
holds (with $T=1$). 
We assume $d\ge 2$. 
The case $d=1$ is similar.

By \eqref{e:lower-2-10} and Lemma \ref{l:survival-probability-t-3-counterexample},
we get that for all $t \in (0,t_0']$ and $x,y \in \R^d_0$ satisfying $|x|<\eta_1(t/4)^{1/(2+\beta)}\le |y|$,
\begin{align*}
	p(t,x,y) & \ge \lim_{\eps \to 0} \P_x \left(\tau_{B(0,\eta_1 (t/4)^{1/(2+\beta)})\setminus B(0,\eps)} <\zeta \wedge (t/3) \right)\inf_{t/6 \le s\le t,\,|z|=\eta_1 (t/4)^{1/(2+\beta)}} p^{\overline B(0,\eps)^c} (s, z,y) \nn\\
	&\ge  \frac{c_1|x|^{-a_1}h_{\beta,\kappa}(|x|)}{ t^{-a_1/(2+\beta)}h_{\beta,\kappa}(\eta_1 (t/4)^{1/(2+\beta)})} t^{-d/2}\exp \left(- \frac{c_2|x-y|^2}{t}-\frac{c_3 t}{|y|^{2+2\beta}}\right).
\end{align*}
Thus, for all $t \in (0,t_0']$ and  $|y| \ge \eta_1(t/4)^{1/(2+\beta)}$, we get $\liminf_{|x|\to 0} p(t,x,y) |x|^{a_1}/h_{\beta,\kappa}(|x|)>0$.
This contradicts \eqref{e:smalltime-upper}.

\noindent (ii) If the assertions in  Theorem \ref{t:smalltime}(ii) 
hold (with $T=1$),
 then by Lemma \ref{l:consequence-HKE}(ii), we get that for all $t \in (0,1]$, $\liminf_{|x|\to 0} \P_x(\zeta>t) /h_{\beta,\kappa}(|x|) >0$. This contradicts Lemma \ref{l:survival-probability-t-counterexample}. \qed

\section{Large time  heat kernel estimates}\label{s:large}

Recall that the function $H_{d, \beta, \kappa}$ is defined in \eqref{e:H},  that the function $\psi_{d,R}$ is defined in \eqref{e:def-psi-R}, and that $\Lg r=\log(e-1+r)$. 
Observe that
\begin{align}
	\log r \le \Lg r \le 2\log r \quad &\text{ for all} \;\, r \ge 2,\label{e:Log-prop-0}\\
		\Lg sr \le \Lg s + \Lg r \quad &\text{ for all} \;\, r,s\ge 0,\label{e:Log-prop-1}\\
		\Lg r \le 
		(r/a)\Lg a
		 \le r \quad &\text{ for all} \;\, r \ge a \ge 1.\label{e:Log-prop-2}
\end{align}
Further, for any $\eps>0$, there exists $C=C(\eps)\ge 1$ such that
\begin{align}\label{e:Log-scaling}
{\Lg r}/{\Lg s} \le C ({r}/{s})^{\eps} \quad \text{for all} \;\, r\ge s>0.
\end{align}
Note also that when $d=1$ or $d=2$, for any  $R>0$,
there exist comparison constants depending on $R$  such that for all $t>0$ and  $r>0$,
\begin{align}\label{e:H-arbitrary-R}
		H_{d,\beta,\kappa}(t,r)\asymp  \begin{cases}
		\displaystyle   1 \wedge \frac{h_{\beta,\kappa}(r)\1_{(0,R)}(r) + 			( \Lg r ) 
			\1_{[R, \infty)}(r) }{\Lg \sqrt t} &\mbox{ if $d=2$}, \\[12pt]
		\displaystyle  1 \wedge \frac{h_{\beta,\kappa}(r) \1_{(0,R)}(r) + r \1_{[R, \infty)}(r)}{\sqrt t} &\mbox{ if $d= 1$}.
	\end{cases}
\end{align}
Using \eqref{e:H-arbitrary-R} (if $d=1$ or $d=2$), we see that for any  $R>0$, there exists $C=C(R)>1$ such that
\begin{align}\label{e:compare-H-psi}
 C^{-1} \psi_{d,R}(t,r)\le 	H_{d,\beta,\kappa}(t,r) \le  C \psi_{d,R}(t,r) \quad \text{for all $t\ge R^2$ and $r \ge 2R$}.
\end{align}

The goal of this section is to establish the following large time estimates for $p(t,x,y)$. 

\begin{thm}\label{t:largetime}
Suppose  $V\in \sK(\beta,\kappa)$. 

\noindent (i) If $d\ge 2$, then	 there exist $c_i>0$, $1\le i \le 4$, such that  for all $t \ge 4$ and $x,y \in \R^d_0$,
\begin{equation}\label{e:largetime}
\!\!\!c_1 H_{d,\beta,\kappa}(t,|x|) H_{d,\beta,\kappa}(t,|y|)q(c_2t,x,y\!)\!\le	\!
p(t,x,y)
\!\le\! c_3 H_{d,\beta,\kappa}(t,|x|) H_{d,\beta,\kappa}(t,|y|) q(c_4t,x,y).
\end{equation}
(ii) If  $d=1$, then \eqref{e:largetime}  holds for all $t \ge 4$ and $x,y \in \R^1_0$ with $xy>0$.
\end{thm}

\medskip

Throughout Subsections \ref{ss:Large-lower} and \ref{ss:Large-upper}, we assume that $V\in \sK(\beta,\kappa)$ and  that $\gamma>0$ is
the constant in \eqref{e:V-cond-3}. We simply write $h$ and $H_d$ instead of $h_{\beta,\kappa}$ and $H_{d,\beta,\kappa}$ respectively.

We introduce 
the following auxiliary function $\wt q_2$:
\begin{align*}
	\wt q_2(t,x,y):=t^{-d/2}\bigg( 1\wedge \frac{ t^{1/2}}{|x-y|}\bigg)^{d+2}, \qquad t>0,\; x,y \in \R^d.
\end{align*}
By  \eqref{e:exp-poly}, we see that there exists $C=C(d)>0$ such that
\begin{align}\label{e:q-q2}
	q(t,x,y) \le C\,\wt q_2(t,x,y) \quad \text{for all $t>0$ and $x,y \in \R^d$}.
\end{align}

\subsection{Large time lower estimates}\label{ss:Large-lower}

\begin{lem}\label{l:one-step-boundary}
	Assume $d=1$ or $d=2$. There exists  $C>0$ such that for all $R\ge 1$, $t\ge R^2$ and $w \in B(0,2R)^c$,
	\begin{align*}
		\int_0^t \int_{B(0,R)^c} \psi_{d,R}(s,|w|) 
		\psi_{d,R}(s,|z|) \psi_{d,R}(t,|z|)  \frac{dz}{|z|^{-2-\gamma}} ds \le C\frac{ t^{d/2} \psi_{d,R}(t,|w|)}{R^{\gamma}}
	\end{align*}
	and
	\begin{align*}
		\int_0^t \int_{B(0,R)^c} \psi_{d,R}(s,|w|) \psi_{d,R}(s,|z|) \psi_{d,R}(t,|z|) \wt q_2(s,w,z) \frac{dz}{|z|^{2+\gamma}} ds \le C\frac{ \psi_{d,R}(t,|w|)}{R^\gamma}.
	\end{align*}
\end{lem}
\pf 
When $d=1$, the  results can be proved by following the arguments in \cite[Lemmas 5.5 and 5.6]{CS24}, with $\alpha$ replaced by $2$. 

 Suppose $d=2$. Define 
\begin{align*}
	\Psi_R(s,r):=  1 \wedge \frac{ \Lg(r/R)  }{  	\Lg(\sqrt s/R) }, \quad s>0, \; r \in [R, \infty).
\end{align*}
Note that
\begin{align}\label{e:Psi-1}
 	\Psi_R(s,r) \ge \psi_{2,R}(s,r) \quad \text{for all $s>0$ and $r \in [R,\infty)$}.
\end{align}
Indeed, if $s\ge R^2$, then  \eqref{e:Psi-1} is evident. If $s<R^2$, then   
$\Psi_R(s,r)=1 \ge \psi_{2,R}(s,r)$. 
For all   $0<s\le t$ and  $w, z \in B(0,R)^c$,  using \eqref{e:Psi-1} and the fact that $u \mapsto \Psi_R(u,|w|) \Lg (\sqrt u/R)$ is non-decreasing,  we obtain
\begin{align}\label{e:Psi-3}
	\psi_{2,R}(s,|w|) \psi_{2,R}(t,|z|) \le  \frac{\Psi_{R}(s,|w|)\Lg (|z|/R) }{\Lg(\sqrt t/R)}  \le  \frac{\Psi_{R}(t,|w|)\Lg (|z|/R) }{\Lg(\sqrt s/R)} .
\end{align}
Using \eqref{e:Psi-3} and $\psi_{2,R}\le 1$ in the first inequality below,  and \eqref{e:Log-scaling} with $\eps=1$ and $\eps=\gamma/2$ in the second, we get
\begin{align}\label{e:one-step-boundary-1}
	&	\int_0^t \int_{B(0,R)^c} \psi_{2,R}(s,|w|) \psi_{2,R}(s,|z|) \psi_{2,R}(t,|z|)  \frac{dz}{|z|^{2+\gamma}} ds\nn\\
	& \le 
	\Psi_R(t,|w|)
		\int_0^{t} \int_{B(0,R)^c} \frac{\Lg(|z|/R)}{|z|^{2+\gamma} \,\Lg(\sqrt s/R)}dzds\nn\\
	& \le \frac{c_1 \sqrt t 
		\Psi_R(t,|w|) 
		\,\Lg 1}{R^{\gamma/2}\Lg (\sqrt t/R)}	\int_0^{t} \frac{ds}{\sqrt s} \int_{B(0,R)^c} \frac{dz}{|z|^{2+\gamma/2}} = \frac{c_2  t
		\Psi_R(t,|w|) 
		}{R^{\gamma}\Lg (\sqrt t/R)}  \le \frac{c_2  t
		\Psi_R(t,|w|) 
		}{R^{\gamma}} . \qquad 
\end{align}
Further, using \eqref{e:Psi-1} and  \eqref{e:Psi-3},
 we also get
\begin{align*}
	&	\int_0^t \int_{B(0,R)^c} \psi_{2,R}(s,|w|) \psi_{2,R}(s,|z|) \psi_{2,R}(t,|z|) \wt q_2(s,w,z) \frac{dz}{|z|^{2+\gamma}} ds\nn\\
	&\le 
	 \Psi_R(t,|w|)
	\bigg( \int_0^{R^2} + 	\int_{R^2}^t \bigg)  \int_{B(0,R)^c}     \Psi_{R}(s,|z|)   \wt q_2(s,w,z) \frac{\Lg(|z|/R)}{|z|^{2+\gamma} \Lg (\sqrt s/R)} dz\,ds\nn\\
	&=: \Psi_R(t,|w|) (I_1+I_2).
\end{align*}
Note that there exists $c_3>0$ independent of $w$ such that $\int_{\R^2} \wt q_2(s,w,z) dz = c_3$ for all $s>0$. Using this, $\Psi_R(s,|z|) \le 1$ and \eqref{e:Log-prop-2},  we get
\begin{align}\label{e:one-step-boundary-2}
	I_1\le \frac{ \Lg 1}{R^{2+\gamma} \Lg 0}\int_0^{R^2} \int_{B(0,R)^c}   \wt q_2(s,w,z) dz\, ds \le \frac{c_3
		}{R^{\gamma} 
		}.
\end{align}
For $I_2$, using \eqref{e:Log-scaling} with $\eps= \gamma/4$, we see that
\begin{align}\label{e:one-step-boundary-3}
	I_2&\le 
	\int_{R^2}^t \int_{B(0,R)^c}    \frac{(\Lg(|z|/R))^2 }{|z|^{2+\gamma}\, s(\Lg (\sqrt s/R))^2} dz\,ds\nn\\
	&\le \frac{c_4 (\Log 1)^2 
		}{R^{\gamma/2}}\int_{B(0,R)^c}    \frac{dz}{|z|^{2+\gamma/2} }  \int_{R^2}^t \frac{ds}{ s(\Lg (\sqrt s/R))^2 } \le \frac{c_5 
		  }{R^{\gamma} } .
\end{align}
Since $t\ge R^2$ and $|w| \ge 2R$, by \eqref{e:Log-prop-2}, we have
\begin{align*}
	\psi_{2,R}(t,|w|) \ge  1 \wedge \frac{\Lg ( (|w|-R)/R)}{ \Lg ( \sqrt t/R)}   \ge 1 \wedge \frac{\Lg ( |w|/(2R))}{ \Lg ( \sqrt t/R)}  \ge \frac{1}{2} \Psi_R(t,|w|). 
\end{align*}
Thus, by \eqref{e:one-step-boundary-1}, \eqref{e:one-step-boundary-2} and \eqref{e:one-step-boundary-3}, we conclude that the results are  also valid for $d=2$.
\qed

Note that, following the arguments for \cite[Lemma 5.3]{CS24}, we get that for $d\ge 3$,
\begin{align}\label{e:sdB}
	\int_{B(0,R)^c} \frac{dz}{|w-z|^{d-2}|z|^{2+\gamma}} \le \frac{c_2}{R^{\gamma}} \quad \text{for all $w \in B(0,2R)^c$}.
\end{align}
Following the arguments in \cite[Proposition 5.2]{CS24}, with $\alpha$ replaced by $2$, and using Lemma \ref{l:one-step-boundary} and \eqref{e:sdB}, we can obtain the next lemma.  We omit the details.

\begin{lem}\label{l:one-step}
There exists $C>0$ such that for all  $R\ge 1$, $t\ge R^2$  and  $x,y \in B(0,2R)^c$,
 \begin{align*}
&\int_0^t \int_{B(0,R)^c}\psi_{d,R}(s,|x|) \psi_{d,R}(t-s,|y|) \psi_{d,R}(s,|z|) 	\psi_{d,R}(t-s,|z|)  \wt q_2(s,x,z) 
	\wt q_2(t-s,z,y)   
	V(z)dzds \nn\\
&	\le CR^{-\gamma}\psi_{d,R}(t,|x|) \psi_{d,R}(t,|y|)  \wt q_2(t,x,y).
 \end{align*}
\end{lem}

\begin{lem}\label{l:largetime-lower-bound}
(i) Assume $d\ge 2$. There exist $R_1\ge 2$ and $C>0$ such that  for all $t\ge R_1^2$ and  $x,y \in B(0,2R_1)^c$ with $|x-y| \le 6t^{1/2}$,
	\begin{align}\label{e:largetime-lower-bound-claim}
	p^{\overline B(0,R_1)^c}(t,x,y) \ge  C \psi_{d,R_1}(t,|x|)  \psi_{d,R_1}(t,|y|) t^{-d/2}.
	\end{align}
\noindent (ii) If $d=1$, then  \eqref{e:largetime-lower-bound-claim} holds for  all $t\ge R_1^2$ and  $x,y \in B(0,2R_1)^c$ with $xy>0$ and $|x-y| \le 6t^{1/2}$.
\end{lem}
\pf 
By \cite[(A.17)]{BJM24},  originally due to \cite[(41)]{BHJ08} and \cite[Appendix]{JW}, we obtain for all $R>0$, $t>0$ and $x,y \in \overline B(0,R)^c$,
\begin{equation}\label{e:largetime-lower-bound}
	\frac{p^{\overline B(0,R)^c} (t,x,y)}{q^{\overline B(0,R)^c}(t,x,y)} \ge  \exp \bigg[ -   \int_0^t \int_{\overline B(0,R)^c} \frac{q^{\overline B(0,R)^c}(s,x,z) q^{\overline B(0,R)^c}(t-s,z,y)} {q^{\overline B(0,R)^c}(t,x,y)} V(z) dz  ds \bigg].  
\end{equation}
By Proposition \ref{p:DHKE}, we have for all $R>0$, $t>0$ and $x,y \in \overline B(0,R)^c$ with $|x-y|\le 6t^{1/2}$,
\begin{align}
	q^{\overline B(0,R)^c}(t,x,y) &\ge c_1^{-1}\psi_{d,R}(t,|x|)\psi_{d,R}(t,|y|) t^{-d/2}\label{e:largetime-lower-bound-0}\\
	& \ge c_1^{-1}\psi_{d,R}(t,|x|)\psi_{d,R}(t,|y|) \wt q_2(t,x,y).\label{e:largetime-lower-bound-0+}
\end{align}
Further, it follows from Proposition \ref{p:DHKE} and \eqref{e:q-q2}
that for all $R>0$, $t>0$ and $x,y \in  \overline B(0,R)^c$,
\begin{align}\label{e:largetime-lower-bound-1}
	q^{\overline B(0,R)^c}(t,x,y) \le c_2 \psi_{d,R}(t,|x|) \psi_{d,R}(t,|y|) \wt q_2(t,x,y).
\end{align}
 By \eqref{e:largetime-lower-bound-0+}, \eqref{e:largetime-lower-bound-1} and  Lemma \ref{l:one-step}, there exists $R_1\ge 2$ such that for all $t \ge R_1^2$ and $x,y \in B(0,2R_1)^c$ with $|x-y|\le 6t^{1/2}$,
\begin{align}\label{e:largetime-lower-bound-2}
	& \int_0^t \int_{\overline B(0,R_1)^c} \frac{q^{\overline B(0,R_1)^c}(s,x,z) q^{\overline B(0,R_1)^c}(t-s,z,y)} {q^{\overline B(0,R_1)^c}(t,x,y)} V(z) dz  ds\nn\\
	&\le c_3\int_0^t \int_{B(0,R_1)^c} \frac{\psi_{d,R_1}(s,|x|) \psi_{d,R_1}(t-s,|y|) \psi_{d,R_1}(s,|z|) 	\psi_{d,R_1}(t-s,|z|)}{ \psi_{d,R_1}(t,|x|) \psi_{d,R_1}(t,|y|)  \wt q_2(t,x,y) } \nn\\
	&\qquad \qquad \qquad \quad  \times   \wt q_2(s,x,z) 
	\wt q_2(t-s,z,y)   	V(z)dzds	\le c_4R_1^{-\gamma} \le 1.
\end{align}
Combining \eqref{e:largetime-lower-bound}, \eqref{e:largetime-lower-bound-0} and  \eqref{e:largetime-lower-bound-2},   we arrive at the desired result.\qed

\begin{lem}\label{l:largetime-lower-bound-full}
	(i) Assume $d\ge 2$. There exist $c_1,c_2>0$ such that 
	\begin{align}\label{e:largetime-lower-bound-off-diagonal}
		p(t,x,y) \ge  c_1 \psi_{d,R_1}(t,|x|)  \psi_{d,R_1}(t,|y|) t^{-d/2} e^{- {c_2|x-y|^2}/{t}}
	\end{align}
for all $t\ge R_1^2$ and  $x,y \in B(0,2R_1)^c$,
	where $R_1\ge 2$ is the constant in Lemma \ref{l:largetime-lower-bound}.
	
	\noindent (ii) If $d=1$, then  \eqref{e:largetime-lower-bound-off-diagonal} holds for  all $t\ge R_1^2$ and  $x,y \in B(0,2R_1)^c$ with $xy>0$. 
\end{lem}
\pf Let  $t\ge R_1^2$ and  $x,y \in B(0,2R_1)^c$. When $d=1$, we  assume  $xy>0$. If $|x-y|\le 6t^{1/2}$,  then the result follows from Lemma \ref{l:largetime-lower-bound}. 
Suppose $|x-y| >6t^{1/2}$. Without loss of generality, we assume $|x|\le |y|$.
Note that $|y| \ge |x-y|/2>3t^{1/2} \ge 2R_1 + t^{1/2}$. Let $n\ge 5$ be such that $\sqrt{n-1}<|x-y|/t^{1/2} \le \sqrt n$. Observe that there exist  $c_1> 1$ and $c_2\in (0,1)$ depending on $d$ only such that for any $\delta \in (0,1)$, there exist
    $k_\delta \ge 3$ and $z^\delta_0,z^\delta_1,\cdots, z^\delta_k \in \R^d_0$ with the following properties: (1) $k_\delta\le c_1|x-y|/(\delta(t/n)^{1/2})$, (2) $z^\delta_0=x$ and $z^\delta_{k_\delta}=y$, (3) $ |z^\delta_i| \ge 2R_1+c_2\delta(t/n)^{1/2}$ for all $1\le i\le k_\delta$ and  
 (4) $|z^\delta_{i}-z^\delta_{i-1}| \le \delta(t/n)^{1/2}$ for all $1\le i \le k_\delta$.
 Set $\delta:=1/c_1$. For any $2\le i\le k_\delta$, $w_{i-1} \in B(z^\delta_{i-1}, c_2\delta(t/n)^{1/2}/2)$ and  $w_{i} \in B(z^\delta_{i}, c_2\delta(t/n)^{1/2}/2)$, we have
\begin{align}\label{e:chain-1}
	|w_{i-1}| \wedge  |w_i| \ge 2R_1+ c_2\delta(t/n)^{1/2}/2 \quad \text{and} \quad |w_{i-1}-w_i| \le 2\delta(t/n)^{1/2}.
\end{align} 
Note that  $k_\delta \ge |x-y|/ (\delta(t/n)^{1/2}) \ge \sqrt{n(n-1)} /\delta $ and  $k_\delta \le c_1 |x-y| /(\delta(t/n)^{1/2}) \le c_1  n/\delta$. Hence, we have $t/k_\delta \le 2\delta t/n$ and 
\begin{align}\label{e:chain-2}
	2(t/k_\delta)^{1/2} \ge 2(\delta t/(c_1n))^{1/2} = 2\delta(t/n)^{1/2}.
\end{align}
In particular, for all $1\le i\le k_\delta$ and  $w_{i} \in B(z^\delta_{i}, c_2\delta(t/n)^{1/2}/2)$,  since $|w_i| -R_1  \ge   
R_1+c_2(t/n)^{1/2}/2$, 
by  the definition of $\psi_{d,R_1}$ and \eqref{e:Log-prop-2}, we see that 
\begin{align}\label{e:chain-3}
	\psi_{d,R_1}(t/k_\delta, |w_i|) &\ge  \psi_{d,R_1}( 2\delta t/n, |w_i|) \nn\\
	&\ge \begin{cases}
		\displaystyle 1 \wedge \frac{R_1}{(2\delta t/n)^{1/2}\wedge R_1}
		&\mbox{ if $d\ge 3$}, \\[8pt]
		\displaystyle   1 \wedge 
		\frac{R_1 \,\Lg ( R_1+ (c_2(t/n)^{1/2}/2)/R_1)}{((2\delta t/n)^{1/2} \wedge R_1)\, \Lg ( (2\delta t/n)^{1/2}/R_1)} 
		&\mbox{ if $d=2$}, \\[8pt]
		\displaystyle  1 \wedge \frac{(c_2(t/n)^{1/2}/2)}{(2\delta t/n)^{1/2}} &\mbox{ if $d= 1$}.
	\end{cases}\nn\\
	&\ge c_3.
\end{align}
Applying Lemma \ref{l:largetime-lower-bound}, from \eqref{e:chain-1}, \eqref{e:chain-2} and \eqref{e:chain-3}, we deduce that for any $1\le i\le k_\delta$, $w_{i-1} \in B(z^\delta_{i-1}, c_2\delta(t/n)^{1/2}/2)$ and  $w_{i} \in B(z^\delta_{i}, c_2\delta(t/n)^{1/2}/2)$,
\begin{align}\label{e:chian-4}
	p(t/k_\delta,w_{i-1}, |w_i|)  \ge c_4 \begin{cases}
		\psi_{d,R_1}(t/k_\delta, |w_{i-1}|)	(t/k_\delta)^{-d/2} &\mbox{ if $i = 1$},\\[2pt]
		(t/k_\delta)^{-d/2} &\mbox{ if $i \ge 2$}.
	\end{cases} 
\end{align}
Using the semigroup property,  \eqref{e:chian-4}, $k_\delta \le c_1n/\delta$ and $n-1<|x-y|^2/t$, we arrive at
\begin{align*}
	&p(t,x,y)\nn\\
	& \ge \int_{B(z^\delta_{1}, c_2\delta(t/n)^{1/2}/2)} \cdots \int_{B(z^\delta_{k_\delta-1}, c_2\delta(t/n)^{1/2}/2)} p( t/k_\delta, x, w_1) \cdots p(t/k_\delta, w_{k_\delta-1},y) \, dw_1 \cdots dw_{k_\delta-1} \nn\\
	&\ge  \psi_{d,R_1}(t/k_\delta,|x|)\,(c_4(t/k_\delta)^{-d/2})^{k_\delta} (c_5 \delta(t/n)^{1/2})^{d(k_\delta-1) } \nn\\
	& = \psi_{d,R_1}(t/k_\delta,|x|) (c_5\delta)^{-d}
	(c_4c_5^d \delta^dk_\delta^{d/2}/n^{d/2})^{k_\delta} 
	n^{d/2} t^{-d/2}\nn\\
	&  \ge c_6 \psi_{d,R_1}(t,|x|)t^{-d/2} \exp\left( - \frac{c_1c_7n}{\delta}\right) \ge c_6  e^{-c_1c_7/\delta}\psi_{d,R_1}(t,|x|) t^{-d/2} \exp \bigg( - \frac{c_1c_7 |x-y|^2}{\delta t}\bigg).
\end{align*}
This completes the proof. \qed

\noindent \textbf{Proof of  Theorem \ref{t:largetime} (Lower estimates).} 
 Let $t\ge 4$ and $x,y \in \R^d_0$. 
  Without loss of generality, we assume that $|x|\le |y|$. 
  When $d=1$, we additionally assume that $xy>0$.   
  Let  $R_1\ge 2$ be the constant in Lemma \ref{l:largetime-lower-bound}.
  If  $4\le t\le 4R_1^2$, then 
$$h( \eta_1 t^{1/(2+\beta)}) \asymp   \exp \Big(- \frac{t}{(|x| \vee |y| \vee t^{1/(2+\beta)})^{2+2\beta}}\Big)\asymp   \Lg t \asymp \sqrt t \asymp 1.$$
Thus,  applying Theorem \ref{t:smalltime}(ii) with $T=4R_1^2$, we get
 \begin{align*}
	p(t,x,y) &\ge c_1( 1 \wedge  h(|x| )) ( 1 \wedge  h(|y| ))q(c_2t,x,y) \ge  c_3( 1 \wedge  H_{d}(t,|x| )) ( 1 \wedge  H_{d}(t,|y| ))q(c_2t,x,y),
\end{align*}
proving the desired lower bound.

 Assume  $t > 4R_1^2$. 
 Define 
$
 	f(t):= \1_{\{d\ge 3\}} + 	(\Lg t)^{-1}  \1_{\{d=2\}} + t^{-1/2}\1_{\{d=1\}}.
$
 Let $z_0 \in \R^d$ be such that $|z_0|=3R_1$. For any $z \in B(z_0,R_1)$, we have
\begin{align}\label{e:largetime-LHK-0}
	\psi_{d,R_1}(t,|z|) \ge  c_4f(t).
\end{align}
Further, for any $v \in \{x,y\}$ and $z \in B(z_0,R_1)$, since $|z| >2R_1$ and $|v-z| \le |v| + 4R_1$, we get from  Theorem \ref{t:smalltime} that
\begin{align}\label{e:largetime-LHK-1}
	p(1,v,z)=p(1,z,v)&\ge c_5(1 \wedge h(|v|)) q(c_6, v,z) \ge c_7 (1 \wedge h(|v|)) e^{-c_8|v|^2}.
\end{align}
We distinguish between  three  cases:

\smallskip

\noindent
{\it Case 1:  $|y| <8R_1$.}  For any $z,w \in B(z_0,R_1)$, 
since $t-2 \ge 3R_1^2$, 
$|z| \wedge |w| >2R_1$ and $|z-w|<2R_1$, using Lemma \ref{l:largetime-lower-bound-full} and \eqref{e:largetime-LHK-0}, we get
\begin{equation}\label{e:largetime-LHK-2}
	p(t-2,z,w) \ge c_9\psi_{d,R_1}(t,|z|)  \psi_{d,R_1}(t,|w|) t^{-d/2}  \ge c_{10} f(t)^2t^{-d/2} 
\end{equation}
By the semigroup property, the almost monotonicity of $h$, \eqref{e:largetime-LHK-1} and \eqref{e:largetime-LHK-2}, we obtain
\begin{align*}
	p(t,x,y) &\ge \int_{B(z_0,R_1)} \int_{B(z_0,R_1)} p(1,x,z) p(t-2,z,w) p(1,w,y) dzdw\\
	&\ge  c_{11} |B(0,R_1)|^2  e^{-2c_8(8R_1)^2}  h(|x|) h(|y|)  f(t)^2t^{-d/2} 
	\ge   c_{12} 	 H_{d}(t,|x|)  H_{d}(t,|y|) t^{-d/2}.
\end{align*}

\noindent
{\it Case 2:  $|x| < 2R_1$ and $|y| \ge 8R_1$.} For all $z \in B(z_0,R_1)$, we have
$
|y-z| < |y| + 4R_1 \le 2|y| - 2|x| \le 2|y-x|.
$
Thus, by Lemma \ref{l:largetime-lower-bound-full}, \eqref{e:compare-H-psi} and \eqref{e:largetime-LHK-0}, we obtain for all $z \in B(z_0,R_1)$, 
\begin{align}\label{e:largetime-LHK-3}
	p(t-1,z,y) &\ge c_{13} \psi_{d,R_1}(t,|z|)  \psi_{d,R_1}(t,|y|) t^{-d/2} e^{ - {c_{14}|z-y|^2}/{t}} \nn\\
	&\ge c_{15}  f(t) H_{d}(t,|y|) t^{-d/2} e^{ - {4c_{14}|x-y|^2}/{t}}.
\end{align}
Using the semigroup  property, \eqref{e:largetime-LHK-1} and \eqref{e:largetime-LHK-3}, we get
\begin{align*}
	p(t,x,y)&\ge \int_{B(z_0,R_1)}  p(1,x,z) p(t-1,z,y) dz\\
	&\ge  c_{16} |B(0,R_1)|    e^{-c_8(2R_1)^2} h(|x|) f(t) H_{d}(t,|y|) t^{-d/2} e^{ - {4c_{14}|x-y|^2}/{t}} \\
	&\ge   c_{17}	 H_{d}(t,|x|)  H_{d}(t,|y|) t^{-d/2}e^{ - {4c_{14}|x-y|^2}/{t}}.
\end{align*}

\noindent
{\it Case 3:  $|x| \ge 2R_1$ and $|y| \ge 8R_1$.} The result follows from Lemma \ref{l:largetime-lower-bound-full} and \eqref{e:compare-H-psi}. 

\smallskip

The proof is complete. \qed

\subsection{Large time upper  estimates}\label{ss:Large-upper}
Recall that  we have assumed that $V\in \sK(\beta,\kappa)$ and  that $\gamma>0$ is
the constant in \eqref{e:V-cond-3}. Note that, when $d\ge 3$,  the  upper bound 
in Theorem \ref{t:largetime}
 follows  from  Corollary \ref{c:UHK-general}. In this section, we focus on the case $d \le 2$ and 
 prove the  upper bound in this case.

For all  $k \ge 1$,  $t>0$, 
$x\in B(0,kt^{1/2})$  and $y \in \R^d$, it holds that
\begin{align*}
	\frac{|x-y|^2}{t} \le \frac{2(|y|^2 + |x|^2)}{t}\le \frac{2|y|^2}{t} + 2k \quad \text{ and } \quad 	\frac{|x-y|^2}{t} \ge \frac{|y|^2/2 - |x|^2}{t}\ge \frac{|y|^2}{2t} -k.
\end{align*}
Thus, for any $k \ge 1$, there exists  a constant $C=C(k)>1$ such  that for all $t>0$, $x\in B(0,kt^{1/2})$ and $y \in \R^d$, 
\begin{align}\label{e:q-estimate-1}
C^{-1} q(t/2,0,y)\le 	q(t,x,y) \le C q(2t,0,y).
\end{align}

We first give  several lemmas that will be used in the proof of the upper estimates in Theorem \ref{t:largetime}.

\begin{lem}\label{l:UHK-general-estimate}
	For any $k \ge 1$,	there exist constants $\lambda_0 \in (0,1/2]$ independent of $k$ and $C=C(k)>0$ such that  for all $t \ge 1$ and $x,y \in \R^d_0$ with $|x| <kt^{1/2}$,
	\begin{align*}
		p(t,x,y) &\le  Ch(|x|) t^{-d/2} e^{-{\lambda_0 |y|^2}/{(4t)}}.
	\end{align*}
\end{lem}
\pf The result follows from Corollary \ref{c:UHK-general} and the second inequality in \eqref{e:q-estimate-1}. \qed 

\begin{lem}\label{l:largetime-2-upper}
 There exist  constants $C,k\ge 1$ such that for all $t\ge 2$ and  $x,y \in \R^1_0$,
	\begin{align}\label{e:largetime-2-upper}
		p(t,x,y) \le C H_{1}(t,|x|) q(kt,x,y).
	\end{align}
\end{lem}
\pf  Let $t \ge 2$ and $x,y \in \R^1_0$.  By Proposition \ref{p:DHKE}, we have  for all $s>0$ and $v,w \in \R^1_0$,
\begin{align}\label{e:largetime-2-upper-1}
&	p(s,v,w)
	=\lim_{n \to \infty}    p^{\overline B(0,1/n)^c}(s,v,w) \nn\\
	&\le \lim_{n \to \infty}    q^{\overline B(0,1/n)^c} (s,v,w)   \le  c_1  \Big( 1 \wedge \frac{ |v|}{  \sqrt s}\Big) q(c_2s,v,w).
\end{align}
If $|x| \ge 1$, then by choosing $k$ larger than $c_2$, \eqref{e:largetime-2-upper} follows from  
\eqref{e:largetime-2-upper-1} and 
\eqref{e:H-arbitrary-R}.

Suppose $|x| <1$.  
By Lemma \ref{l:UHK-general-estimate} and \eqref{e:exp-poly},  we obtain for all $z \in \R^1_0$,
\begin{equation}\label{e:largetime-2-upper-2}
 p(1,x,z) |z| \le  c_3 h(|x|) |z| e^{-c_4 |z|^2}  \le c_5  h(|x|) e^{-c_4|z|^2/2} .
\end{equation}
 Using the semigroup property in the first line below,  \eqref{e:largetime-2-upper-1}-\eqref{e:largetime-2-upper-2}  and the fact that $p(s,w,y) \le q(s,w,y)$ for all $s>0$ and $w \in \R^1_0$ in the second,  and $(t-1)/2 \ge t/4$ in the fourth,  we obtain
\begin{align}\label{e:largetime-2-upper-3}
	p(t,x,y)& =    \int_{\R^1_0} \int_{\R^1_0}  p(1,x,z) \, p((t-1)/2, z,w) \, p( (t-1)/2,w,y) \, dzdw \nn\\
	&\le \frac{c_6 h(|x|)}{\sqrt{(t-1)/2}}  \int_{\R} \int_{\R} e^{-c_4|z|^2/2}  q(c_2(t-1)/2, z,w) \,  q((t-1)/2,w,y) dzdw\nn\\
		&= \frac{ c_6  h(|x|)}{\sqrt{(t-1)/2}}  \int_{\R} e^{-c_4|z|^2/2}  q((c_2+1)(t-1)/2, z,y) dz \nn\\
		&\le \frac{ c_7 h(|x|)}{t}  \int_{\R}  \exp \bigg( - \frac{c_4|z|^2}{2} - \frac{c_8|y-z|^2}{t} \bigg) dz .
\end{align}
Using $t\ge 2$ and the inequality $a^2 + b^2 \ge (a+b)^2/2$ for $a,b \in \R$, we get
\begin{align}\label{e:largetime-2-upper-4}
	\int_{\R}  \exp \bigg( - \frac{c_4|z|^2}{2} - \frac{c_8|y-z|^2}{t} \bigg) dz &  \le \int_{\R}  \exp \bigg( - \frac{c_4|z|^2}{4} - \frac{c_4|z|^2}{2t} - \frac{c_8|y-z|^2}{t} \bigg) dz \nn\\
		&  \le  e^{-c_9|y|^2/t}\int_{\R}   e^{-c_4|z|^2/4} dz  = c_{10}e^{-c_9|y|^2/t}.
\end{align}
Combining \eqref{e:largetime-2-upper-3} with \eqref{e:largetime-2-upper-4} and using \eqref{e:q-estimate-1}, we arrive at
\begin{align*}
	p(t,x,y) \le \frac{c_{11} h(|x|)}{\sqrt t} q(t/(4c_9),0,y) \le \frac{
		c_{12}h(|x|)}{\sqrt t} q(t/(2c_9),x,y).
\end{align*}
The proof is complete. \qed

The next result follows from \cite[Lemma 8.1]{CS24} (with $M=\R^d_0$ and $U=\overline B(0,\eps)^c$).
\begin{lem}\label{l:Dirichlet-upper}
	Let $t>0$ and $\eps>0$. 	If there exists a  continuous function 
	$F_{t,\eps}$ on $\R^d_0$ such that
	\begin{align*}
		\sup_{s\in [t/2,t],\, z \in \overline B(0,\eps)\setminus \{0\}}p(s, z,y) \le F_{t,\eps}(y) \quad \text{for all} \;\, y \in \R^d_0,
	\end{align*}
	then we have
	\begin{align*}
		p(t,x,y) \le p^{\overline B(0,\eps)^c}(t,x,y) + F_{t,\eps}(x) + F_{t,\eps}(y)  \quad \text{for all} \;\, x,y \in \R^d_0.
	\end{align*}
\end{lem}

For $n \ge 1$, we denote by  $\Lg^n  r:=\Lg \circ \cdots \circ \Lg r$ the $n$-th iterate of the function $\Lg$.

\begin{lem}\label{l:criticial-upper-1}
	There exist  $C\ge1$ and $\lambda_1\in (0,1/2]$ such that the following holds: If there exist $n \in \N$, $a_1\ge 1$ and $a_2 \in (0,1/2]$  such that for all $t\ge 1$  and $z,y \in \R^2_0$ with $|z|\le 1$, 
	\begin{equation}\label{e:criticial-upper-1-ass}
		p(t,z,y) \le 
		\frac{a_1 h(|z|)\, \Log^n \sqrt t }{t\,\Lg \sqrt t}  \exp \Big(- \frac{a_2|y|^2}{4t}\Big),
	\end{equation}
	then for all $t\ge 2$ and $x,y \in \R^2_0$ with  $|x| \le \sqrt{2t}$, 
	\begin{equation}\label{e:criticial-upper-1}
		p(t,x,y) \le 
		\frac{C (\Lg a_1)(\Lg |x| + \Log^{n+1} \sqrt t ) }{t\,\Lg \sqrt t} \exp \Big(- \frac{(\lambda_1 \wedge a_2)|y|^2}{4t}\Big).
	\end{equation}
\end{lem}
\pf Let $t\ge 2$ and $x,y\in \R^2_0$ with  $|x| \le \sqrt {2t}$. Set
$$
\eps:=  \frac{\Lg |x| \wedge \Log^n \sqrt t}{a_1\Log^n \sqrt t} \in (0,1].
$$ 
Since $h$ is exponentially increasing at $0$, there exists $c_1>0$ such that 
\begin{align}\label{e:criticial-upper-1-2}
	h(|z|) \le c_1  |z| \quad \text{for all $z \in B(0,1)$}.
\end{align}
If $|x| \le \eps$, then by \eqref{e:criticial-upper-1-ass} and \eqref{e:criticial-upper-1-2}, 
\begin{align*}
	p(t,x,y) \le 
	\frac{c_1a_1 \eps\, \Log^n \sqrt t }{t\,\Lg \sqrt t}  \exp \Big(- \frac{a_2|y|^2}{4t}\Big) \le 
	\frac{c_1 \Lg |x| }{t\,\Lg \sqrt t}  \exp \Big(- \frac{a_2|y|^2}{4t}\Big).
\end{align*}
Thus,  since $\Lg a_1 \ge \Lg 1= 1$,  taking $C$ larger than $c_1$, we get \eqref{e:criticial-upper-1} in this case.

Assume $|x|>\eps$. By \eqref{e:q-estimate-1} and Proposition \ref{p:DHKE}, there exist constants $c_2,c_3>0$ and $\lambda_1\in (0,1/2]$ independent of $t,x$ and $a_1$  such that for all $w\in \R^2_0$,
\begin{align}\label{e:criticial-upper-1-0}
	\frac{ p^{\overline B(0,\eps)^c}(t,x,w)}{q(t/\lambda_1 ,0,w)}
	&	\le 	\frac{
		c_2 q^{\overline B(0,\eps)^c}(t,x,w)}{q(t/(2\lambda_1) ,x,w)}\le \frac{c_3 \eps \, \Lg (|x|/\eps)}{(\sqrt t \wedge \eps) \Lg (\sqrt t/\eps)}= \frac{c_3   \Lg (|x|/\eps)}{ \Lg (\sqrt t/\eps)}.
\end{align}
Using \eqref{e:Log-prop-1}  and \eqref{e:Log-prop-2}, since $\Lg a_1 \ge  1$ and $\Lg r\ge \log(e-1)$ for all $r\ge 0$,  we see that
\begin{align}\label{e:criticial-upper-1-1}
	\frac{\Lg (|x|/\eps)}{ \Lg (\sqrt t/\eps)}  &  \le \frac{\Lg |x| + \Lg a_1 + \Lg (  (\Log^n \sqrt t)/\log(e-1) ) }{\Lg \sqrt t} \nn\\
	&\le  \frac{ (\Lg a_1) (\Lg |x|  + 1 + (\Log^{n+1} \sqrt t)/\log (e-1)  )}{\Lg \sqrt t}\nn\\
	&\le  \frac{c_4 (\Lg a_1)(\Lg |x|  + \Log^{n+1} \sqrt t)}{\Lg \sqrt t}.
\end{align}
We deal with the cases $|y|\le \sqrt{6t/a_2}$ and $|y|>\sqrt{6t/a_2}$ separately.

\smallskip

\noindent
{\it Case 1: $|y|\le \sqrt{6t/a_2}$.} By \eqref{e:criticial-upper-1-ass}, \eqref{e:criticial-upper-1-2} and  \eqref{e:Log-prop-2},  we have  for all $w \in \R^2_0$,
\begin{equation}\label{e:criticial-upper-1-3}
	\sup_{s\in [t/2,t],\, z\in \overline B(0,\eps)\setminus \{0\}} p (s,z,w)  \le	\frac{c_5 a_1  \eps \, \Log^n \sqrt t }{(t/2)\,\Lg \sqrt{t/2}}  \exp \Big(- \frac{a_1|w|^2}{4t}\Big)  \le  \frac{2^{3/2}c_5 \Lg |x| }{t\,\Lg \sqrt t}  .
\end{equation}
Applying Lemma \ref{l:Dirichlet-upper}, and using \eqref{e:criticial-upper-1-0}, \eqref{e:criticial-upper-1-1}, \eqref{e:criticial-upper-1-3} and $a_2|y|^2\le 6t$,  we deduce that
\begin{align*}
	p(t,x,y) &\le  \frac{c_6 (\Lg a_1)(\Lg |x|  + \Log^{n+1} \sqrt t)}{t\,\Lg \sqrt t}  \exp \Big( - \frac{\lambda_1|y|^2}{4t}\Big)   + \frac{2^{5/2}c_5 \Lg |x| }{t\,\Lg \sqrt t}\\
	&\le  \frac{(c_6 + 2^{5/2}e^{3/2} c_5)(\Lg a_1)(\Lg |x|  + \Log^{n+1} \sqrt t)}{t\,\Lg \sqrt t}  \exp \Big( - \frac{(\lambda_1\wedge a_2)|y|^2}{4t} \Big) .
\end{align*}

\noindent
{\it Case 2: $|y|> \sqrt{6t/a_2}$.} By the strong Markov property, we have for a.e. $w \in B(0,\sqrt{6t/a_2})^c$,
\begin{align*}
	p(t,x,w)=
	p^{\overline B(0,\eps)^c}(t,x,w) + \E_x[ p(t- \tau_{\overline B(0,\eps)^c} , X_{ \tau_{\overline B(0,\eps)^c} },w):\tau_{\overline B(0,\eps)^c}<t]=:I_1+I_2.
\end{align*}
For $I_1$, using \eqref{e:criticial-upper-1-0} and \eqref{e:criticial-upper-1-1}, we get
\begin{align}\label{e:criticial-upper-1-I1}
	I_1&\le  \frac{c_7 (\Lg a_1)(\Lg |x|  + \Log^{n+1} \sqrt t)}{t\,\Lg \sqrt t}  \exp \Big(- \frac{\lambda_1|w|^2}{4t}\Big).
\end{align}
For $I_2$, using $|w|\ge \sqrt{6t/a_2} \ge  2\sqrt 6$ and $(r-1)^2 \ge r^2 /2 \ge a_2r^2$ for all $r \ge 2\sqrt 6$, we see that
\begin{align*}
	& \sup_{s\in (0,1],\, z \in \overline B(0,\eps)\setminus \{0\} }p(s,z,w) \le  \sup_{s\in (0,1],\, z \in \overline B(0,1)\setminus \{0\} }q(s,z,w)\\
	& \le\frac{1}{4\pi}\sup_{s\in (0,1]} s^{-1} \exp  \Big( - \frac{(|w|-1)^2}{4s}\Big) \le \frac{\sqrt t \, \Lg |x| }{4\pi \log(e-1)	\Lg \sqrt t}\sup_{s\in (0,1]}s^{-3/2} \exp  \Big( - \frac{a_2|w|^2}{4s} \Big).
\end{align*}
where we used \eqref{e:Log-prop-2} and the fact that $\Lg |x| \ge \log(e-1)$ in the last inequality. 
Further, using \eqref{e:criticial-upper-1-ass}  in the first inequality below,  \eqref{e:Log-prop-2} and \eqref{e:criticial-upper-1-2} in the second,  we obtain
\begin{align*}
	& \sup_{s\in [1,t],\, z \in \overline B(0,\eps)\setminus \{0\} }p(s,z,w) \le    \sup_{s\in [1,t],\, z \in \overline B(0,\eps)\setminus \{0\} }	\frac{a_1  h(|z|)\, \Log^n \sqrt s }{s\,\Lg \sqrt s}  \exp \Big(- \frac{a_2|w|^2}{4s}\Big)\nn\\
	&\le   \frac{c_1a_1 \eps \sqrt t\, \Log^n \sqrt t}{\Lg \sqrt t}  \sup_{s\in [1,t]} s^{-3/2}  \exp \Big(- \frac{a_2|w|^2}{4s}\Big)
	\le   \frac{c_1 \sqrt t \, \Log |x|}{\Lg \sqrt t}  \sup_{s\in [1,t]} s^{-3/2}  \exp \Big(- \frac{a_2|w|^2}{4s}\Big).
\end{align*}
For any $\lambda \ge 3t/2$, the map $s\mapsto s^{-3/2} e^{-\lambda/s}$  is increasing in $(0,t]$. It follows that
\begin{align*}
	I_2 &\le \sup_{s\in (0,t],\, z \in \overline B(0,\eps)\setminus \{0\} }p(s,z,w) \le \frac{c_8 \Lg |x|}{t\,\Lg \sqrt t}   \exp \Big(- \frac{a_2|w|^2}{4t}\Big).
\end{align*}
Combining this with \eqref{e:criticial-upper-1-I1}, we deduce that for a.e. $w \in B(0,\sqrt{6t/a_1})^c$,
\begin{align*}
	p(t,x,w) \le \frac{c_9(\Lg a_1) (\Lg |x|  + \Log^{n+1} \sqrt t)}{t\,\Lg \sqrt t}  \exp \Big(- \frac{(\lambda_1 \wedge a_2)|w|^2}{4t}\Big).
\end{align*}
By the lower semi-continuity of $p$, this completes the proof. \qed

\begin{lem}\label{l:criticial-upper-2}
	There exists  $C\ge 1$ 
	such that the following holds:  If there exist $n \in \N$, $a_1\ge 1$ and $a_2 \in (0,1/2]$  such that for all $t\ge 2$ and $x,y \in \R^2_0$ with $|x|<\sqrt{2t}$, 
	\begin{equation}\label{e:criticial-upper-2-ass}
		p(t,x,y) \le 
		\frac{ a_1(\Lg |x| + \Log^{n} \sqrt t ) }{t\,\Lg \sqrt t} \exp \Big(- \frac{a_2|y|^2}{4t}\Big),
	\end{equation}
	then for all $t\ge 1$ 
	  and $z,y \in \R^2_0$ with $|z|\le 1$, 
	\begin{equation}\label{e:criticial-upper-2-result}
		p(t,z,y) \le 
		\frac{C a_1  h(|z|)\, \Log^n \sqrt t }{t\,\Lg \sqrt t}  \exp \Big(- \frac{(\lambda_0 \wedge a_2)|y|^2}{4t}\Big),
	\end{equation}
	where $\lambda_0 \in (0,1/2]$ is the constant in Lemma \ref{l:UHK-general-estimate}.
\end{lem}
\pf  Let $t\ge 1$ and $z,y \in \R^2_0$ with $|z|\le 1$. When $t\le 4$, since $1\le \Log^n \sqrt t \le  \Lg \sqrt t \le \Lg 2$,  using Lemma \ref{l:UHK-general-estimate},   we get
\begin{align*}
	p(t,z,y) &\le  	\frac{c_1  h(|z|)}{t}  \exp \Big(- \frac{\lambda_0 |y|^2}{4t}\Big) \le 	\frac{c_2 h(|z|)\, \Log^n \sqrt t}{t\, \Log \sqrt t}  \exp \Big(- \frac{\lambda_0 |y|^2}{4t}\Big).
\end{align*}
Thus, by taking $C$ larger than $c_2$, \eqref{e:criticial-upper-2-result} holds.

Assume $t \ge 4$. 
By the semigroup property and Lemma \ref{l:UHK-general-estimate}, we have
\begin{align*}
	p(t,z,y)&=  \int_{\R^2_0} p(2^{-1}\Log^n  \sqrt t,z,w)\,p(t- 2^{-1}\Log^n \sqrt t,w,y) \,dw\nn\\
	&\le   c_3h(|z|) \bigg( \int_{B(0, \sqrt t)}  + \int_{B(0, \sqrt t)^c}\bigg)   q(\frac{\Log^n \sqrt t}{2\lambda_0}, 0, w)  \, p(t- 2^{-1}\Log^n \sqrt t,w,y)\, dw\nn\\
	&=:c_3 h(|z|)  (I_{1}+I_{2}).
\end{align*}
By  \eqref{e:Log-prop-2},  we get  $1\le \Log^n  \sqrt t  \le \sqrt t \le  
t/2$. Using this, \eqref{e:criticial-upper-2-ass} and \eqref{e:Log-prop-2}, we obtain
\begin{align}\label{e:criticial-upper-2-3}
	I_1	\le  \frac{(4/3)^{3/2} a_1}{ t\,\Lg \sqrt{t}}\exp \Big(- \frac{a_2|y|^2}{4t}\Big)  \int_{B(0, \sqrt t)}   q(\frac{\Log^n \sqrt t}{2\lambda_0}, 0, w)  (\Lg |w| + \Log^n \sqrt t) \, dw.
\end{align}
Note that
\begin{align}\label{e:criticial-upper-2-4}
	&\int_{B(0, \Log^{n-1} \sqrt t)}   q(\frac{\Log^n \sqrt t}{2\lambda_0}, 0, w)  (\Lg |w| + \Log^n \sqrt t) \, dw \nn\\
	&\le  2\Log^n \sqrt t \int_{B(0, \Log^{n-1} \sqrt t)}   q(\frac{\Log^n \sqrt t}{2\lambda_0}, 0, w)  dw \le 2\Log^n \sqrt t.
\end{align}
Further, using \eqref{e:q-q2} in the first inequality below,  \eqref{e:Log-prop-2} in the second and  $\Log^{n-1} \sqrt t\ge \Log^n \sqrt t \ge 1$ in the third, we obtain
\begin{align*}
	&\int_{B(0, \Log^{n-1} \sqrt t)^c}   q(\frac{\Log^n \sqrt t}{2\lambda_0}, 0, w)  (\Lg |w| + \Log^n \sqrt t) \, dw  
	\le \frac{c_4\Log^n \sqrt t}{\lambda_0} \int_{B(0, \Log^{n-1} \sqrt t)^c} \frac{\Lg |w|}{|w|^4}  dw\nn\\
	&\le \frac{c_4 (\Log^n \sqrt t)^2}{\lambda_0 \Log^{n-1}\sqrt t } \int_{B(0, \Log^{n-1} \sqrt t)^c} \frac{dw}{|w|^{3}}  = \frac{c_5 (\Log^n \sqrt t)^2}{\lambda_0(\Log^{n-1}\sqrt t)^{2}} \le  \frac{c_5 \Log^n \sqrt t}{\lambda_0}.
\end{align*}
Combining this, \eqref{e:criticial-upper-2-3} and \eqref{e:criticial-upper-2-4}, we get
\begin{align}\label{e:criticial-upper-I1}
	I_1&\le  \frac{c_6  a_1 \Log^n \sqrt t}{t\,\Lg \sqrt{t}}\exp \Big(- \frac{a_2|y|^2}{4t}\Big) .
\end{align} 
On the other hand, note that
\begin{align*}
	\frac{s^2}{a} + \frac{r^2}{b} \ge \frac{(s+r)^2}{a+b} \quad \text{for all $a,b>0$ and $s,r \ge 0$}.
\end{align*}
Using this,  $t-2^{-1}\Log^n \sqrt t \ge 3t/4$ and $\lambda_0\le 1/2$, we see that
\begin{align*}
	I_2 &\le
	\int_{B(0, \sqrt t)^c} q(\frac{\Log^n \sqrt t}{2\lambda_0}, 0, w)  \, q(t-2^{-1} \Log^n \sqrt t,w,y)\, dw\nn\\
	&\le \frac{c_{7}}{t\, \Log^n \sqrt t}  \int_{B(0, \sqrt t)^c}   \exp \Big(- \frac{\lambda_0 |w|^2}{2\Log^n \sqrt t} - \frac{\lambda_0|w-y|^2}{4(t-2^{-1} \Log^n \sqrt t)} \Big) dw \nn\\
	&\le \frac{c_{7}}{t\, \Log^n \sqrt t}  \exp \Big( - \frac{\lambda_0|y|^2}{4t} \Big)\int_{B(0, \sqrt t)^c}  \exp \Big( - \frac{3\lambda_0|w|^2}{8\Log^n \sqrt t} \Big)   dw.
\end{align*}
By \eqref{e:exp-poly}, we have
\begin{align*}
	\int_{B(0, \sqrt t)^c}  \exp \Big( - \frac{3\lambda_0|w|^2}{8\Log^n \sqrt t} \Big)   dw \le c_{8} \Big(  \frac{8\Log^n \sqrt t}{3\lambda_0} \Big)^2 	\int_{B(0, \sqrt t)^c}  \frac{dw}{|w|^4}  = \frac{c_{9} (\Log^n \sqrt t)^2}{t}  .
\end{align*}
Thus, since $\Log \sqrt t \le t/2$, we deduce that
\begin{align}\label{e:criticial-upper-I2}
	I_2 \le \frac{c_{10}\Log^n \sqrt t}{t^2}  \exp \Big( - \frac{\lambda_0|y|^2}{4t} \Big) \le \frac{c_{10}\Log^n \sqrt t}{2t \,\Log \sqrt t}  \exp \Big( - \frac{\lambda_0|y|^2}{4t} \Big).
\end{align}
Combining \eqref{e:criticial-upper-I1} and \eqref{e:criticial-upper-I2}, we arrive at the desired result. \qed

\begin{lem}\label{l:largetime-3-upper}
 There exist constants $C,k\ge 1$ such that for all $t\ge 2$ and  $x,y \in \R^2_0$,
\begin{align}\label{e:largetime-3-upper}
	p(t,x,y) \le C H_{2}(t,|x|) q(kt,x,y).
\end{align}
\end{lem}
\pf 
Set $\lambda_2:= \lambda_0 \wedge \lambda_1$, where  $\lambda_0$ and $\lambda_1$ are  the constants from Lemmas \ref{l:UHK-general-estimate} and \ref{l:criticial-upper-1} respectively. Since Lemma \ref{l:UHK-general-estimate} implies \eqref{e:criticial-upper-1-ass} with $a_2=\lambda_2$  for $n=1$,
 by Lemma \ref{l:criticial-upper-1}, there exist $c_1,c_2\ge1$ such that	 for all $s\ge 2$ and $v,w \in \R^2_0$ with $|v| \le \sqrt {2s}$,
\begin{equation*}
	p(s,v,w) \le 
	\frac{c_1 (\Lg c_2)(\Lg |v| + \Log^{2} \sqrt s ) }{s\,\Lg \sqrt s} \exp \Big(- \frac{\lambda_2|w|^2}{4s}\Big),
\end{equation*}
implying \eqref{e:criticial-upper-2-ass} with $a_2=\lambda_2$ for $n=2$. Thus, by Lemma \ref{l:criticial-upper-2}, there exists $c_3 \ge 1$ such that for all $s \ge 1$ and $z,w\in \R^2_0$ with $|z|\le 1$,
\begin{align*}
	p(s,z,w) &\le \frac{c_1  c_3 (\Lg c_2) h(|z|) \Log^{2} \sqrt s}{s \,\Lg \sqrt s} \exp \Big( - \frac{\lambda_2|w|^2}{4s}\Big).
\end{align*}By iterating this procedure, we deduce that the following inequalities hold  for all  $n \ge 1$:
\begin{equation}\label{e:largetime-3-upper-iteration-1}
	p(s,v,w) \le 
\frac{b_n
	(\Lg |v| + \Log^{n} \sqrt s ) }{s\,\Lg \sqrt s} \exp \Big(- \frac{\lambda_2|w|^2}{4s}\Big) 
\end{equation}
for all $s\ge 2$ and  $v,w \in \R^2_0$ with $|v|\le \sqrt {2s}$, and
\begin{equation}\label{e:largetime-3-upper-iteration-2}
		p(s,z,w) \le \frac{d_n
			h(|z|) \Log^{n} \sqrt s}{s \,\Lg \sqrt s} \exp \Big( - \frac{\lambda_2|w|^2}{4s}\Big)
\end{equation}
for all $s\ge 1$ and  $z,w \in \R^2_0$ with $|z|\le 1$, 
where the sequences $(b_n)_{n\ge 1}$ and  $(d_n)_{n\ge 1}$ are inductively defined as 
$$b_{1}:=c_1 \Lg c_2, \quad d_n:= c_3 b_n\quad \text{and} \quad b_{n+1}:=c_1 \Lg d_n. $$
By \eqref{e:Log-scaling} with $\eps=1/2$, there is $c_4\ge 1$ such that $\Lg r \le c_4\sqrt r$ for all $r \ge 1$. It follows that 
\begin{align*}
	b_{n+1} = c_1 \Lg (c_3 b_n) \le c_1 c_3^{1/2}c_4 \sqrt b_n \quad \text{for all $n \ge 1$}.
\end{align*}
Thus,  $\limsup_{n\to \infty} b_n \le c_1^2 c_3c_4^2=:c_5$. 
 Using this and the fact that  $\lim_{n\to \infty} \Log^{n} \sqrt s = 1$ for all $s \ge 1$, by taking the limit superior  $\limsup_{n \to \infty}$ in \eqref{e:largetime-3-upper-iteration-1} and \eqref{e:largetime-3-upper-iteration-2},  we conclude that  
 \begin{equation}\label{e:largetime-3-upper-conclusion-1}
 p(s,v,w) \le 	\frac{c_5(\Lg |v| +1 ) }{s\,\Lg \sqrt s} \exp \Big(- \frac{\lambda_2|w|^2}{4s}\Big)  
 \end{equation}
 for all $s\ge 2$ and $v,w \in \R^2_0$ with $|v|\le \sqrt {2s}$, and
 \begin{equation}\label{e:largetime-3-upper-conclusion-2}
 p(s,z,w) \le \frac{c_3 c_5h(|z|) }{s \,\Lg \sqrt s} \exp \Big( - \frac{\lambda_2|w|^2}{4s}\Big)
 \end{equation}
 for all $s\ge 1$ and $z,w \in \R^2_0$ with $|z|\le 1$.

Set $k:=2/\lambda_2$. Let   $t\ge 2$ and $x,y \in \R^2_0$.  If $|x| >\sqrt {2t}$, then since $H_{2}(t,|x|)=1 $, \eqref{e:largetime-3-upper} follows
 from the fact that $p(t,x,y) \le q(t,x,y) \le k^{d/2} q(k t,x,y)$. If $1 < |x|\le \sqrt {2t}$, then  using \eqref{e:largetime-3-upper-conclusion-1}, 
 $\Lg |x| \ge 1$ and \eqref{e:q-estimate-1}, we get that
\begin{align*}
	p(t,x,y) &\le 	\frac{4\pi k c_5\Lg |x|  }{\Lg \sqrt t} q(kt/2, 0,y) \le c_6 H_{2}(t,|x|) q(kt, x,y).
\end{align*}
 If $|x| <1$, then by \eqref{e:largetime-3-upper-conclusion-2} and \eqref{e:q-estimate-1}, we obtain
 \begin{align*}
 	p(t,x,y) \le \frac{2\pi k c_3c_5 h(|x|) }{\Lg \sqrt t} q(kt/2, 0,y) \le c_7 H_{2}(t,|x|)q(kt,x,y).
 \end{align*}
 The proof is complete. \qed

\noindent \textbf{Proof of Theorem \ref{t:largetime} (Upper estimates).}  If $d\ge 3$, then  the  upper bound follows  from  Corollary \ref{c:UHK-general}.  Assume $d=1$ or $d=2$.  Using the semigroup property and symmetry  in the first line below, and  Lemma \ref{l:largetime-2-upper} if $d=1$ and Lemma \ref{l:largetime-3-upper} if $d=2$ in the second, we get that  for all $t \ge 4$ and $x,y \in \R^d_0$,
\begin{align*}
	p(t,x,y) &= \int_{\R^d_0} 	p(t/2,x,z)	p (t/2,y,z) dz\\
	&\le c_1^2 H_{d}(t/2,|x|) H_{d}(t/2,|y|)\int_{\R^d_0} 	q(kt/2,x,z)	 q(kt/2,y,z) dz \\
	&= c_1^2 H_{d}(t/2,|x|) H_{d}(t/2,|y|) q(kt,x,y).
\end{align*}
Since $H_{d}(s/2,|z|)\le 2 H_{d}(s,|z|)$ for all $s>0$ and $z \in \R^d_0$, this completes the proof.  \qed

\section{Green function estimates}\label{s:7}

 Let
\begin{align*}
G(x,y):=\int_0^\infty p(t,x,y)dt, \quad x,y \in \R^d_0.
\end{align*} 
Define for $x,y \in \R^d$,
\begin{align*}
g_0(x,y):=	\begin{cases}
		|x-y|^{2-d} &\mbox{ if $d\ge 3$},\\[2pt]
		\displaystyle \Lg \Big( \frac{|x| \wedge |y|}{|x-y| \wedge1}\Big)  &\mbox{ if $d= 2$},\\[7pt]
		|x| \wedge |y|&\mbox{ if $d= 1$}
	\end{cases}
\end{align*}
and
\begin{align*}
	f_0(x,y):= \begin{cases}
		|x-y|^{2-d} &\mbox{ if $d\ge 3$},\\
		\displaystyle	\Lg \Big( \frac{\big( (|x| \vee |y|) \wedge 1\big)^{1+\beta}}{|x-y|} \Big)  &\mbox{ if $d= 2$},\\[7pt]
	|x-y|	\vee \big( (|x| \vee |y|) \wedge 1\big)^{1+\beta}  &\mbox{ if $d= 1$}.
	\end{cases}
\end{align*}
Note that  in any dimension $d$, we have
\begin{align}\label{e:f0-lower}	f_0(x,y) \ge    \log(e-1)\,|x-y|^{2-d}  \quad \text{for all $x,y \in \R^d$}.\end{align}

In this section, we  prove the following two-sided Green function estimates.

\begin{thm}\label{t:Green}
Suppose   $V \in \sK(\beta,\kappa)$.  Then there exist constants $\eta_2>0$, $c_1\ge c_2>0$ and $C>1$ such that   the following estimates hold:

\smallskip

\noindent (i) Assume $d\ge 2$. For all $x,y \in \R^d_0$,  if $|x| \wedge |y| \le 2$, then 
\begin{align*}
 &C^{-1}\Big( 1\wedge \frac{h_{\beta,\kappa}(|x|\wedge |y|)}{h_{\beta,\kappa}(\eta_2|x-y|)}\Big) f_0(x,y) \exp \Big( - \frac{c_1|x-y|}{(|x| \vee |y|)^{1+\beta}}\Big)\\
 &\le 		G(x,y)  \le  C\Big( 1\wedge \frac{h_{\beta,\kappa}(|x|\wedge |y|)}{h_{\beta,\kappa}(\eta_2|x-y|)}\Big) f_0(x,y) \exp \Big( - \frac{c_2|x-y|}{(|x| \vee |y|)^{1+\beta}}\Big),
\end{align*}
and if $|x| \wedge |y| \ge 2$, then
\begin{align*}
C^{-1}g_0(x,y)\le 	G(x,y) \le Cg_0(x,y).
\end{align*}

\noindent (ii) If $d=1$, then the estimates in (i) hold for all $x,y \in \R^1_0$ with $xy>0$.
\end{thm}

The proof of Theorem \ref{t:Green} will be given at the end of this section.

\begin{remark}\label{r:Green}
\rm For all $|x| \wedge |y| \ge 2$, we have $h_{\beta,\kappa}(|x| \wedge |y|) \asymp1$ and $(|x| \vee |y|)^{1+\beta} \ge 2^\beta (|x| \vee |y|) \ge 2^{\beta-1}|x-y|$. Using these, we deduce that for $d\ge 3$, the Green function estimates  in Theorem \ref{t:Green} 
can be simplified to the  following form: For all $x,y \in \R^d_0$,
\begin{align*}
	&C^{-1}\Big( 1\wedge \frac{h_{\beta,\kappa}(|x|\wedge |y|)}{h_{\beta,\kappa}(\eta_2|x-y|)}\Big) |x-y|^{2-d} \exp \Big( - \frac{c_1|x-y|}{(|x| \vee |y|)^{1+\beta}}\Big)\\
	&\le 		G(x,y)  \le  C\Big( 1\wedge \frac{h_{\beta,\kappa}(|x|\wedge |y|)}{h_{\beta,\kappa}(\eta_2|x-y|)}\Big) |x-y|^{2-d}\exp \Big( - \frac{c_2|x-y|}{(|x| \vee |y|)^{1+\beta}}\Big).
\end{align*}
\end{remark}

For an open subset $U\subset \R^d$, define 
$$g^U(x,y):=\int_0^\infty q^U(t,x,y)dt, \quad x,y \in \R^d.$$
We begin with  two-sided estimates for $g^{\overline B(0,1)^c}(x,y)$. 
Note that  by \eqref{e:exp-poly}, for any $a>0$, there exists $C=C(a)>0$ such that for aall $r,R>0$,
\begin{align}\label{e:Green-elementary}
	\int_0^{r} t^{-a/2} e^{-R/t} dt \le \frac{C}{(R/2)^{a/2}} \int_0^{r} e^{-R/(2t)} dt  \le  \frac{C r e^{-R/(2r)}}{(R/2)^{a/2}}.
\end{align}

\begin{lem}\label{l:Green-exterior}
	(i) Assume $d\ge 2$.	There   exist comparison constants  such that  for all $x,y \in \R^d_0$ with $|x| \wedge  |y|  \ge 2$,
	\begin{align}\label{e:Green-exterior}
	g^{\overline B(0,1)^c}(x,y)   \asymp g_0(x,y).
	\end{align}
	(ii) If $d=1$, then  \eqref{e:Green-exterior} holds for all $x,y \in \R^1_0$ with $xy>0$ and $|x| \wedge  |y|  \ge 2$. 
\end{lem}
\pf Let $x,y \in \R^d_0$ with $|x| \wedge |y| \ge 2$. Without loss of generality, we assume $|x|\le |y|$.  When $d=1$, we 
also assume that $xy>0$.   We consider the following three cases separately.

\smallskip

\noindent
{\it
Case 1: $d\ge 3$.} We have $	g^{B(0,1)^c}(x,y)  \le g^{\R^d}(x,y)  =c_1 |x-y|^{2-d}$. On the other hand, using  Proposition \ref{p:DHKE}, we obtain
\begin{align*}
		g^{\overline B(0,1)^c}(x,y)  \ge \int_{|x-y|^2}^{2|x-y|^2} q^{\overline B(0,1)^c}(t,x,y)dt \ge     c_2\int_{|x-y|^2}^{2|x-y|^2}  t^{-d/2} dt  = c_3 |x-y|^{2-d}.
\end{align*}

\noindent
{\it
Case 2: $d=2$.}  By Proposition \ref{p:DHKE} and \eqref{e:Log-prop-0},  we have $g^{\overline B(0,1)^c}(x,y)  \ge c_4 \int_{4|y|^2}^\infty  \frac{\Lg |x|  \log |y|}{t(\log \sqrt t)^2} dt =  c_5\Lg |x|.$
Further, when $|x-y|<1$, using \eqref{e:Log-prop-0}, we see that
\begin{align*}
			g^{\overline B(0,1)^c}(x,y)  \ge c_6	\int_{|x-y|^2}^4 \frac{dt}{t} =2c_6 \log \Big(\frac{2}{|x-y|}\Big) \ge c_7\Lg \Big(\frac{1}{|x-y|}\Big).
\end{align*}
Hence, using the inequality $a\vee b \ge (a+b)/2$  for all $a,b\ge 0$ and  \eqref{e:Log-prop-1}, we obtain
\begin{align*}
	&		g^{\overline B(0,1)^c}(x,y)  \ge   \frac{c_5}{2}\Log |x| + \frac{c_7}{2}\1_{\{|x-y|<1\}} \Lg \Big(\frac{1}{|x-y|}\Big) \ge c_8 \Lg \Big( \frac{|x|}{|x-y| \wedge 1}\Big) .
\end{align*}
For the upper bound,  using Proposition \ref{p:DHKE} and  \eqref{e:Green-elementary}, we get
\begin{align*}
&	g^{\overline B(0,1)^c}(x,y) \le  \int_0^{4|y|^2} q(t,x,y)dt + \int_{4|y|^2}^\infty q^{\overline B(0,1)^c}(t,x,y) dt\\
&\le  \frac{1}{4\pi}\int_0^{|x-y|^2}  \frac{1}{ t} \exp \Big(- \frac{|x-y|^2}{4t}\Big) dt  + \frac{1}{4\pi}\int_{|x-y|^2}^{4|y|^2}  \frac{1}{t}dt +  c_{9} \int_{4|y|^2}^\infty  \frac{\log |x| \log |y|}{t(\log \sqrt t)^2} dt \nn\\
	&\le  c_{10} +  \frac{1}{2\pi}\log \Big( \frac{2|y|}{|x-y|}\Big)  +  \frac{2c_{9}\log |x| \log|y|}{\log (2|y|)}\le 	c_{11} \Big( 1 + \log \Big( \frac{2|y|}{|x-y|}\Big)  + \log |x| \Big).
\end{align*}
Since  $ 1 \vee \log |x| \le \Lg |x|$ and  
\begin{align*}
	\log \Big( \frac{2|y|}{|x-y|}\Big) & \le 	\log \Big( 2+\frac{2|x|}{|x-y|}  \Big)   \le 	\log \Big( (e-1)^2+\frac{2(e-1)|x|}{|x-y|\wedge 1}  \Big)  \le  2 \Lg  \Big ( \frac{|x|}{|x-y| \wedge 1}\Big),
\end{align*}
we arrive at the desired result in this case.

\smallskip

\noindent
{\it
Case 3: $d=1$.}
By Proposition \ref{p:DHKE}, we obtain $		g^{\overline B(0,1)^c}(x,y)  \ge c_{12} |x|   |y| \int_{4|y|^2}^\infty  t^{-3/2} dt =  c_{13}|x|$. 
For the upper bound, if $|x| \ge |x-y|$, using Proposition \ref{p:DHKE} and the fact that  $|y| \le |x| +|x-y|\le 2|x|$, we get
\begin{align*}
		g^{\overline B(0,1)^c}(x,y)  &\le \int_0^{|y|^2} q(t,x,y) dt  + \int_{|y|^2}^\infty q^{\overline B(0,1)^c}(t,x,y)dt\\
	&\le  \frac{1}{4\pi} \int_0^{|y|^2} \frac{1}{t^{1/2}}dt  + c_{14}\int_{|y|^2}^\infty  \frac{ |x|  |y|}{t^{3/2}} dt    = \frac{|y|}{2\pi} + 2c_{14} |x| \le \Big(\frac{1}{\pi} + 2c_{14}\Big) |x|.
\end{align*}
If $|x|<|x-y|$, then using Proposition \ref{p:DHKE}, $|y| \le |x-y| + |x| \le 2|x-y|$ and  \eqref{e:Green-elementary}, we get
\begin{align*}
		g^{\overline B(0,1)^c}(x,y)  &\le   c_{15}\int_0^{|x-y|^2}  \frac{|x|}{ t} \exp \Big(- \frac{c_{16}|x-y|^2}{t}\Big) dt  + c_{17} \int_{|x-y|^2}^\infty  \frac{ |x|  |y|}{t^{3/2}} dt \nn\\
	&\le  c_{18}|x|  + \frac{c_{19}|x||y|}{|x-y|} \le c_{20}|x|.
\end{align*}

The proof is complete.\qed 

In the remainder of this section, we assume that  $V \in \sK(\beta,\kappa)$ and write $h$  instead of $h_{\beta,\kappa}$.

\begin{lem}\label{l:Green-boundary-0}
	There   exist comparison constants  such that  for all $x,y \in \R^d_0$ with $|x| \wedge  |y|  \ge 2$,
	\begin{align*}
		G(x,y)   \asymp 	g^{\overline B(0,1)^c}(x,y) .
	\end{align*}
\end{lem}
\pf  By Theorems \ref{t:smalltime} and \ref{t:largetime}, \eqref{e:compare-H-psi} and Proposition \ref{p:DHKE}, we get
\begin{align*}
	G(x,y) &\le c_1\bigg( \int_0^4 q(c_2t,x,y)dt  +  \int_4^\infty \psi_{d,1/2}(t,x)\psi_{d,1/2}(t,x) q(c_2t,x,y) dt \bigg) \\
	&\le c_3 \int_0^\infty q^{\overline B(0,1)^c}(c_4t,x,y)dt =  c_3c_4^{-1} g^{\overline B(0,1)^c}(x,y) 
\end{align*}
and $G(x,y) \ge c_5 \int_0^\infty q^{\overline B(0,1)^c}(c_6t,x,y)dt = c_5 c_6^{-1}g^{\overline B(0,1)^c}(x,y)$.  \qed

Lemmas \ref{l:Green-exterior} and \ref{l:Green-boundary-0}
give the two-sided Green function estimates for $|x| \wedge  |y|  \ge 2$.

We will use the following elementary lemma. The proof is straightforward.
\begin{lem}\label{l:Green-elementary}
	For any $a,b>0$, there exists a constant $c=c(a/b,\beta)>0$  such that 
	\begin{align*}
		\sup_{t>0} \bigg( - \frac{ar^2}{t}  + \frac{b}{t^{\beta/(2+\beta)}}\bigg) = - \frac{a}{cr^\beta}  + \frac{b}{c^{\beta/(2+\beta)} r^\beta}.
	\end{align*}
\end{lem}

In the next two lemmas, we establish the two-sided Green function estimates for $(|x| \wedge  |y|) \vee |x-y|  \le 2$.

\begin{lem}\label{l:Green-boundary-1-upper}
	There   exist constants  $\eta_2,c_0>0$ and $C>0$  such that  for all $x,y \in \R^d_0$ with $(|x| \wedge  |y|) \vee |x-y|  \le 2$,
	\begin{align*}
		G(x,y)  \le  C\Big( 1\wedge \frac{h(|x|\wedge |y|)}{h(\eta_2|x-y|)}\Big) f_0(x,y) \exp \Big( - \frac{c_0|x-y|}{(|x| \vee |y|)^{1+\beta}}\Big).
	\end{align*}
\end{lem}
\pf We assume, without loss of generality, that  $|x|\le |y|$.  Note that $|y| \le 4$. By Theorems \ref{t:smalltime} and \ref{t:largetime}, 
 \eqref{e:Log-prop-0} and 
 \eqref{e:H-arbitrary-R} (if $d=1$ or $d=2$),  we have
\begin{align*}
	\frac{G(x,y)}{h(|x|)} &\le  c_1  \int_{0}^{4}    \frac{e^{-c_2|x-y|^2/t}}{t^{d/2} h(\eta_1t^{1/(2+\beta)})}  dt  + c_1   h(|y|) \int_{4}^\infty \Big( \frac{\1_{\{d\ge 3\}}}{t^{d/2}}  + \frac{\1_{\{d=2\}}}{t (\log t)^2} + \frac{\1_{\{d=1\}}}{t^{3/2}} \Big) dt\\
	&=:c_1(I_1+I_2).
\end{align*}
  By Lemma \ref{l:Green-elementary}, there exists $c_3>0$ such that
\begin{align*}
	\sup_{t>0}   \Big(- \frac{c_2|x-y|^2}{2t} +\frac{\sqrt \kappa}{ \beta \eta_1^\beta t^{\beta/(2+\beta)}} \Big)  = - \frac{c_2}{2c_3|x-y|^\beta} +\frac{\sqrt \kappa}{ \beta  (\eta_1c_3^{1/(2+\beta)}|x-y|)^{\beta} }.
\end{align*}
Using this, \eqref{e:Green-elementary} and \eqref{e:exp-poly} (if $\beta(d-2-\beta)/(4+2\beta) + d>0$), we get
\begin{align}\label{e:Green-1-I1}
	I_1&\le  \frac{c_4|x-y|^{-(d-2-\beta)/2}}{ h(\eta_1c_3^{1/(2+\beta)}|x-y|)}   \exp \Big(  - \frac{c_2}{2c_3|x-y|^\beta}\Big) \int_{0}^{4}    t^{(d-2-\beta)/(4+2\beta) -d/2}  \exp \Big(- \frac{c_2|x-y|^2}{2t} \Big) dt\nn\\
	&\le  \frac{c_5 }{|x-y|^{\beta(d-2-\beta)/(4+2\beta) +d} \, h(\eta_1c_3^{1/(2+\beta)}|x-y|)} \exp \Big(- \frac{c_2}{2c_3|x-y|^\beta} - \frac{c_2|x-y|^2}{16}\Big)   \nn\\
	&\le \frac{c_6}{ h(\eta_1c_3^{1/(2+\beta)}|x-y|)} \exp \Big(- \frac{c_2}{3c_3|x-y|^\beta}\Big).
\end{align}
Set $\eta_2:=\eta_1c_3^{1/(2+\beta)}$.  We have $h(\eta_2 |x-y|) \le c_7 h(2\eta_2)$.  Thus, using Lemma \ref{l:h-ratio}, we obtain
\begin{align}\label{e:Green-1-I2}
	I_2 \le 	c_8h(|y|)\le  \frac{c_{9}}{h(\eta_2|x-y|)} \exp \Big(- \frac{\sqrt \kappa}{2\beta |y|^\beta} \Big).
\end{align} 
Combining \eqref{e:Green-1-I1} and \eqref{e:Green-1-I2}, and  using $|x-y| \le 2|y|$, we deduce that
\begin{align}\label{e:Green-1-case1}
	G(x,y) &\le \frac{c_6 h(|x|)}{ h(\eta_2|x-y|)} \exp \Big(- \frac{c_2}{3c_3|x-y|^\beta}\Big) + \frac{c_9 h(|x|)}{h(\eta_2|x-y|)} \exp \Big(- \frac{\sqrt \kappa}{2\beta |y|^\beta} \Big) \nn\\
	&\le
	\frac{c_{10} h(|x|)}{ h(\eta_2|x-y|)} \exp \Big(- \Big(  \frac{c_2}{3\cdot 2^{1+\beta} c_3} \wedge \frac{\sqrt \kappa}{4\beta}\Big) \frac{|x-y|}{|y|^{1+\beta}}\Big).
\end{align}
Assume $|x| < \eta_2|x-y|$. Then $|y|< (1+\eta_2)|x-y|$. Hence, by \eqref{e:f0-lower} and \eqref{e:exp-poly},  it holds that
\begin{align*}
	f_0(x,y) & \ge \log(e-1)|x-y|^{2-d} \ge 2^{1-d}\log(e-1)|x-y|\\
	& \ge c_{11} \exp \Big(- \frac{1}{2(1+\eta_2)^{1+\beta}} \Big(  \frac{c_2}{3\cdot 2^{1+\beta} c_3} \wedge \frac{\sqrt \kappa}{4\beta}\Big) \frac{1}{|x-y|^{\beta}} \Big) \nn\\
	& \ge  c_{11}\exp \Big(- \frac{1}{2}\Big(  \frac{c_2}{3\cdot 2^{1+\beta} c_3} \wedge \frac{\sqrt \kappa}{4\beta}\Big) \frac{|x-y|}{|y|^{1+\beta}}\Big).
\end{align*}
Combining this with \eqref{e:Green-1-case1}, we get the desired result in this case.

We now assume that $|x| \ge \eta_2|x-y|$.  Using Theorems  \ref{t:smalltime} and  \ref{t:largetime},  we get
\begin{align*}
	G(x,y) &\le c_{12} \int_{0}^{(2|y|/\eta_1)^{2+\beta}}    t^{-d/2}  \exp \Big(- \frac{c_{13}t}{|y|^{2+2\beta}} - \frac{c_{14}|x-y|^2}{t}\Big)  dt  \\
	&\quad + c_{12} \int_{(2|y|/\eta_1)^{2+\beta}}^{(2|y|/\eta_1)^{2+\beta} \vee 4}   \frac{h(|y|) }{t^{d/2} h(\eta_1 t^{1/(2+\beta)})}  \exp \Big(-\frac{c_{13}}{t^{\beta/(2+\beta)}}\Big)  dt \\
		&\quad +  c_{12}   h(|y|) \int_{4}^\infty \Big( \frac{\1_{\{d\ge 3\}}}{t^{d/2}}  + \frac{\1_{\{d=2\}}}{t (\Lg t)^2} + \frac{\1_{\{d=1\}}}{t^2} \Big) dt\\
	&=:c_{12}(I_1'+I_2'+I_3').
\end{align*}
For $I_2'$, by Lemma \ref{l:h-ratio} and \eqref{e:exp-poly}, we see that for all $t\in [ (2|y|/\eta_1)^{2+\beta},4]$,
\begin{equation}\label{e:Green-1-case2}
	\frac{h(|y|)}{h(\eta_1t^{1/(2+\beta)})} \le  c_{15}\exp \Big (- \frac{ 9^\beta \sqrt \kappa}{10^\beta \beta|y|^\beta} +  \frac{3^\beta \sqrt \kappa}{5^\beta \beta|y|^\beta} \Big)  \le  c_{16}|y|^{1+\beta}\exp \Big (- \frac{ (4^\beta - 3^\beta )\sqrt \kappa}{5^\beta \beta|y|^\beta} \Big) .
\end{equation}
Further, since $|x-y|\le 2$ and $|y|\le 4$,  we have
\begin{align}\label{e:f0-lower2}
	f_0(x,y) \ge 2^{2-d} \1_{\{d\ge 3\}} + \log(e-1) \1_{\{d= 2\}} + (|y|\wedge 1)^{1+\beta} \1_{\{d= 1\}} \ge c_{17} |y|^{1+\beta}.
\end{align}
Using \eqref{e:Green-1-case2}, \eqref{e:f0-lower2} and $|y| \ge |x-y|/2$, we obtain
\begin{align*}
	I_2'& \le c_{18}f_0(x,y)\exp \Big (- \frac{ (4^\beta - 3^\beta )\sqrt \kappa}{5^\beta \beta|y|^\beta} \Big) \int_{0}^{4}   t^{-d/2} \exp \Big(- \frac{c_{13}}{t^{\beta/(2+\beta)}} \Big)   dt \nn\\
	& \le c_{19} f_0(x,y)\exp \Big (- \frac{ (4^\beta - 3^\beta )\sqrt \kappa \,|x-y| }{2\cdot 5^\beta \beta|y|^{1+\beta}} \Big) .
\end{align*}
For $I_3'$,  by Lemma \ref{l:h-ratio} and \eqref{e:exp-poly}, we have
\begin{align*}
	I_3' &\le c_{20} h(|y|)\le c_{21}|y|^{1+\beta}\exp \Big (- \frac{ 4^\beta \sqrt \kappa  }{5^\beta \beta|y|^{\beta}} \Big) \le c_{22} f_0(x,y)\exp \Big (- \frac{ 4^\beta \sqrt \kappa \,|x-y| }{2\cdot 5^\beta \beta|y|^{1+\beta}} \Big).
\end{align*}
For $I_1'$, we consider the cases $|x-y| \ge |y|^{1+\beta}$ and $|x-y|<  |y|^{1+\beta}$ separately. Suppose that $|x-y| \ge |y|^{1+\beta}$. Using \eqref{e:exp-poly} in the first and the fourth inequalities below, and \eqref{e:Green-elementary} in the second, we obtain
\begin{align*}
	I_1' &\le   \int_{0}^{(2/\eta_1)^{2+\beta}|y|^{1+\beta}|x-y|/2}    t^{-d/2}  \exp \Big( - \frac{c_{14}|x-y|^2}{t}\Big)  dt \nn\\
	&\quad + c_{23}\int_{(2/\eta_1)^{2+\beta}|y|^{1+\beta}|x-y|/2}^{(2/\eta_1)^{2+\beta}|y|^{2+\beta}}    \frac{|y|^{2+2\beta}}{t^{(d+2)/2}}  \exp \Big( - \frac{c_{13}t}{2|y|^{2+2\beta}}\Big)  dt \nn\\
	&\le \frac{c_{24} |y|^{1+\beta}}{|x-y|^{d-1}} \exp \Big( - \frac{c_{14}|x-y|}{(2/\eta_1)^{2+\beta} |y|^{1+\beta}}\Big)    \nn\\
	&\quad + c_{23}\exp \Big( - \frac{c_{13}(2/\eta_1)^{2+\beta}|x-y|}{4|y|^{1+\beta}}\Big) \int_{(2/\eta_1)^{2+\beta}|y|^{1+\beta}|x-y|/2}^{(2/\eta_1)^{2+\beta}|y|^{2+\beta}}    \frac{|y|^{2+2\beta}}{t^{1+d/2}}  dt\nn\\
		&\le \frac{c_{24} |y|^{1+\beta}}{|x-y|^{d-1}} \exp \Big( - \frac{c_{14}|x-y|}{(2/\eta_1)^{2+\beta} |y|^{1+\beta}}\Big)    + \frac{c_{25} |y|^{2+2\beta} }{(|y|^{1+\beta}|x-y|)^{d/2}}  \exp \Big( - \frac{c_{13}(2/\eta_1)^{2+\beta}|x-y|}{4|y|^{1+\beta}}\Big)\nn\\
			&\le \frac{c_{24} |y|^{1+\beta}}{|x-y|^{d-1}} \exp \Big( - \frac{c_{14}|x-y|}{(2/\eta_1)^{2+\beta} |y|^{1+\beta}}\Big)    + \frac{c_{26} |y|^{2+2\beta} }{|x-y|^d}  \exp \Big( - \frac{c_{13}(2/\eta_1)^{2+\beta}|x-y|}{5|y|^{1+\beta}}\Big).
\end{align*}
Since $|y|^{2+2\beta} |x-y|^{-d} \le |y|^{1+\beta} |x-y|^{1-d} \le |x-y|^{2-d}\le f_0(x,y)/\log(e-1)$ by  the assumption $|y|^{1+\beta}\le |x-y|$ and \eqref{e:f0-lower}, this yields the desired result.   Suppose  that $|x-y| < |y|^{1+\beta}$. Observe that
\begin{align*}
	\inf_{t>0} \Big( \frac{c_{13}t}{2|y|^{2+2\beta}} + \frac{c_{14}|x-y|^2}{2t}\Big)  = \frac{(c_{13}c_{14})^{1/2}|x-y|}{|y|^{1+\beta}}.
\end{align*}
Using this, we get
\begin{align*}
	I_1' & \le \Big[  \int_{0}^{(1/\eta_1)^{2+\beta} |x-y|^2}  t^{-d/2}  \exp \Big(- \frac{c_{14}|x-y|^2}{2t}\Big)  dt    +  \int_{(1/\eta_1)^{2+\beta} |x-y|^2} ^{(1/\eta_1)^{2+\beta}|y|^{2+2\beta}}    t^{-d/2}   dt      \nn\\
	&\qquad    +  \int_{(1/\eta_1)^{2+\beta} |y|^{2+2\beta}} ^{(2/\eta_1)^{2+\beta}|y|^{2+\beta} }    t^{-d/2}\exp \Big(- \frac{c_{13}t}{2|y|^{2+2\beta}}  \Big) dt   \Big]  \,	\exp \Big(- \frac{(c_{13}c_{14})^{1/2}|x-y|}{|y|^{1+\beta}}\Big) \nn\\
	&=:(I_{1,1}' + I_{1,2}'  + I_{1,3}')	\exp \Big(- \frac{(c_{13}c_{14})^{1/2}|x-y|}{|y|^{1+\beta}}\Big) .
\end{align*}
Using $|y|\le 4$, we see that $I_{1,2}' \le c_{27}f_0(x,y).$
By \eqref{e:Green-elementary} and \eqref{e:f0-lower}, we have $I_{1,1}' \le c_{28} |x-y|^{d-2} \le c_{29}f_0(x,y).$ Further, using \eqref{e:exp-poly} and \eqref{e:f0-lower} if $d=2$, we obtain
\begin{align*}
I_{1,3}'&\le   c_{30}\int_{(1/\eta_1)^{2+\beta} |y|^{2+2\beta}} ^{\infty}    \frac{|y|^{2+2\beta}}{t^{1+d/2}} dt\le c_{31}|y|^{(1+\beta)(2-d)} \le c_{31} \begin{cases}
	4^{1+\beta}(|y|\wedge 1)^{1+\beta} &\mbox{ if $d=1$},\\
	|x-y|^{2-d} &\mbox{ if $d\ge 2$}
\end{cases} \nn\\
&  \le c_{32} f_0(x,y).
\end{align*}
This completes the proof.  \qed

\begin{lem}\label{l:Green-boundary-1-lower}
(i) Assume $d\ge 2$.	There   exist constants  $c_1>0$ and $C>0$  such that  for all $x,y \in \R^d_0$ with $(|x| \wedge  |y|) \vee |x-y|  \le 2$,
	\begin{align}\label{e:Green-boundary-1-lower}
		G(x,y)  \ge  C\Big( 1\wedge \frac{h(|x|\wedge |y|)}{h(\eta_2|x-y|)}\Big) f_0(x,y) \exp \Big( - \frac{c_1|x-y|}{(|x| \vee |y|)^{1+\beta}}\Big),
	\end{align}
	where $\eta_2>0$ is the constant in Lemma \ref{l:Green-boundary-1-upper}.
	
	\noindent (ii) If $d=1$, then  \eqref{e:Green-boundary-1-lower} holds for all $x,y \in \R^1_0$ with  $xy>0$ and $(|x| \wedge  |y|) \vee |x-y|  \le 2$.
\end{lem}
\pf We assume, without loss of generality, that  $|x|\le |y|$. When $d=1$, we additionally assume that $xy>0$.  We consider the following three cases separately:

\smallskip

\noindent
{\it
Case 1:  $|x| \le \eta_2|x-y|$.} Note that $|x-y|/2\le |y| \le (1+\eta_2)|x-y|$ in this case. Using Theorem \ref{t:smalltime} with $T= (4\eta_2/\eta_1)^{2+\beta}$ in the first inequality below,  the almost monotonicity of $h$ in the second, and  $|x-y|/2\le |y| \le (1+\eta_2)|x-y|$ and  Lemma \ref{l:h-ratio} in the third, we obtain
\begin{align}\label{e:Green-boundary-1-lower-1}
&	G(x,y) \ge  c_1  \int_{(\eta_2|x-y|/\eta_1)^{2+\beta}/2}^{(\eta_2|x-y|/\eta_1)^{2+\beta}}  \Big(1 \wedge \frac{h(|x|)}{h(\eta_1t^{1/(2+\beta)}) }\Big)  \Big(1 \wedge \frac{h(|y|)}{h(\eta_1t^{1/(2+\beta)}) }\Big)  \nn\\
	&\qquad \qquad \qquad \qquad\qquad \qquad\quad \times t^{-d/2} \exp \Big(- \frac{c_2t}{|y|^{2+2\beta}} - \frac{c_3|x-y|^2}{t}\Big)   dt \nn\\
	&\ge  \frac{c_4 h(|x|)}{h(\eta_2|x-y|) } \Big(1 \wedge \frac{h(|x-y|/2)}{h(\eta_2|x-y|) }\Big)    \nn\\
	&\qquad \times  \exp \Big(- \frac{c_2(\eta_2|x-y|/\eta_1)^{2+\beta} }{|y|^{2+2\beta}} - \frac{2c_3}{(\eta_2/\eta_1)^{2+\beta} |x-y|^\beta}\Big) \int_{(\eta_2|x-y|/\eta_1)^{2+\beta}/2}^{(\eta_2|x-y|/\eta_1)^{2+\beta}} \frac{1}{t^{d/2} }    dt  \nn\\
		&\ge  \frac{c_5 |x-y|^{(2+\beta)(2-d)/2}h(|x|)}{h(\eta_2|x-y|) }  \exp \Big( - \frac{ (1+\eta_2)^{1+\beta}( 2^{1+\beta}  - 2^{-1}\eta_2^{-\beta} )_+  \sqrt \kappa \,|x-y|}{\beta |y|^{1+\beta}}  - \frac{c_6|x-y|}{|y|^{1+\beta}} \Big)  .
\end{align}
When $d\ge 2$, since $|x-y|\le 2$ and  $|y| \le (1+\eta_2)|x-y|$, we see
\begin{align*}
	|x-y|^{(2+\beta)(2-d)/2}   \ge (2^{\beta(2-d)/2}\1_{\{d\ge 3\}} + (\Lg (2^\beta(1+\eta_2)^{1+\beta}))^{-1}\1_{\{d=2\}}) f_0(x,y).
\end{align*}
When $d=1$,   using  $|x-y|\le 2$,  \eqref{e:exp-poly} and $|y| \le (1+\eta_2)|x-y|$, we obtain
\begin{align*}
&	|x-y|^{1+\beta/2} \ge \Big( \frac{|x-y|^{1+\beta}}{2^{\beta/2}} \vee  c_7|x-y| \Big) \exp \Big( - \frac{c_6}{(1+\eta_2)^{1+\beta}|x-y|^{\beta}}\Big)\nn\\
	&\ge \Big( \frac{|y|^{1+\beta}}{2^{\beta/2}(1+\eta_2)^{1+\beta}} \vee  c_7|x-y| \Big) \exp \Big( - \frac{c_6|x-y|}{|y|^{1+\beta}}\Big)\ge c_8 f_0(x,y)\exp \Big( - \frac{c_6|x-y|}{|y|^{1+\beta}}\Big).
\end{align*}
Combining these estimates for $f_0(x,y)$ with \eqref{e:Green-boundary-1-lower-1}, we are done in this case.

\smallskip

\noindent
{\it
Case 2:  $|x|>\eta_2|x-y|$ and $|x-y| \ge |y|^{1+\beta}$.} Then $|y|< (1+\eta_2^{-1})|x|$.  Applying Theorem \ref{t:smalltime}  with $T= (4\eta_2/\eta_1^2)^{(2+\beta)/2}$,  we get
\begin{align}\label{e:Green-boundary-1-lower-2}
	G(x,y) &\ge c_9 \int_{2^{-1}\eta_2 |x-y| |x|^{1+\beta}/\eta_1^{2+\beta}}^{\eta_2 |x-y| |x|^{1+\beta}/\eta_1^{2+\beta}}    t^{-d/2}  \exp \Big(- \frac{c_{10}t}{|y|^{2+2\beta}} - \frac{c_{11}|x-y|^2}{t}\Big)  dt \nn\\
	 &\ge   c_{12}(|x-y| |x|^{1+\beta})^{-(d-2 )/2}  \exp \Big(- \frac{c_{10}\eta_2 |x-y| |x|^{1+\beta}}{\eta_1^{2+\beta} |y|^{2+2\beta}} - \frac{2c_{11}\eta_1^{2+\beta}|x-y|}{\eta_2  |x|^{1+\beta}}\Big) \nn\\
	  &\ge      c_{12}(|x-y| |x|^{1+\beta})^{-(d-2 )/2}   \exp \Big(- \frac{c_{13}|x-y|}{ |y|^{1+\beta}}\Big)  .
\end{align}
When $d\ge 2$, since $|x|^{1+\beta} \le |y|^{1+\beta} \le |x-y|$,  we see $(|x-y| |x|^{1+\beta})^{-(d-2 )/2} \ge |x-y|^{2-d} \ge \frac{f_0(x,y)}{\Lg 2}.$ When $d=1$, by \eqref{e:exp-poly}, we get
\begin{align*}
	 (|x-y| |x|^{1+\beta})^{1/2} & \ge c_{14}|x-y| \exp \Big( - \frac{c_{15}|x-y|}{(1+\eta_2^{-1})^{1+\beta}|x|^{1+\beta}}\Big) \ge c_{14}f_0(x,y) \exp \Big( - \frac{c_{15}|x-y|}{|y|^{1+\beta}}\Big) .
\end{align*}
Thus, from \eqref{e:Green-boundary-1-lower-2}, the desired result follows.

\smallskip
\noindent
{\it
Case 3:  $|x|>\eta_2|x-y|$ and $|x-y| < |y|^{1+\beta}$.} Applying Theorem \ref{t:smalltime}  with $T=  2^{2+2\beta}/(2\eta_1^{2+\beta})$ and using $|y|< (1+\eta_2^{-1})|x|$,  we get
\begin{align*}
	G(x,y) &\ge c_{16} \int_{ |x-y|^{2}/(4(1+\eta_2^{-1})^{2+2\beta}\eta_1^{2+\beta})}^{ |x|^{2+2\beta}/(2\eta_1^{2+\beta})}    t^{-d/2}   dt
	\ge c_{17}f_0(x,y).
\end{align*} 
The proof is complete. \qed

\begin{lem}\label{l:Green-boundary-2}
(i) If $d\ge 2$, then	there   exist comparison constants  such that 
	\begin{align*}
		G(x,y)     \asymp h(|x|) |x-y|^{2-d} \quad \text{ for all $x,y \in \R^d_0$ with $|x| \wedge  |y|  \le  2\le |x-y|$} .
	\end{align*}

	\noindent (ii) If $d=1$, then there   exist comparison constants  such that 
	\begin{align*}
		G(x,y)     \asymp h(|x|) \quad \text{  for all $x,y \in \R^1_0$ with $xy>0$ and  $|x| \wedge  |y|  \le  2\le |x-y|$}. 
	\end{align*}
\end{lem}
\pf Let $x,y \in \R^d_0$ with $|x| \wedge  |y|  \le  2\le |x-y|$. Without loss of generality, we assume $|x|\le |y|$.  When $d=1$, we additionally assume that $xy>0$. Note that $|y| \ge |x-y|/2 \ge 1$ and $|y| \le |x-y| + |x| \le 2|x-y|$.   For the lower bound, using Theorem \ref{t:largetime}, 
 \eqref{e:Log-prop-0} and 
 \eqref{e:H-arbitrary-R} (if $d=1$ or $d=2$), we see
\begin{align*}
	G(x,y) &\ge \int_{4|x-y|^2}^{5|x-y|^2}  p(t,x,y)dt \ge  c_1h(|x|) \int_{4|x-y|^2}^{5|x-y|^2} \Big( \frac{\1_{\{d\ge 3\}}}{t^{d/2}}  + \frac{\1_{\{d=2\}} \Lg |y|}{t (\log t)^2} + \frac{\1_{\{d=1\}}|y|}{t^{3/2}} \Big) dt \\
	&\ge   c_2 h(|x|) \Big( \1_{\{d\ge 3\}}|x-y|^{2-d} + \frac{\1_{\{d=2\}} \Lg |y|}{\log |x-y|} + \frac{\1_{\{d=1\}}|y|}{|x-y|} \Big) \\
	&\ge   c_3 h(|x|) (\1_{\{d\ge 2\}}|x-y|^{2-d}   + \1_{\{d= 1\}}) .
\end{align*}
For the upper bound, by  Theorems \ref{t:smalltime} and \ref{t:largetime},  \eqref{e:Log-prop-0} and
\eqref{e:H-arbitrary-R} (if $d=1$ or $d=2$), we obtain
\begin{align*}
\frac{	G(x,y)}{h(|x|)} & \le c_4 \int_0^{4}  \frac{1}{t^{d/2}h(\eta_1 t^{1/(2+\beta)})}\exp \Big(- \frac{c_5|x-y|^2}{t}\Big)  dt \\
&\quad + c_4\int_{4}^{|x-y|^2} \Big( \frac{\1_{\{d\ge 3\}}}{t^{d/2}}  + \frac{\1_{\{d=2\}}}{t (\log t)} + \frac{\1_{\{d=1\}}}{t} \Big)  \exp \Big(- \frac{c_5|x-y|^2}{t}\Big)  dt \\
&\quad  + c_4  \int_{|x-y|^2}^\infty \Big( \frac{\1_{\{d\ge 3\}}}{t^{d/2}}  + \frac{\1_{\{d=2\}} \Lg |y|}{t (\log t)^2} + \frac{\1_{\{d=1\}}|y|}{t^{3/2}} \Big) dt =:c_4(I_1+I_2+I_3).
\end{align*}
 For $I_1$, using \eqref{e:exp-poly} and Lemma \ref{l:h-ratio} in the first inequality below,  Lemma \ref{l:Green-elementary} in the third and \eqref{e:exp-poly} in the fifth,  we get
\begin{align*}
	I_1&\le c_6	 \int_0^{4}  \exp \Big( \frac{2\sqrt \kappa}{\beta (\eta_1t^{1/(2+\beta)})^\beta} - \frac{c_5|x-y|^2}{t}  \Big)  dt\\
	&\le c_6	\Big[ \sup_{0<t\le 4} \exp \Big( - \frac{c_5|x-y|^2}{2t}\Big) \Big] \Big[ \sup_{t>0}  \exp \Big( \frac{2\sqrt \kappa}{\beta (\eta_1t^{1/(2+\beta)})^\beta} - \frac{c_5|x-y|^2}{2t}\Big)  \Big]\int_0^{4}   dt \\
	&\le 4c_6	   \exp \Big( - \frac{c_5|x-y|^2}{8} + \frac{c_7}{|x-y|^\beta}\Big) \le c_6 e^{c_7/2^\beta} |x-y|^2	   \exp \Big(- \frac{c_5|x-y|^2}{8}\Big) \\
	& \le c_8 (\1_{\{d\ge 2\}}|x-y|^{2-d}   + \1_{\{d= 1\}}).
\end{align*}
For $I_2$, by \eqref{e:Green-elementary}, we have
\begin{align*}
	I_2 \le\int_{4}^{|x-y|^2} \Big( \frac{\1_{\{d\ge 2\}}}{t^{d/2}}   + \frac{\1_{\{d= 1\}}}{t} \Big)  \exp \Big(- \frac{c_5|x-y|^2}{t}\Big)  dt \le c_9(\1_{\{d\ge 2\}}|x-y|^{2-d}   + \1_{\{d= 1\}}).
\end{align*}
For $I_3$, using  $|x-y|/2\le |y| \le 2|x-y|$, we obtain
\begin{align*}
	I_3 \le c_{10}  \Big( \1_{\{d\ge 3\}} |x-y|^{d-2}  + \frac{\1_{\{d=2\}} \Lg |y|}{\log |x-y|} + \frac{\1_{\{d=1\}}|y|}{|x-y|} \Big) \le c_{11}(\1_{\{d\ge 2\}}|x-y|^{2-d}   + \1_{\{d= 1\}}).
\end{align*}
The proof is complete. \qed

\noindent \textbf{Proof of Theorem \ref{t:Green}.}  Let $x,y \in \R^d_0$. Without loss of generality, we assume $|x| \le |y|$. When $d=1$, we additionally assume that $xy>0$. If $|x| \ge 2$, then  
the result follows from Lemmas \ref{l:Green-exterior} and \ref{l:Green-boundary-0}. If $|x| \vee |x-y|\le 2$, then the result follows from Lemmas \ref{l:Green-boundary-1-upper} and \ref{l:Green-boundary-1-lower}. Assume $|x| \le 2 \le |x-y|$. Then $ |y| \ge |x-y|/2  \ge 1$ and $|x-y| \le 2^{-\beta}|x-y|^{1+\beta} \le  2 |y|^{1+\beta}$. Thus, $\exp(-|x-y|/|y|^{1+\beta}) \asymp 1$ and by the almost monotonicity of $h$,  we obtain
\begin{align*}
	\Big( 1\wedge \frac{h(|x|)}{h(\eta_2|x-y|)}\Big) f_0(x,y)\asymp h(|x|) (\1_{\{d\ge 2\}}|x-y|^{2-d}  + \1_{\{d=1\}} ).
\end{align*}
Now the result follows from Lemma \ref{l:Green-boundary-2}. \qed

 \bigskip
\noindent
{\bf Acknowledgements:}
Part of the research for this paper was done while the third-named author was visiting Jiangsu Normal University, where he was partially supported by a grant from the National Natural Science Foundation of China (11931004, Yingchao Xie).

	\vskip 0.2truein

\noindent {\bf Soobin Cho:} Department of Mathematics,
University of Illinois Urbana-Champaign,
Urbana, IL 61801, U.S.A.
Email: \texttt{soobinc@illinois.edu}
	
	\medskip
	
\noindent {\bf Panki Kim:}
Department of Mathematical Sciences and Research Institute of Mathematics,
	Seoul National University,	Seoul 08826, Republic of Korea.
	Email: \texttt{pkim@snu.ac.kr}

	\medskip

\noindent {\bf Renming Song:} Department of Mathematics,
University of Illinois Urbana-Champaign,
Urbana, IL 61801, U.S.A.
Email: \texttt{rsong@illinois.edu}

\end{document}